        \documentclass[a4paper, 11pt, twoside,leqno]{article}

        \usepackage{bm}
        \usepackage{geometry}
        \usepackage{amssymb,amsmath}

        \usepackage[utf8]{inputenc}
        \usepackage[english]{babel}
        \usepackage{amsthm,thmtools}
        \usepackage{tikz}
            \usetikzlibrary{cd}
            \usetikzlibrary{arrows.meta, positioning, quotes}

        \usepackage{mathtools}
        \usepackage{titlesec}
        \usepackage{url}
        
        \usepackage{afterpage}
        \usepackage{enumitem}
        \setlist{nolistsep}
        \usepackage[fulladjust]{marginnote}
        \usepackage{stmaryrd} 
        \usepackage[textwidth=0.89in,textsize=scriptsize]{todonotes}

        \usepackage{babel,csquotes,xpatch}
        \usepackage{cite}

        \usepackage{graphicx}
        \usepackage{array}
        \usepackage{quiver}

        \usepackage{hyperref} 
        \hypersetup{
        colorlinks=true,
        linkcolor=blue,
        filecolor=magenta,      
        urlcolor=cyan,
        }

        \usepackage{cleveref} 

        \usepackage{lmodern}
        \usepackage[mathcal]{eucal}

        \newcommand{\cS}{{\mathcal S}}
        \newcommand{\bGamma}{\mathbf \Gamma} 

        \tikzset{
            symbol/.style={
            draw=none,
            every to/.append style={
            edge node={node [sloped, allow upside down, auto=false]{$#1$}}}
            }
        }

        \usetikzlibrary{decorations.markings}

            \makeatletter
            \tikzcdset{
              open/.code     = {\tikzcdset{hook, circled};},
              closed/.code   = {\tikzcdset{hook, slashed};},
              open'/.code    = {\tikzcdset{hook', circled};},
              closed'/.code  = {\tikzcdset{hook', slashed};},
              circled/.code  = {\tikzcdset{markwith = {\draw (0,0) circle (.375ex);}};},
              slashed/.code  = {\tikzcdset{markwith = {\draw[-] (-.4ex,-.4ex) -- (.4ex,.4ex);}};},
              markwith/.code ={
                \pgfutil@ifundefined%
                {tikz@library@decorations.markings@loaded}%
                {\pgfutil@packageerror{tikz-cd}{You need to say %
                  \string\usetikzlibrary{decorations.markings} to use arrows with markings}{}}{}%
                \pgfkeysalso{/tikz/postaction = {
                  /tikz/decorate,
                  /tikz/decoration={markings, mark = at position 0.5 with {#1}}}
                }
              },
            }
            \makeatother
    
        \renewcommand{\ss}{{\text{ss}}}

        \newcommand{\N}{{\ensuremath{{\mathbb N}}}}
        \newcommand{\C}{{{\mathbb C}}}
        
        \newcommand{\R}{{\ensuremath{{\mathbb R}}}}

        \newcommand{\Z}{{\ensuremath{{\mathbb Z}}}}
        \newcommand{\Q}{{\ensuremath{{\mathbb Q}}}}
        \newcommand{\G}{{\ensuremath{{\mathbb G}}}}
        \newcommand{\A}{{\ensuremath{{\mathbb A}}}}
        \let\P\relax
        \newcommand{\P}{{\ensuremath{{\mathbb P}}}}

        \newcommand{\bA}{{\mathbb A}}

        \newcommand{\cA}{{\mathcal A}}
        \newcommand{\cB}{{\mathcal B}}
        \newcommand{\cC}{{\mathcal C}}
        
        \newcommand{\cE}{{\mathcal E}}
        \newcommand{\cF}{{\mathcal F}}

        \newcommand{\cI}{{\mathcal I}}

        \newcommand{\cL}{{\mathcal L}}
        \newcommand{\cM}{{\mathcal M}}
        \newcommand{\cN}{{\mathcal N}}
        \newcommand{\cO}{{\mathcal O}}

        \newcommand{\cR}{{\mathcal R}}

        \newcommand{\cU}{{\mathcal U}}
        \newcommand{\cV}{{\mathcal V}}
        \newcommand{\cW}{{\mathcal W}}
        \newcommand{\cX}{{\mathcal X}}
        \newcommand{\cY}{{\mathcal Y}}
        \newcommand{\cZ}{{\mathcal Z}}

        \newcommand{\fl}{{\mathfrak l}}
        
        \newcommand{\fk}{{\mathfrak k}}
        \newcommand{\fg}{{\mathfrak g}}

        \newcommand{\bXi}{\mathbf{\Xi}}
        
        \newcommand{\abs}[1]{\left|  #1\right| } 

        \makeatletter
        \newcommand{\colim@}[2]{%
          \vtop{\m@th\ialign{##\cr
            \hfil$#1\operator@font colim$\hfil\cr
            \noalign{\nointerlineskip\kern1.5\ex@}#2\cr
            \noalign{\nointerlineskip\kern-\ex@}\cr}}%
        }
        \newcommand{\colim}{%
          \mathop{\mathpalette\colim@{\rightarrowfill@\scriptscriptstyle}}\nmlimits@
        }
        \renewcommand{\varprojlim}{%
          \mathop{\mathpalette\varlim@{\leftarrowfill@\scriptscriptstyle}}\nmlimits@
        }
        \renewcommand{\varinjlim}{%
          \mathop{\mathpalette\varlim@{\rightarrowfill@\scriptscriptstyle}}\nmlimits@
        }
        \makeatother

        \newcommand{\id}{\text{id}}
        
        \newcommand{\fppf}{\text{fppf}}

        \newcommand{\git}{\mathbin{
          \mathchoice{/\mkern-6mu/}
            {/\mkern-6mu/}
            {/\mkern-5mu/}
            {/\mkern-5mu/}}}
        \newcommand{\lex}{{\text{lex}}}

        \newcommand{\cat}[1]{\mathbf{\operatorname{#1}}}
        \newcommand{\pre}{{\text{pre}}}

        \DeclareMathOperator{\im}{im}

        \DeclareMathOperator{\ev}{ev}
        
        \DeclareMathOperator{\Hom}{Hom}
        \DeclareMathOperator{\uHom}{\underline{Hom}}

        \DeclareMathOperator{\Aut}{\text{Aut}}
        \DeclareMathOperator{\Sym}{\text{Sym}}
        \DeclareMathOperator{\GL}{\text{GL}}
        \DeclareMathOperator{\Lie}{\text{Lie}}

        \DeclareMathOperator{\Stab}{\text{Stab}}
        \DeclareMathOperator{\Spec}{Spec}
        
        \DeclareMathOperator{\Proj}{Proj}

        \DeclareMathOperator{\st}{\vert\ }

        \DeclareMathOperator{\Filt}{Filt}
        \DeclareMathOperator{\Grad}{Grad}
        \let\ev\relax
        \DeclareMathOperator{\ev}{ev}
        \DeclareMathOperator{\gr}{gr}
        \DeclareMathOperator{\conv}{conv}

        \newcommand{\norm}[1]{\lVert #1\rVert}

        \DeclareMathOperator{\cone}{cone}
        \let\Re\relax
        \DeclareMathOperator{\Re}{\mathfrak{Re}}
        \let\Im\relax
        \DeclareMathOperator{\Im}{\mathfrak{Im}}
        \DeclareMathOperator{\Bl}{Bl} 
        
        \DeclareMathOperator{\Pic}{Pic}
        
        \DeclareMathOperator{\Flag}{Flag}
        
        \DeclareMathOperator{\rk}{rk}

        \newcommand{\dash}{\,\text{-}}
        \newcommand{\red}{\text{red}}

        \DeclareMathOperator{\ind}{ind} 
        \DeclareMathOperator{\Bun}{Bun} 

        \DeclareMathOperator{\wt}{\text{wt}}

        \DeclareMathOperator{\qinfilt}{\Q^\infty\text{-}\Filt}
        \DeclareMathOperator{\qfilt}{\Q\,\text{-}\Filt}
        
        \renewcommand{\phi}{\varphi}

        \makeatletter 
        \DeclareFontEncoding{LS1}{}{}
        \DeclareFontSubstitution{LS1}{stix}{m}{n}
        \DeclareMathAlphabet{\mathcaldos}{LS1}{stixscr}{m}{n}
        \makeatother

        \DeclareMathOperator{\cRep}{\mathcaldos{R\mkern-4mu e\mkern-4.5mu p\mkern-1mu}}

        \DeclareMathOperator{\Db}{D^{\text b}} 
        \DeclareMathOperator{\Ind}{Ind} 

        \declaretheorem[numberwithin=subsection, name=Theorem]{theorem}
        \declaretheorem[sibling=theorem, style=definition, name=Definition]{definition}
        \declaretheorem[sibling=theorem, name=Lemma]{lemma}
        \declaretheorem[sibling=theorem, name=Corollary]{corollary}
        \declaretheorem[sibling=theorem, name=Proposition]{proposition}
        \declaretheorem[sibling=theorem, name=Conjecture]{conjecture}
        \declaretheorem[sibling=theorem, style=remark, name=Remark]{remark}
        \declaretheorem[sibling=theorem, style=remark, name=Construction]{construction}
        
        \declaretheorem[sibling=theorem, style=remark, name=Claim]{claim}
        \declaretheorem[sibling=theorem, style=remark, name=Example]{example}

        \declaretheorem[sibling=theorem, style=remark, name=Assumption]{assumption}
        
        \declaretheorem[sibling=theorem, style=definition, name=Proposition and Definition]{propositionanddefinition}

\tikzset{
    labl/.style={anchor=south, rotate=90, inner sep=.5mm}
}
\tikzset{
    lablh/.style={anchor=south, inner sep=.5mm}
}

\title{\textbf{Refined Harder-Narasimhan filtrations in \\ moduli theory}}
\author{Andrés Ibáñez Núñez}
\date{November 2023}

\begin{document}
\maketitle

\begin{abstract}
We define canonical refinements of Harder-Narasimhan filtrations and stratifications in moduli theory, generalising and relating work of Haiden-Katzarkov-Kontsevich-Pandit and Kirwan. More precisely, we define a canonical stratification on any noetherian algebraic stack $\cX$ with affine diagonal that admits a good moduli space and is endowed with a norm on graded points. The strata live in a newly defined \emph{stack of sequential filtrations} of $\cX$. Therefore the stratification gives a canonical sequential filtration, the \emph{iterated balanced filtration}, for each point of $\cX$. Examples of applicability include moduli of principal bundles on a curve, moduli of objects at the heart of a Bridgeland stability condition and moduli of K-semistable Fano varieties. We conjecture that the iterated balanced filtration describes the asymptotics of the Kempf-Ness flow in Geometric Invariant Theory, as part of a larger project aiming to describe the asymptotics of natural flows in moduli theory. In the case of quotient stacks by diagonalisable algebraic groups, we give an explicit description of the iterated balanced filtration in terms of convex geometry.

\end{abstract}
\tableofcontents

\section{Introduction}
There is a fascinating relation in moduli theory between three seemingly quite different phenomena: stratifications, filtrations and flows. To illustrate this, consider the moduli stack $\Bun(C)$ of vector bundles on a smooth projective curve $C$ over the complex numbers.
Every vector bundle $E\in \Bun(C)(\C)$ has a canonical filtration, its Harder-Narasimhan filtration $\lambda_{\text{HN}}(E)$ \cite{_Harder_Onthecohomologygroupsofmodulispacesofvectorbundlesoncurves}, and the assignment of $\lambda_{\text{HN}}(E)$ for every $E$ defines a stratification of $\Bun(C)$ by Harder-Narasimhan type. 
This stratification was used by Atiyah and Bott to study the cohomology of moduli spaces of semistable vector bundles \cite{_Atiyah_TheYangMillsEquationsoverRiemannSurfaces}. 
On the other hand, one can consider the Yang-Mills flow on the space $M_E$ of hermitian metrics on the vector bundle $E$. The Harder-Narasimhan filtration of $E$ can then be realised as the slope at infinity of this flow from any initial metric.

However, the Harder-Narasimhan filtration is insufficient to detect more subtle asymptotic properties of the Yang-Mills flow. In recent work, Haiden-Katzarkov-Kontsevich-Pandit gave a complete description of these asymptotics in terms of a refinement of the Harder-Narasimhan filtration that we call in this work the \emph{HKKP filtration} \cite{_Haiden_Semistabilitymodularlatticesanditeratedlogarithms,_Haiden_Iteratedlogarithmsandgradientflows}. Its construction relies solely on the lattice of subbundles of $E$, and in fact the HKKP filtration is available more generally for moduli of objects in a linear category, where a lattice of subobjects of a given object can be considered. There are however other moduli problems where there is no such lattice, and the HKKP filtration is not defined in these cases. Examples include principal bundles, K-semistable Fano varieties, and GIT quotients. There are analogues of the Yang-Mills flow in these moduli problems whose asymptotics are not yet well understood.

Coming back to the case of $\Bun(C)$, and following the parallelism between flows, filtrations and stratifications, one would expect that the assignment of the HKKP filtration for every $E$ should define a stratification of $\Bun(C)$ by type of HKKP filtration refining the stratification by Harder-Narasimhan type. Refinements of this kind had previously been considered by Kirwan for GIT quotients of smooth projective varieties \cite{_Kirwan_RefinementsoftheMorsestratificationofthenormsquareofthemomentmap} and for $\Bun(C)$ \cite{_Kirwan_ModulispacesofbundlesoverRiemannsurfacesandtheYangMillsstratificationrevisited}, but the relation to the HKKP filtration has remained unclear. A caveat of Kirwan's construction in the GIT setting is that it does not give an analogue of the HKKP filtration that could be used to describe asymptotics of flows. On the other hand, many moduli problems of interest, like Bridgeland semistable objects or K-semistable Fano varieties, cannot be realised as GIT quotients.

In this work, we propose a definition of refined Harder-Narasimhan filtrations in abstract moduli theory, valid beyond the case of moduli objects in a linear category.
Our construction is linked to stratifications from the outset, and in particular it provides a precise link between HKKP filtrations and Kirwan's refined stratification in the case of $\Bun(C)$. We work with general algebraic stacks admitting good moduli spaces, and in this way we cover a wide range of examples, like GIT quotients over general bases, moduli stacks of principal bundles on a curve, moduli of K-semistable Fano varieties, and moduli of objects at the heart of a Bridgeland stability condition (\Cref{section: examples}).

The definition of our filtration involves as main ingredients blow-ups of stacks and $\Theta$-stratifications \cite{_HalpernLeistner_Onthestructureofinstabilityinmodulitheory}, as well as the construction of a \emph{stack of sequential filtrations} $\Filt_{\Q^\infty}(\cX)$ for an algebraic stack $\cX$ (\Cref{definition: Q^infty filtrations and gradings}). When $\cX$ is noetherian, has affine diagonal, admits a good moduli space, and is endowed with a norm on graded points, we define a canonical \emph{balancing stratification} of $\cX$ (\Cref{theorem: sequential stratifications for good moduli stacks} and \Cref{definition: balancing stratification}). It is a \emph{sequential stratification} of $\cX$ in the sense that the strata are also substacks of $\Filt_{\Q^\infty}(\cX)$ (\Cref{definition: sequential stratification}). The balancing stratification assigns to every point $x$ of $\cX$ a canonical sequential filtration of $x$ that we call the \emph{iterated balanced filtration}.

Despite the apparently complicated nature of its definition, the iterated balanced filtration is often computable in terms of combinatorial data. In this paper, we establish a convex geometric algorithm to compute the filtration in the case of quotient stacks by the action of a diagonalisable algebraic group (\Cref{corollary: iterated balanced filtration for states equals that for stacks}). The proof involves the theory of chains of stacks, developed in \Cref{section: chains of stacks}. In a sequel to this paper, we will use these techniques to show that the iterated balanced filtration agrees with the HKKP filtration when the latter makes sense (\Cref{theorem: correspondence HKKP filtration and iterated balanced filtration}).

Our motivation for this work is the expectation that the iterated balanced filtration describes the asymptotics of natural flows in moduli theory, like the Yang-Mills flow for principal bundles, the Calabi flow for varieties and the gradient flow of the Kempf-Ness potential in Geometric Invariant Theory.
We provide a conjectural statement for this expectation in the case of GIT for affine spaces (\Cref{conjecture: on asymptotics of gradient flows affine GIT}).

\subsection{Normed good moduli stacks}
The stratification of $\Bun(C)$ by Harder-Narasimhan type is a $\Theta$-stratification. In particular, every stratum $\cS$ has an $\A^1$-retraction $\cS\to \cZ$ onto what is called its centre $\cZ$. In this case, $\cZ$ enjoys the property of admitting a good moduli space $\cZ\to Z$, that is, a map to an algebraic space $Z$ that best approximates the stack $\cZ$ (we recall the precise definition from \cite{_Alper_GoodmodulispacesforArtinstacks} in \ref{definition: good moduli space}). Having a good moduli space is a strong condition on the stack that implies many desirable properties.

We may hope that algebraic stacks admitting good moduli spaces have canonical stratifications and, in a suitable sense, canonical filtrations for every point. These could then be pulled back along the retraction $\cS\to \cZ$ from each centre $\cZ$ and produce, in the case of $\Bun(C)$, the sought-after stratification by HKKP type, as well as recovering the HKKP filtration of every point. This is close to being true. What is missing is stability type data on $\cZ$ on which this stratification depends. The correct notion for our purposes turns out to be a \emph{norm on graded points} of $\cZ$, a concept from the Beyond GIT programme \cite{_HalpernLeistner_Onthestructureofinstabilityinmodulitheory}.

A \emph{graded point} of a stack $\cX$ is a map $g\colon B\G_{m,k}\to \cX$, where $B\G_{m,k}$ is the classifying stack of the multiplicative group over a field $k$. Graded points on $\cX$ can also be seen as ordinary points on the mapping stack $\Grad(\cX)=\uHom(B\G_m,\cX)$. A norm on graded points of $\cX$ is roughly speaking the data of a positive real number $\norm{g}$ for every graded point $g\colon B\G_{m,k}\to \cX$, and this data is required to satisfy some local constancy and nondegeneracy conditions (\Cref{definition: norm on graded points}).
For $\Bun(C)$, a natural norm is given by the rank of vector bundles. In this case, a graded point $g\colon B\G_{m,\C}\to \Bun(C)$ corresponds to a vector bundle $E$ together with a direct sum decomposition $E=\bigoplus_{c\in \Z} E_c$, and its norm is then defined by the formula $\norm{g}^2=\sum_{c\in \Z}c^2 \rk(E_c)$. This restricts to a norm on graded points on each of the centres $\cZ$ of the Harder-Narasimhan stratification.

Many other moduli stacks $\cX$ are also naturally endowed with a $\Theta$-stratification, related to the assignment, for every point $x$ of $\cX$, of a Harder-Narasimhan filtration of $x$ in a generalised sense. It is often the case that each centre $\cZ$ of the stratification has a good moduli space and is naturally endowed with a norm on graded points (see \Cref{section: examples} for examples). Therefore, in order to define refined Harder-Narasimhan filtrations in moduli theory in great generality, it is enough to deal with algebraic stacks admitting a good moduli space and endowed with a norm on graded points. We
call these \emph{normed good moduli stacks} for simplicity. Our aim thus becomes to stratify and produce canonical filtrations for normed good moduli stacks.

\subsection{Sequential filtrations and stratifications}

The first obstacle we encounter is the very meaning of filtration in this generality. The HKKP filtration of a vector bundle $E$ consists of both a chain
\begin{equation*}
0\neq E_0\subset E_1 \subset \cdots \subset E_n=E
\end{equation*}
of subbundles and a chain
\begin{equation*}
c_0>c_1>\cdots>c_n\end{equation*}
of labels $c_i\in \Q^\infty$. Here, $\Q^\infty$ is the set of eventually zero sequences of rational numbers, ordered lexicographically. In this sense the HKKP filtration is a \emph{sequential filtration}, or \emph{$\Q^\infty$-filtration}, of $E$. 

We first look at the well-studied case of $\Z$-filtrations, that is, when the labels $c_i$ are integers. These are closely related to the quotient stack $\A^1/\G_m$ of the affine line $\A^1$ by the scaling action of the multiplicative group, also denoted $\Theta=\A^1/\G_m$. Indeed, the mapping stack $\Filt(\Bun(C))=\uHom(\A^1/\G_m,\Bun(C))$ parametrises vector bundles endowed with a $\Z$-filtration \cite{_Alper_Existenceofmodulispacesforalgebraicstacks,_Heinloth_HilbertMumfordstabilityonalgebraicstacksandapplicationstoGbundlesoncurves}. Therefore it makes sense to define a $\Z$-filtration of a $k$-point $x$ in a general stack $\cX$ to be a map $f\colon \Theta_k\to \cX$ together with an isomorphism $x\sim f(1)$, and to define the stack of filtrations on $\cX$ to be the mapping stack $\Filt(\cX)=\uHom(\Theta,\cX)$. This idea lies at the heart of Halpern-Leistner's approach to moduli theory beyond GIT \cite{_HalpernLeistner_Onthestructureofinstabilityinmodulitheory}. There is an \emph{associated graded map $\gr\colon \Filt(\cX)\to \Grad(\cX)$} and a forgetful map $\ev_1\colon \Filt(\cX)\to \cX$. 
For our purposes, it is better to consider $\Q$-filtrations, and we work instead with the stacks of $\Q$-filtrations $\Filt_\Q(\cX)$ and of $\Q$-gradings $\Grad_\Q(\cX)$, that we construct by formally localising with respect to the natural action of the monoid $(\Z_{>0},\cdot,1)$ on $\Filt(\cX)$ and $\Grad(\cX)$ (\Cref{definition: rational filtrations and gradings}). Rational filtrations of a point $x\in \cX(k)$ are defined in the same way and they form a set $\qfilt(\cX,x)$.

Our first goal is to construct an algebraic stack $\Filt_{\Q^\infty}(\cX)$ that, in the case $\cX=\Bun(C)$, parametrises vector bundles $E$ endowed with a sequential filtration.
To this aim, we observe that giving a $\Q^\infty$-filtration on a vector bundle $E$ on $C$ is equivalent to first giving a $\Q$-filtration $F_\bullet$ of $E$, then a $\Q$-filtration of the associated graded object $\gr F_\bullet$, and so on, until the process finishes in finitely many steps. For general stacks, we formalise this idea by defining the stack $\Filt_{\Q^n_{\text{lex}}}(\cX)$ of $\Q^n_{\text{lex}}$-filtrations inductively as a fibre product \[\Filt_{\Q^n_{\text{lex}}}(\cX)=\Filt_{\Q^{n-1}_{\text{lex}}}(\cX)\times_{\gr,\Grad_{\Q^{n-1}}(\cX),\ev_1} \Filt_\Q\left(\Grad_{\Q^{n-1}}(\cX)\right).\] Here, $\Q^n_\lex$ is $\Q^n$ endowed with the lexicographic order, and $\Grad_{\Q^l}(\cX)$ is defined inductively as $\Grad_{\Q^l}(\cX)=\Grad_\Q\left(\Grad_\Q^{l-1}(\cX)\right)$. The associated graded map $\gr\colon \Filt_{\Q^l_{\lex}}(\cX)\to \Grad_{\Q^l}(\cX)$ is also defined by induction. Then we set:

\begin{definition}[\Cref{definition: Q^infty filtrations and gradings}]
The stack $\Filt_{\Q^\infty}(\cX)$ of \emph{sequential filtrations}, or \emph{$\Q^\infty$-filtrations}, of $\cX$ is the colimit of the stacks $\Filt_{\Q^n_{\text{lex}}}(\cX)$ when $n$ tends to $\infty$.
\end{definition}

Under reasonable conditions on $\cX$, the stack of $\Q^\infty$-filtrations $\Filt_{\Q^\infty}(\cX)$ is algebraic (\Cref{proposition: algebraicity of Q^infty filtrations}). There is also a map $\ev_1\colon \Filt_{\Q^\infty}(\cX)\to \cX$ corresponding to “forgetting the filtration”, so it now makes sense to define a \emph{$\Q^\infty$-filtration} of a field-valued point $x\colon \Spec k\to \cX$ to be a $k$-point $\lambda\in \Filt_{\Q^\infty}(\cX)(k)$ together with an isomorphism $x\sim \ev_1(\lambda)$. We denote $\qinfilt(\cX,x)$ the set of $\Q^\infty$-filtrations of $x$ (\Cref{definition:Q^infty filtrations of a point}).

The relation between stratifications of $\cX$ and sequential filtrations is encapsulated in the definition of \emph{sequential stratification}.

\begin{definition}[\Cref{definition: sequential stratification}]
A \emph{sequential stratification} of an algebraic stack $\cX$ is a family $(\cS_\alpha)_{\alpha\in \Gamma}$ of locally closed substacks of $\Filt_{\Q^\infty}(\cX)$, indexed by a partially ordered set $\Gamma$, such that
\begin{enumerate}
\item each composition $\cS_c\to \Filt_{\Q^\infty}(\cX)\to \cX$ is a locally closed immersion,
\item the topological spaces $\abs{\cS_c}$ are pairwise disjoint and cover $\abs{\cX}$, and
\item for every $c\in \Gamma$, the union $\bigcup_{c'\leq c} \abs{\cS_{c'}}$ is open in $\cX$.
\end{enumerate}
\end{definition}
Thus the strata $\cS_c$ are locally closed substacks of $\cX$ together with a choice of lift to $\Filt_{\Q^\infty}(\cX)$. Therefore a sequential stratification provides each point $x$ of $\cX$ with a choice of sequential filtration of $x$.

The definition of sequential stratification is inspired in Halpern-Leistner's definition of $\Theta$-stratification of a stack $\cX$ \cite{_HalpernLeistner_Onthestructureofinstabilityinmodulitheory} (recalled in \Cref{definition: Theta-stratification my version}), which roughly speaking is a partition of $\cX$ into locally closed substacks $\cS_c$ of $\cX$ that are also open substacks of $\Filt_\Q(\cX)$. Very importantly, each stratum $\cS_c$ retracts onto what is called its \emph{centre} $\cZ_c$, which is an open substack of $\Grad_\Q(\cX)$.
Establishing existence and nice properties of certain $\Theta$-stratifications is a crucial ingredient of our construction of the balancing stratification.

If $\cX$ is a normed good moduli stack and $f\colon \cY\to \cX$ is a representable polarised projective morphism it is a result of Halpern-Leistner \cite[Theorem 5.6.1]{_HalpernLeistner_Onthestructureofinstabilityinmodulitheory}
that $\cY$ has a natural $\Theta$-stratification induced by the norm on $\cX$ and the polarisation. This generalises Kirwan's construction of the instability stratification in GIT \cite{_Kirwan_Cohomologyofquotientsinsymplecticandalgebraicgeometry}. Our first result establishes a property of this stratification that will be fundamental for our purposes:
\begin{theorem}[\Cref{theorem: theta stratification proper over gms}]\label{theorem: theta-stratifications introduction}
Let $\cX$ be a normed noetherian good moduli stack with affine diagonal. Consider a representable projective morphism $f\colon \cY\to \cX$ and the $\Theta$-stratification $(\cS_c)_{c\in \Q_{\geq 0}}$ of $\cY$ induced by an $f$-ample line bundle $\cL$ and the norm on graded points. Then the centre $\cZ_c$ of every stratum $\cS_c$ has a good moduli space.
\end{theorem}
For the proof, we need to consider the more general case when $f$ is proper and the line bundle $\cL$ is replaced by the weaker notion of \emph{linear form on graded points}, and then use the concepts of $\Theta$-monotonicity and $S$-monotonicity developed in \cite{_HalpernLeistner_Onthestructureofinstabilityinmodulitheory} to check $\Theta$-reductivity and $S$-completeness of $\cZ_c$, which implies the existence of a good moduli space \cite{_Alper_Existenceofmodulispacesforalgebraicstacks}.

\subsection{The balancing stratification}

Our main construction (\Cref{theorem: sequential stratifications for good moduli stacks} and \Cref{definition: balancing stratification}) produces a canonical sequential stratification for every normed noetherian good moduli stack $\cX$ with affine diagonal. We call it the \emph{balancing stratification}, since the adjective \emph{balanced} was used both by Kirwan and Haiden-Katzarkov-Kontsevich-Pandit to describe the filtrations they studied. The balancing stratification produces a canonical sequential stratification, the \emph{iterated balanced filtration}, for every point of the stack $\cX$.

Recall that a point $p$ in a stack $\cX$ with good moduli space $\pi\colon \cX\to X$ is said to be \emph{polystable} if it is closed in the fibre of  $\pi$ containing $p$. Every fibre of $\pi$ contains a unique polystable point. If $x\in \cX(k)$ is a non-polystable geometric point, it follows from a result of Kempf \cite{_Kempf_InstabilityinInvariantTheory}
that there is a filtration $\lambda\colon \Theta_k \to \cX$ of $x$ such that $\lambda(0)=y$ is polystable. The question arises whether one can choose a canonical such polystable degeneration $\lambda$. For this problem, it is useful to endow $\cX$ with a norm on graded points. After replacing $\cX$ with a $\pi$-saturated open substack containing $x$, we may assume that $y$ lies in the closed substack $\cX^{\max}$ of points with maximal stabiliser dimension, defined by Edidin-Rydh \cite{_Edidin_CanonicalreductionofstabilizersforArtinstackswithgoodmodulispaces}. Kempf defines a natural number $\langle \lambda, \cX^{\max}\rangle\in \N$ (\Cref{definition: Kempfs intersection number}) that can be thought of as measuring the velocity at which $\lambda$ converges to $y$. The definition extends to rational filtrations $\lambda$, giving a rational number $\langle \lambda, \cX^{\max}\rangle$. It turns out that there is a unique rational filtration $\lambda_{\text{b}}(x)\in \qfilt(\cX,x)$ such that $\langle \lambda, \cX^{\max}\rangle\geq 1$ and $\norm{\lambda}$ is minimal among rational filtrations with this property (\Cref{theorem: Kempf}). We call $\lambda_{\text{b}}(x)$ the \emph{balanced filtration} of $x$. We can think of the balanced filtration as degenerating $x$ to its associated polystable point with optimal velocity and minimal cost. 

If we lift $x$ to the blow-up $\cB=\Bl_{\cX^{\max}}\cX$, then the balanced filtration of $x$ coincides with the filtration of $x$ associated to the natural $\Theta$-stratification $(\cS_c)_{c\in \Q_{\geq 0}}$ on $\cB$ from \Cref{theorem: theta-stratifications introduction}. If $\cE\subset \cB$ is the exceptional divisor, then we have a stratification of $\cX$ by type of balanced filtration where the strata are $\cX^{\max}$ and the $\cS_c\setminus \cE$. This is the first approximation to the iterated balanced filtration, that we get by iterating this procedure from the centres of the strata $\cS_c$.

The balancing stratification is indexed by a totally ordered set $\bGamma$ defined explicitly (\Cref{definition: totally ordered poset bGamma}). It consists of sequences $\left((d_0,c_0),\ldots,(d_n,c_n)\right)$ with $d_0, d_1, \cdots, d_n$ in $\N$, $c_0,\ldots,c_{n-1}\in \Q_{>0}$, $c_n=\infty$, and satisfying some other conditions, and the poset structure is given by lexicographic order. The balancing stratification is uniquely characterised by the following theorem. See \Cref{theorem: theta stratification proper over gms} for a precise formulation.


\begin{theorem}[Existence and characterisation of the balancing stratification]\label{theorem: main construction introduction}
There is a unique way of assigning, to every normed noetherian good moduli stack $\cX$ with affine diagonal, a sequential stratification $(\cS^\cX_\alpha)_{\alpha\in \bGamma}$ of $\cX$, called the \emph{balancing stratification} of $\cX$, in such a way that for every such $\cX$ the following holds. Let $(\cS_c)_{c\in \Q_{\geq 0}}$ be the $\Theta$-stratification of the blow-up $\cB=\Bl_{\cX^{\max}}\cX$ from \Cref{theorem: theta-stratifications introduction}, let $\cZ_c$ be the centres of the strata $\cS_c$, let $\cE\subset \cB$ be the exceptional divisor, let $d$ be the maximal stabiliser dimension of $\cX$, denote $\pi\colon \cX\to X$ the good moduli space, and let $\cU=\cX\setminus \pi^{-1}\pi(\cX^{\max})$. Then 
\begin{enumerate}
\item the highest stratum is $\cS^\cX_{(d,\infty)}=\cX^{\max}$, embedded in $\Filt_{\Q^\infty}(\cX)$ via the trivial filtration map $\cX\to \Filt_{\Q^\infty}(\cX)$;
\item for every $c\in \Q_{>0}$, and every $\alpha\in \bGamma$ with $\cS^{\cZ_c}_{\alpha}\neq \varnothing$, we have 
\[\cS^\cX_{((d,c),\alpha)}=\left(\cS_c\times_{\cZ_c}\cS^{\cZ_c}_\alpha\right)\setminus \cE\] 
with its natural structure of locally closed substack of $\Filt_{\Q^{\infty}}(\cX)$; and
\item for every $\alpha\in \bGamma$ with $\cS^{\cU}_\alpha\neq \varnothing$, we have $\cS_\alpha^\cX=\cS_\alpha^\cU$, as a substack of $\Filt_{\Q^\infty}(\cX)$.
\end{enumerate}
\end{theorem}

In the statement we are implicitly using that the $\cZ_c$ have good moduli spaces, which is \Cref{theorem: theta-stratifications introduction} above. Let us explain how to realise $\cS^\cX_{((d,c),\alpha)}$ inside $\Filt_{\Q^\infty}(\cX)$. For $c\in \Q_{>c}$, we can pullback each nonempty $\cS_\alpha^{\cZ_c}$ along the associated graded map $\cS_c\to \cZ_c$ to get a stack $\cV_{c,\alpha}=\cS_c\times_{\cZ_c} \cS^{\cZ_c}_\alpha$. The centre $\cZ_c$ lives as an open substack in $\Grad_\Q(\cB)$, so $\cS_\alpha^{\cZ_c}$ lives in $\Filt_{\Q^\infty}\left(\Grad_\Q(\cB)\right)$. The iterative definition of $\Filt_{\Q^\infty}(\cB)$ gives a cartesian square
\[\begin{tikzcd}
    {\Filt_{\Q^\infty}(\cB)} & {\Filt_{\Q^\infty}\left(\Grad_\Q(\cB)\right)} \\
    {\Filt_{\Q}(\cB)} & {\Grad_\Q(\cB)}
    \arrow[""{name=0, anchor=center, inner sep=0}, "\gr", from=2-1, to=2-2]
    \arrow[from=1-2, to=2-2]
    \arrow[from=1-1, to=1-2]
    \arrow[from=1-1, to=2-1]
    \arrow["\ulcorner"{anchor=center, pos=0.125}, draw=none, from=1-1, to=0]
\end{tikzcd}\]
and hence $\cV_{c,\alpha}$ lives in $\Filt_{\Q^\infty}(\cB)$. After subtracting the exceptional divisor $\cE$, we get a locally closed substack $\cS^\cX_{((d,c),\alpha)}\coloneqq \cV_{c,\alpha}\setminus \cE$  of $\Filt_{\Q^\infty}(\cX)$.

In order to prove \Cref{theorem: main construction introduction}, we introduce the concept of \emph{central rank} of a stack. The central rank of $\cX$, denoted $z(\cX)$, is the biggest natural number $n$ such that $B\G_m^n$ acts on $\cX$ in a nondegenerate way (\Cref{definition: central rank}). The condition can be thought of as every stabiliser of $\cX$ containing a copy of $\G_m^n$ in its centre. The maximal dimension of a stabiliser of $\cX$ is denoted $d(\cX)$. The proof is by induction on $N(\cX)\coloneqq d(\cX)-z(\cX)$.

The main observation is that, whenever we have a blow-up $f\colon \cY=\Bl_{\cR}\cX\to \cX$ of $\cX$ along some closed substack $\cR$, if the $\Theta$-stratification of $\cY$ is $(\cS_c)_{c\in \Q_{\geq 0}}$, then for every unstable stratum $\cS_c$, its centre $\cZ_c$ has bigger central rank than $\cX$: $z(\cZ_c)>z(\cX)$ (\Cref{lemma: centres of unstable strata have bigger central rank}). We also have $d(\cZ_c)\leq d(\cX)$ by representability of $f$, so that $N(\cZ_c)<N(\cX)$.

The balancing stratification defines a canonical sequential filtration for every point:

\begin{definition}[\Cref{definition: iterated balanced filtration}]
Let $\cX$ be a normed noetherian good moduli stack with affine diagonal. The \emph{iterated balanced filtration} of a field-valued point $x\in \cX(k)$ is the element $\lambda_{\text{ib}}(x)\in \qinfilt(\cX,x)$ determined by the balancing stratification of $\cX$.
\end{definition}

The iterated balanced filtration is defined over $k$ even if $k$ is not perfect. This has to do with the assumption that $\cX$ has a good moduli space instead of just an adequate moduli space \cite{_Alper_Adequatemodulispacesandgeometricallyreductivegroupschemes}, which is a more general notion in positive characteristic.

For a stack $\cY$ endowed with a $\Theta$-stratification $(\cS_c)$ such that all the centres $\cZ_c$ of the strata are normed good moduli stacks, the balancing stratification of each $\cZ_c$ can be pulled back to the $\cS_c$ to define a sequential stratification of $\cY$ (\Cref{remark: refining theta-stratifications}). This produces a \emph{refined Harder-Narasimhan filtration} for every point of $\cY$.

We show that the balancing stratification has nice functorial properties. 

\begin{proposition}[\Cref{proposition: compatibility of balancing stratification with pullback}]
Let $f\colon \cX\to \cX'$ be a morphism between noetherian normed good moduli stacks with affine diagonal. If $f$ is either a closed immersion or a base change from a map $X\to X'$ between the good moduli spaces of $\cX$ and $\cX'$, then for all $\alpha\in\bGamma$, the stratum $\cS^\cX_\alpha$ equals the pullback $\cS^{\cX'}_\alpha\times_{\cX',f}\cX$ of $\cS^{\cX'}_\alpha$ along $f$, with its natural structure of locally closed substack of $\Filt_{\Q^\infty}(\cX)$.
\end{proposition}
\subsection{Relation to convex geometry and artinian lattices}
Despite its seemingly convoluted definition involving several blow-ups and $\Theta$-stratifications, the iterated balanced filtration has a particularly simple description for a point $x$ in a quotient stack of the form $\Spec A/G$ with $G$ a diagonalisable algebraic group over a field $k$ (for example a split torus $\G_{m,k}^n$) and $A$ a finite type $k$-algebra.
For simplicity, in this introduction we will assume that $\Spec A$ is the total space of a vector space $V=k^l$, that $G$ acts on $V$ via the characters $\chi_1,\ldots,\chi_l\in \Gamma_\Z(G)=\Hom(G,\G_{m,k})$ and that we are interested in computing the iterated balanced filtration of the point $x=(1,\ldots,1)$ in the quotient stack $V/G$. We denote $N=\Gamma^\Q(G)$ the set of rational cocharacters of $G$. The norm on graded points of $V/G$ corresponds to a rational inner product $(-,-)$ on $N=N_0$, and we identify $N$ and its dual via this inner product. 
  
We now describe a sequence of elements $\lambda_0,\ldots,\lambda_n\in N$ in terms of the \emph{state} $\Xi_0=\{\chi_1,\ldots,\chi_l\}$ of $x$. Let $F_0$ be the smallest face containing $0$ of the cone $\cone(\Xi_0)$ generated by $\Xi_0$ inside $N$. Then we let $\lambda_0$ be the unique element of the orthogonal complement $N_1\coloneqq F_0^\bot$ such that $(\lambda_0,\alpha)\geq 1$ for all $\alpha\in \Xi_0\setminus F_0$ and $\norm{\lambda_0}$ is minimal. This is a convex optimisation problem.

To compute $\lambda_1,\cdots,\lambda_n$, we proceed  
as follows. We let $\Xi_1=\{p_1(\alpha)\in \Xi_0\st (\lambda_0,\alpha)=1\}\subset N_1$, where $p_1\colon N\to N_1$ is the orthogonal projection. It will always be the case that $\lambda_0\in \cone\left(\Xi_1\right)\subset N_1$ (\Cref{theorem: Recognition of the balanced filtration for states}), and thus there is a smallest face $F_1$ of $\cone\left(\Xi_1\right)$ containing $\lambda_0$. Let $N_2$ be the orthogonal complement of $F_1$ inside $N_1$. Then we let $\lambda_1$ be the unique element in $N_2$ such that $(\lambda_1,\alpha)\geq 1$ for all $\alpha\in \Xi_1\setminus F_1$ and $\norm{\lambda_1}$ is minimal. We define $\Xi_2=\{p_2(\alpha)\in \Xi_1\st (\lambda_1,\alpha)=1\}\subset N_2 $, where $p_2\colon N_1\to N_2$ is the orthogonal projection. Repeating this process, we get $\lambda_2,\ldots,\lambda_n$. The algorithm terminates when we get to $\Xi_{n+1}\subset F_{n+1}$.

A $\Q^\infty$-filtration of $x$ in $V/G$ is determined by a finite sequence of elements of $N$ (\Cref{remark: description sequential filtrations quotient stacks}), and under this correspondence we have:

\begin{theorem}[\Cref{theorem: balancing chain of states corresponds to torsor chain} and \Cref{corollary: iterated balanced filtration for states equals that for stacks}]
The iterated balanced filtration of $x$ in $V/G$ is given by the sequence $\lambda_0,\ldots,\lambda_n$ described above.
\end{theorem}

For the proof, we use the machinery of \emph{chains of stacks}. A chain of stacks is the data of a sequence of pointed $k$-stacks $(\cX_n,x_n)$, together with a $\Q$-filtration $\lambda_n$ of each $x_n$ and link maps $(\cX_{n+1},x_{n+1})\to \left(\Grad(\cX_n),\gr \lambda_n\right)$ (\Cref{definition: chain of stacks}). Associated to every chain there is a $\Q^\infty$-filtration of the point $x_0$ in $\cX_0$ (\Cref{definition: Qinfinity filtration of a chain}). For a pointed normed good moduli stack $(\cX,x)$, we give two different constructions of chains computing sequential filtrations of $(\cX,x)$, the \emph{balancing chain} (\Cref{construction: balancing chain}) and the \emph{torsor chain} (\Cref{construction: torsor chain}). The former is closely related to the balancing stratification, while the latter tends to be closer to combinatorial structures, like states or lattices. We show that both chains compute the iterated balanced filtration of $(\cX,x)$ (\Cref{proposition: balancing chain computes iterated balanced filtration} and \Cref{theorem: torsor computes the iterated balanced filtration}).

In order to relate the convex-geometrical picture of states to chains of stacks, we define a category of normed semistable polarised states (\Cref{definition: polarised state}), combinatorial in nature. We then define an analogue of the notion of chain in this category, and a canonical \emph{balancing chain} for every object. There is a functor from normed semistable polarised states to pointed normed good moduli stacks (\Cref{definition: stack associated to a state,definition: map of states defines map between stacks}), and we show that it sends the balancing chain of a state to the torsor chain of the corresponding stack (\Cref{theorem: balancing chain of states corresponds to torsor chain}).

The theory of chains of stacks will be used to establish, in a sequel to this paper \cite{IbanezNunez_Blowupsandlattices}, a correspondence between the iterated balanced filtration and the HKKP filtration for normed artinian lattices as defined by Haiden-Katzarkov-Kontsevich-Pandit \cite{_Haiden_Semistabilitymodularlatticesanditeratedlogarithms}. 

\begin{theorem}[see \cite{IbanezNunez_Blowupsandlattices}]\label{theorem: correspondence HKKP filtration and iterated balanced filtration}
Let $k$ be an algebraically closed field and let $\cA$ be a locally noetherian $k$-linear Grothendieck abelian category. Suppose that the moduli stack of objects $\cM_\cA$ is an algebraic stack locally of finite type over $k$, and let $\cX$ be a quasi-compact open substack of $\cM_\cA$ admitting a good moduli space $\pi\colon \cX\to X$, and endowed with a linear norm on graded points. For any $k$-point $x\in \cX(k)$, there is a canonically defined normed artinian lattice $L_x$ and a canonical bijection $\qinfilt(L_x)\cong\qinfilt(\cX,x)$ under which the HKKP filtration of $L_x$ and the iterated balanced filtration of $(\cX,x)$ agree.
\end{theorem}
Here, by an artinian lattice $L$ we mean a lattice that is modular and of finite length. The set $\qinfilt(L)$ of sequential filtrations of $L$ is defined as in the case of vector bundles above.
The norm being linear means it is compatible with the underlying abelian category in a precise sense. Examples of this setup include moduli spaces of Bridgeland semistable objects, moduli spaces of semistable vector bundles on a curve $C$, and moduli of quiver representations (\Cref{section: examples}). For abelian categories $\cA$ satisfying suitable finiteness condition, algebraicity of $\cM_\cA$ has been established in \cite{_Fernandez-Herrero_Lennen_Makarva_Modulifinitelengthabeliancategories}. In that context, choices of K-theory classes of the category $\cA$ produce open substacks $\cX\subset \cM_\cA$ having a good moduli space and a norm on graded points \cite[Theorem 4.17]{_Fernandez-Herrero_Lennen_Makarva_Modulifinitelengthabeliancategories}, so \Cref{theorem: correspondence HKKP filtration and iterated balanced filtration} applies.
 For the proof of \Cref{theorem: correspondence HKKP filtration and iterated balanced filtration}, we will need a new characterisation of the HKKP filtration for artinian lattices established in \cite{IbanezNunez_CharacterisationHKKP}.

\subsection{Asymptotics of flows} 
In view of \Cref{theorem: correspondence HKKP filtration and iterated balanced filtration} and the results of Haiden-Katzarkov-Kontsevich-Pandit on asymptotics of flows in the case of quiver representations and vector bundles on a smooth projective complex algebraic curve \cite{_Haiden_Semistabilitymodularlatticesanditeratedlogarithms,_Haiden_Iteratedlogarithmsandgradientflows}, we expect the iterated balanced filtration to play a role in describing asymptotics of natural flows in moduli theory. We now make this expectation precise in the case of Geometric Invariant Theory on affine spaces.

We begin by recalling the Kempf-Ness potential in a more general framework. Let $G$ be a connected reductive algebraic group over $\C$, endowed with a norm on cocharacters, that is, the data of a Weyl-invariant rational inner product on the set $\Gamma^\Z(T)$ of cocharacters of a maximal torus $T$ of $G$ (\Cref{definition: norm on cocharacters of a group}). Let $K$ be a maximal compact subgroup of $G$. We denote $\fg$ the Lie algebra of $G$ and $\fk$ the Lie algebra of $K$. We consider a smooth projective-over-affine scheme $X$ over $\C$, endowed with an action of $G$. An ample line bundle $\cL$ on $X/G$ (that is, an ample line bundle $L$ on $X$ with a $G$-equivariant structure) defines an open semistable locus $(X/G)^\ss=X^\ss/G$, which is the open stratum of a $\Theta$-stratification of $X/G$ (\Cref{theorem: theta stratification proper over gms}).
The line bundle $\cL$ restricts to a line bundle $\overline\cL$ on the differentiable stack “of metrics” $X/K$, and we endowed $\overline\cL$ with a hermitian norm $\norm{-}$ (that is, we endow $L$ with a $K$-invariant hermitian norm). For a point $x\in X(\C)$, after a choice of nonzero lift $x^*$ of $x$ to the total space of $\overline \cL$, the \emph{Kempf-Ness potential} is defined to be
\[p_{ x}\colon G/K\to \R\colon Kg\mapsto \log\norm{x^*}-\log \norm{gx^*}.\]
Here, $G/K$ denotes the quotient by the action of $K$ on $G$ on the left. Note that the Kempf-Ness potential is independent of the choice of lift $x^*$ of $x$. 
Kempf-Ness type theorems state that, under some conditions, the following equivalences hold:
\begin{enumerate}
\item $x$ is semistable if and only if $p_x$ is bounded below.
\item $x$ is polystable if and only if $p_x$ attains a minimum.
\end{enumerate}
This is the case, for example, if either $X$ is projective \cite{_Kempf_Ness,_Kirwan_Cohomologyofquotientsinsymplecticandalgebraicgeometry,the-moment-weight-inequality} or $X$ is affine and some additional conditions are satisfied \cite{_King_Moduliofrepresentationsoffinitedimensionalalgebras,Hoskins_AffineGIT,Mayrand_AffineKempfNess}. In these cases, if $x$ is polystable, then the negative gradient flow of $p_x$ converges to a minimum from any starting point. In the strictly semistable case, we are interested in understanding the asymptotic behaviour of this flow.

In order to define the gradient flow, we need a Riemannian metric on $G/K$. This comes from the norm on cocharacters of $G$, which gives a $K$-invariant euclidean inner product on $\fk$. By Hadamard's theorem, the map 
\[\phi\colon \fk\to G/K\colon v\mapsto K\exp\left(\dfrac{1}{2i}v\right)\]
is a diffeomorphism (see \cite[Appendix A]{the-moment-weight-inequality}). The isomorphism $d_0\phi\colon \fk\to T_e (G/K)$ induces an inner product on $T_e (G/K)$, and we extend it to a Riemannian metric on $G/K$ using the action of $G$ on $G/K$ by right translations (see \cite[Appendix C]{the-moment-weight-inequality} for details). This allows us to define the gradient vector field $\nabla p_x$ on $G/K$.

The semistable locus $X^\ss/G$ is a good moduli stack, endowed with a norm on graded points coming from the norm on cocharacters of $G$. If $x\in X(\C)$ is a semistable point, then the iterated balanced filtration $\lambda_{\text{ib}}(x)$ of $x$ is defined, and by \Cref{remark: description sequential filtrations quotient stacks} it is identified with an equivalence class of sequences of commuting rational one-parameter subgroups $\lambda_1,\ldots,\lambda_n$.
We may choose representatives $\lambda_j$ that are compatible with $K$, in the sense that for some power $\lambda_j^l\colon \C^\times \to G$ that is integral (and hence for any such power), the inclusion $\lambda_j^l(S^1)\subset K$ holds. The $\lambda_j$ induce well-defined linear maps on Lie algebras $\Lie(\lambda_j)\colon \C\to \fg$, giving elements $\nu_j=i\Lie(\lambda_j)(1)\in \fk$. 

We now shift our attention to a particular case of the above setting. We will take $X=V$ to be a finite dimensional $G$-representation. Let $\alpha\colon G\to \G_{m,\C}$ be a character, and take the linearisation $\cL=\cO_{X/G}(\alpha)=(X/G\to B\G_{m,\C})^*\alpha$, where we regard $\alpha$ as an element of $\Pic(BG)$. Choose a $K$-invariant hermitian metric $\norm{-}$ on $V$. The total space of $\cL$ is $\left(V\times \C\right)/G$, where the action is $g(x,c)=(gx,\alpha(g)c)$. The norm $\norm{-}$ defines a hermitian norm $\norm{-}_{\overline \cL}$ on $\overline \cL$ by the formula
\[\norm{(v,c)}_{\overline \cL}=e^{-\norm{v}^2}\abs{c}.\]
The associated Kempf-Ness potential for a point $x\in V$ is
\[p_x(Kg)=\norm{gx}^2-\log \abs{\alpha(g)} -\norm{x}^2.\]

Suppose that $x\in V$ is semistable and let $\nu_1,\ldots,\nu_n\in \fk$ represent the iterated balanced filtration of $x$ as above. We denote $\log=\phi^{-1}\colon G/K\to \fk$ the inverse of the map $\phi$ defined above. With this setup, we conjecture:

\begin{conjecture}\label{conjecture: on asymptotics of gradient flows affine GIT}
Let $h\colon (0,\infty)\to G/K$ be a flow line for $-\nabla p_x$. Then the expression
\[\log h(t)+\log(t)\nu_1+\log\log(t)\nu_2+\cdots+\underbrace{\log\cdots \log}_{n}(t)\nu_n\]
in $\fk$ is bounded for $t>>0$.
\end{conjecture}

In the case where $V=\bigoplus_{a\in Q_1}\Hom(\C^{d_{s(a)}},\C^{d_{t(a)}})$ is the representation space of a quiver $Q$ with dimension vector $d$ and $G=\prod_{i\in Q_0}\GL_{d_i,\C}$ with the standard action, the conjecture is true for specific choices of hermitian norm on $V$ and norm on cocharacters of $G$ by \Cref{theorem: correspondence HKKP filtration and iterated balanced filtration} and \cite[Theorem 5.11]{_Haiden_Semistabilitymodularlatticesanditeratedlogarithms}. See \Cref{example: simple case of conjecture on flows} for a simple example where the conjecture is checked beyond the quiver case. We hope to return to \Cref{conjecture: on asymptotics of gradient flows affine GIT} in future work.

Our expectation is that the iterated balanced filtration describes also the asymptotics of natural gradient flows in other moduli problems. Examples include the Calabi flow for a K-semistable Fano variety, the Yang-Mills flow for semistable $G$-bundles on a smooth projective curve and the gradient flow for the Kempf-Ness potential in more general GIT situations. In all these examples, there is an underlying normed good moduli stack of semistable objects, so the iterated balanced filtration is defined (\Cref{section: examples}).

\subsection{Notation and conventions}
The set of nonnegative integers, or natural numbers, is denoted $\N=\Z_{\geq 0}$. We will also denote $\N^*=\Z_{>0}$.

We follow the definitions and conventions of \cite{stacks-project} regarding algebraic stacks and algebraic spaces. We denote $\cat{St_{\fppf}}$ the 2-category of stacks for the fppf site of schemes. For an algebraic stack $\cX$, we denote $\abs{\cX}$ its topological space and $\pi_0(\cX)$ the set of connected components of $\abs{\cX}$. If $x\colon T\to \cX$ a $T$-point, with $T$ a scheme, we denote $\Aut(x)$ the automorphism group of $x$ as a group algebraic space over $T$. The multiplicative group over $\Z$ is denoted $\G_m=\Spec \Z[t,t^{-1}]$ and, for a scheme $T$, we use the notation $\G_{m,T}=\G_m\times T$.

For a group scheme $G$ acting on an algebraic space $X$, we use the notation $X/G$ for the quotient stack, omitting the customary brackets.

For a field $k$ and an algebraic group $G$ over $k$, we denote $\Gamma_\Z(G), \Gamma^\Z(G), \Gamma_\Q(G)$ and $\Gamma^\Q(G)$ the sets of characters, cocharacters, rational characters and rational cocharacters of $G$, respectively. If $\lambda\in \Gamma^\Z(G)$ is a cocharacter, and $G$ acts on a scheme $X$ over $k$, then we denote $X^{\lambda,0}$ the fixed point locus of the induced $\G_m$-action on $X$ and $X^{\lambda,+}$ the \emph{attractor}, defined functorially on $k$-schemes $T$ by the formula $\Hom(T,X^{\lambda,+})=\Hom^{\G_{m,k}}(\A^1_T,X)$, where $\Hom^{\G_{m,k}}$ denotes $\G_{m,k}$-equivariant maps and $\A^1_k$ is endowed with the usual scaling action \cite{_Drinfeld_OnalgebraicspaceswithanactionofGm}. For the particular case of the conjugation action of $G$ on itself, we denote $L(\lambda)=G^{\lambda,0}$ and $P(\lambda)=G^{\lambda,+}$. If $G$ is reductive, then $P(\lambda)$ is a parabolic subgroup with Levi factor $L(\lambda)$. If $g\in G(k)$, we denote $\lambda^g=g\lambda g^{-1}$.

If $\cF$ is a vector bundle on an algebraic stack $\cX$, the total space of $\cF$ is $\A(\cF)\coloneqq\Spec_{\cX}\Sym_{\cO_\cX} \cF^\vee$ and the associated projective bundle is $\P(\cF)\coloneqq \Proj_\cX\Sym_{\cO_\cX} \cF^\vee$.

\subsection{Acknowledgements}
 I wish to express sincere gratitude to my PhD advisor Frances Kirwan for her constant help and encouragement throughout this project, and for introducing me to many interesting mathematics related to this work. I also sincerely thank my second advisor Fabian Haiden for many useful discussions and for his continuous support. Special thanks go to Daniel Halpern-Leistner for explaining many aspects of the Beyond GIT programme. Finally, I would like to thank as well Lukas Brantner, Michel Brion, Chenjing Bu, Ben Davison, Ruadhaí Dervan, Andres Fernandez Herrero, Oscar García-Prada, Tasuki Kinjo and David Rydh for helpful conversations related to this project.

\section{Preliminaries}
In this section we start by recalling, mainly following \cite{_HalpernLeistner_Onthestructureofinstabilityinmodulitheory}, the kind of stability structures on algebraic stacks that we will use in the rest of the paper and how they give rise to stratifications. The two main concepts are that of a \emph{norm on graded points} (\Cref{definition: norm on graded points}) and of a \emph{linear form on graded points} (\Cref{definition: linear form on stack}), and these give rise to \emph{$\Theta$-stratifications} (\Cref{definition: Theta-stratification my version}).
While in \cite{_HalpernLeistner_Onthestructureofinstabilityinmodulitheory} the stacks $\Grad(\cX)$ of graded points and $\Filt(\cX)$ of filtrations of an algebraic stack $\cX$ are used, in this work we will need a generalisation of these, what we call the stack of \emph{rational graded points} $\Grad_\Q(\cX)$ and the stack of $\emph{rational filtrations}$ $\Filt_\Q(\cX)$, that we define in \Cref{subsection: stacks of rational filtrations and graded points}.

The main result in this section is \Cref{theorem: theta stratification proper over gms}, where we prove that, for a representable projective morphism $f\colon \cY\to \cX$ into a stack $\cX$ with a good moduli space and a norm on graded points, the stack $\cY$ carries a natural $\Theta$-stratification $(\cS_c)_{c\in \Q_{\geq 0}}$ whose centres $(\cZ_c)_{c\in \Q_{\geq 0}}$ have good moduli spaces. This is an improvement of Halpern-Leistner's result \cite[Theorem 5.5.10]{_HalpernLeistner_Onthestructureofinstabilityinmodulitheory} that $\cY$ has a weak $\Theta$-stratification whose semistable locus $\cZ_0$ has a good moduli space. Existence of good moduli spaces for all centres $\cZ_c$ will be fundamental in our construction of the balancing stratification (\Cref{theorem: sequential stratifications for good moduli stacks}).

\subsection{Good moduli spaces and local structure theorems}
We start by recalling the definition of good moduli space from \cite[Definition 4.1]{_Alper_GoodmodulispacesforArtinstacks}, with the slightly modified conventions of \cite[1.7.3, 1.7.4]{_Alper_Theetalelocalstructureofalgebraicstacks}.
\begin{definition}[Good moduli space]\label{definition: good moduli space}
A morphism $\pi\colon \cX\to X$ from an algebraic stack $\cX$ to an algebraic space $X$ is said to be an \emph{good moduli space} if
\begin{enumerate}
\item the map $\pi$ is quasi-compact and quasi-separated;
\item the map $\cO_X\to \pi_*\cO_{\cX}$ is an isomorphism; and
\item the pushforward functor $\pi_*$ on quasi-coherent sheaves is exact, and the same is true after any base change $X'\to X$, where $X'$ is an algebraic space.
\end{enumerate}
\end{definition}

If $X$ is quasi-separated, then the third condition can be simplified to $\pi_*$ being exact, the statement for any base change being then automatic by \cite[Proposition 3.10, (vii)]{_Alper_GoodmodulispacesforArtinstacks}. We will often use the term \emph{good moduli stack} meaning an algebraic stack $\cX$ that admits a good moduli space $\pi\colon \cX\to X$.

\begin{remark}
Using good moduli spaces, we can recover the concept of linear reductivity. Indeed, an affine algebraic group $G$ over a field $k$ is linearly reductive precisely when the map $BG\to \Spec k$ is a good moduli space.
\end{remark}

A good moduli space $\pi\colon \cX\to X$ enjoys many special properties, for example:
\begin{enumerate}
\item Any map $\cX\to Y$ with $Y$ an algebraic space factors uniquely through $\pi$ \cite[Theorem 3.12]{_Alper_Theetalelocalstructureofalgebraicstacks}. In particular, the good moduli space $\pi$ is uniquely determined by $\cX$.
\item Any base change of $\pi$ along a morphism $X'\to X$ with $X'$ an algebraic space is a good moduli space \cite[Proposition 4.7, (i)]{_Alper_GoodmodulispacesforArtinstacks}.
\item If $h\colon \cX'\to \cX$ is an affine morphism, then $\cX'$ has a good moduli space $\cX'\to X'$ and the induced map $X'\to X$ is affine with $X'=\Spec_X \pi_*h_*\cO_{\cX'}$ \cite[Lemma 4.14]{_Alper_GoodmodulispacesforArtinstacks}.
\item For every point $p\in \abs{X}$, the fibre $\pi^{-1}(p)$ has a unique closed point $q$, and the dimension of the stabiliser of $q$ is bigger than that of any other point of $\pi^{-1}(p)$ \cite[Proposition 9.1]{_Alper_GoodmodulispacesforArtinstacks}. Moreover, the stabiliser of $q$ is linearly reductive \cite[Proposition 12.14]{_Alper_GoodmodulispacesforArtinstacks}. The points $q\in \abs{\cX}$ that are closed in the fibre of $\pi$ containing $q$ are said to be \emph{polystable}.
\end{enumerate}

Stacks with good moduli spaces are étale locally quotient stacks. More precisely:

\begin{theorem}[Local structure {\cite{_Alper_Theetalelocalstructureofalgebraicstacks}}]\label{theorem: local structure}
Let $\cX$ be an algebraic stack and $\pi\colon \cX\to X$ a good moduli space. Assume that $\cX$ is of finite presentation over a quasi-compact and quasi-separated algebraic space $B$ and that $\cX$ has affine diagonal.

Then there is a natural number $n$, an affine scheme $\Spec A$ endowed with an action of $\GL_n$, and a cartesian square
\[\begin{tikzcd}[ampersand replacement=\&]
    {(\Spec A)/\GL_n} \& \cX \\
    {\Spec(A^{\GL_n})} \& X
    \arrow["\pi", from=1-2, to=2-2]
    \arrow["h", from=2-1, to=2-2]
    \arrow[from=1-1, to=2-1]
    \arrow[from=1-1, to=1-2]
    \arrow["\lrcorner"{anchor=center, pos=0.125}, draw=none, from=1-1, to=2-2]
\end{tikzcd}\]
with $h$ an affine Nisnevich cover (in particular, étale). Here, $A^{\GL_n}$ denotes the ring of invariants. Moreover, $X\to B$ is of finite presentation and $X$ has affine diagonal.
\end{theorem}

The theorem is \cite[Theorem 6.1]{_Alper_Theetalelocalstructureofalgebraicstacks}, together with the argument at the end of the proof of \cite[Theorem 5.3]{_Alper_Theetalelocalstructureofalgebraicstacks} to guarantee that $h$ can be taken to be affine. To see that $X$ has affine diagonal, just take good moduli spaces for the diagonal $\cX\to \cX\times \cX$, which is affine, to obtain the diagonal of $X$.

From \Cref{theorem: local structure}, it follows that stacks whose good moduli space is a point are necessarily quotient stacks.

\begin{corollary}\label{corollary: fibres of good moduli space are quotient stacks}
Let $\cX$ be an algebraic stack of finite presentation over a field $k$ and assume that $\pi\colon \cX\to \Spec k$ is a good moduli space with affine diagonal. Then $\cX\cong (\Spec A)/\GL_n$, where $A$ is a $k$-algebra of finite type and $\Spec A$ is endowed with a $\GL_n$-action.
\end{corollary}

Over an algebraically closed field, there is a stronger result.

\begin{corollary}[of {\cite[Theorem 4.12]{_Alper_ALunaetaleslicetheoremforalgebraicstacks}}]\label{corollary: stack with good moduli space a point is a quotient algebraically closed case}
Let $\cX$ be an algebraic stack of finite presentation over an algebraically closed field $k$ and suppose that $\pi\colon \cX\to \Spec k$ is a good moduli space with affine diagonal. Let $x\in \cX(k)$ be the unique closed $k$-point of $\cX$ and let $G$ be the stabiliser of $x$. Then $\cX\cong (\Spec A)/G$, where $A$ is a finite type $k$-algebra and $\Spec A$ is endowed with an action of $G$.
\end{corollary}
\subsection{Stacks of rational filtrations and graded points}\label{subsection: stacks of rational filtrations and graded points}
In \cite{_HalpernLeistner_Onthestructureofinstabilityinmodulitheory}, Halpern-Leistner defines the stacks $\Grad(\cX)$ of graded points  and $\Filt(\cX)$ of filtrations of an algebraic stack $\cX$ as mapping stacks. If $\cX$ parametrises objects in an abelian category, then $\Grad(\cX)$ parametrises objects endowed with a $\Z$-grading and $\Filt(\cX)$ parametrises objects endowed with a $\Z$-filtration \cite[Proposition 7.12 and Corollary 7.13]{_Alper_Existenceofmodulispacesforalgebraicstacks}, but $\Grad(\cX)$ and $\Filt(\cX)$ can be defined for very general $\cX$. In this section we revisit the construction of $\Grad(\cX)$ and $\Filt(\cX)$, and extend it to consider rational filtrations and gradings.

Following \cite{_HalpernLeistner_Onthestructureofinstabilityinmodulitheory}, we define the stack $\Theta$ over $\Spec(\Z)$ to be the quotient stack $\Theta=\A_\Z^1/\G_{m,\Z}$, the action of $\G_{m,\Z}$ on $\A^1_{\Z}$ being the usual scaling action. For an algebraic space $S$ we denote $\Theta_S=\Theta\times S$. 

Recall that for $\cY$ and $\cZ$ two stacks over a base space $S$, an object of the \emph{mapping stack} $\uHom_S(\cY,\cZ)$ over a scheme $T$ is a map $T\to S$ together with a morphism $T \times_S\cY\to \cX$ over $S$. In the case $S=\Spec(\Z)$, we omit the subindex from the notation. The following definition is in \cite[Section 1.1]{_HalpernLeistner_Onthestructureofinstabilityinmodulitheory}, except that we do not work relative to a base algebraic stack.

\begin{definition}[Stacks of filtrations and graded points]\label{definition: stacks of filtations and graded points}
Let $\cX$ be an algebraic stack and let $n$ be a positive integer. We define the stack $\Grad^n(\cX)$ of \emph{$\Z^n$-graded points} of $\cX$ to be the mapping stack $\Grad^n(\cX)\coloneqq \uHom(B\G_{m,\Z}^n,\cX)$. Similarly, we define the stack $\Filt^n(\cX)$ of \emph{$\Z^n$-filtrations} of $\cX$ to be $\Filt^n(\cX)\coloneqq \uHom(\Theta^n, \cX)$.
\end{definition}

We will simply denote $\Filt(\cX)=\Filt^1(\cX)$ and $\Grad(\cX)=\Grad^1(\cX)$. 

\begin{lemma}[Independence of base for mapping stacks]
Let $\cY$ be an algebraic stack such that $\cY\to \Spec \Z$ is a good moduli space, and let $\cX$ be an algebraic stack defined over an algebraic space $B$. Then there is a canonical isomorphism
\[\uHom(\cY,\cX)\cong\uHom_B(\cY_B,\cX)\]
of mapping stacks. In particular, there is a canonical map $\uHom(\cY,\cX)\to B$.
\end{lemma}
\begin{proof}
For a scheme $T$, an object of the groupoid $\uHom_B(\cY_B,\cX)(T)$ is a pair $(a,b)$ with $b\colon T\to B$ and $a\colon T\times_B \cY_B=T\times \cY\to \cX$ a morphism over $B$. Since $T\times \cY\to T$ is a good moduli space, for any given $a\colon T\times \cY\to \cX$, an object of $\uHom(\cY,\cX)(T)$, the composition $T\times \cY\to \cX\to B$ factors uniquely through $T\times \cY\to T$, giving a unique $b\colon T\to B$ such that $(a,b)$ is an object of $\uHom_B(\cY_B,\cX)(T)$.
\end{proof}

Applying the Lemma for $\cY=\Theta^n$ or $\cY=B\G_m^n$, we see that $\Filt^n(\cX)$ and $\Grad^n(\cX)$ are independent of the base algebraic space considered.

To guarantee that $\Filt^n(\cX)$ and $\Grad^n(\cX)$ are well-behaved, we consider the following assumption on an algebraic stack $\cX$ defined over an algebraic space $B$.

\begin{assumption}\label{assumption: basic assumptions}
The algebraic space $B$ is quasi-separated and locally noetherian, and the map $\cX\to B$ is locally finitely presented and has affine diagonal.
\end{assumption}

\begin{example}\label{example: noetherian good moduli stack with affine diagonal}
Suppose that $\cX$ is a noetherian algebraic stack with affine diagonal and $\pi\colon \cX\to X$ is a good moduli space. Then $X$ is also noetherian by \cite[Theorem 4.16]{_Alper_GoodmodulispacesforArtinstacks}, and $\pi$ is of finite type by \cite[Theorem A.1]{_Alper_ALunaetaleslicetheoremforalgebraicstacks}. The diagonal of $X$ is affine since it is obtained from the diagonal of $\cX$ by taking good moduli spaces. In particular, $\cX$ satisfies \Cref{assumption: basic assumptions} with $B=X$.
\end{example}

Under Assumption \ref{assumption: basic assumptions}, the stacks $\Filt^n(\cX)$ and $\Grad^n(\cX)$ are algebraic and also satisfy \Cref{assumption: basic assumptions} \cite[Theorem~6.22]{_Alper_Theetalelocalstructureofalgebraicstacks}. Note that by \cite[Remark 6.16]{_Alper_Theetalelocalstructureofalgebraicstacks} it is not necessary to assume that $B$ is excellent, since the stacks $B\G_{m,B}^n$ and $\Theta_B^n$ satisfy condition (N) in \cite{_Alper_Theetalelocalstructureofalgebraicstacks}. 
See also \cite[Theorem 5.1.1]{_HalpernLeistner_Mappingstacksandcategoricalnotionsofproperness} for a related algebraicity result.

There are several maps relating $\Grad(\cX)$, $\Filt(\cX)$ and $\cX$:
\begin{enumerate}
    \item The “evaluation at 1” map $\ev_1\colon \Filt(\cX)\to \cX$, defined by precomposition along $\{1\}\to \Theta$. It is representable and separated \cite[Proposition 1.1.13]{_HalpernLeistner_Onthestructureofinstabilityinmodulitheory}.
    \item The “associated graded” map $\gr\colon \Filt(\cX)\to \Grad(\cX)$, defined by precomposition along $B\G_m=\{0\}/\G_m\to \Theta$.
    \item The “forgetful” map $u\colon \Grad(\cX)\to \cX$, defined by precomposition along $\Spec(\Z)\to B\G_m$.
    \item The “evaluation at 0” map $\ev_0\colon \Filt(\cX)\to \cX$, which is the composition $\ev_0=u\circ \gr$.
    \item The “split filtration” map $\sigma\colon \Grad(\cX)\to \Filt(\cX)$, defined by precomposition along the canonical representable morphism $\Theta\to B\G_m$.
    \item The “trivial grading” map $\cX\to \Grad(\cX)$, given by precomposition along $B\G_m\to \Spec \Z$. It is an open and closed immersion.
    \item The “trivial filtration” map $\cX\to \Filt(\cX)$, defined by precomposing along $\Theta\to \Spec \Z$. It is an open and closed immersion.
\end{enumerate}
\begin{remark}
The fact that the “trivial grading”  and “trivial filtration”  maps are closed and open immersions follows from \cite[Proposition 1.3.9]{_HalpernLeistner_Onthestructureofinstabilityinmodulitheory} by the argument in \cite[Proposition 1.3.11]{_HalpernLeistner_Onthestructureofinstabilityinmodulitheory}.
\end{remark}

\begin{remark}
Sometimes, the assumption that $\cX\to B$ has affine diagonal can be relaxed to $\cX$ having affine stabilisers and being quasi-separated over $B$. This is the case for the construction and algebraicity of $\Filt_\Q^n(\cX)$, $\Grad_\Q^n(\cX)$, $\Filt_{\Q^\infty}(\cX)$ and $\Grad_{\Q^\infty}(\cX)$ (\Cref{definition: rational filtrations and gradings} and \Cref{definition: Q^infty filtrations and gradings}). Representability of $\ev_1\colon \Filt(\cX)\to \cX$ follows under the additional hypothesis that $\cX$ has separated inertia.

The assumption that $B$ is locally noetherian guarantees that the topological spaces of the algebraic stacks considered are locally connected, and hence their connected components are open. 
\end{remark}

The monoid $(\N^*,\cdot, 1)$ acts\footnote{Formally, an action of $(\N^*,\cdot, 1)$ on a stack $\cY$ is a pseudofunctor $B(\N^*,\cdot, 1)\to \cat{St_{\fppf}}$ sending the unique object of $B(\N^*,\cdot, 1)$ to $\cY$. Here, we are denoting $B(\N^*,\cdot, 1)$ the category with one object and endomorphism monoid equal to $(\N^*,\cdot, 1)$, and $\cat{St_\fppf}$ is the 2-category of stacks on the category of schemes with the fppf topology.} on the stacks $\Filt(\cX)$ and $\Grad(\cX)$. A natural number $n>0$ acts on $\Filt(\cX)$ by the map 
\begin{tikzcd}[column sep=large]
    {\Filt(\cX)} & {\Filt(\cX)}
    \arrow["{(\bullet)^n}", from=1-1, to=1-2]
\end{tikzcd}
given by precomposition along the $n$th power map $\Theta_T\to \Theta_T$, and similarly in the case of $\Grad(\cX)$. Now denote $\cY$ one of the stacks $\Grad(\cX)$ or $\Filt(\cX)$ with its $(\N^*,\cdot, 1)$-action. We define a diagram $D_\cY\colon (\N^*,\vert)\to\cat{St_\fppf}$, i.e. a pseudofunctor, on the 2-category $\cat{St_\fppf}$ of stacks for the fppf site of schemes. The index category is the filtered poset $(\N^*,\vert)$ of positive integers with the divisibility order, and $D_{\cY}$ is defined by setting $D_\cY(n)=\cY$ for all $n$, and $D_{\cY}(n\vert m)$ to be the “rising to the $\frac{m}{n}$th power” map \begin{tikzcd}[column sep=large]
    {\cY} & {\cY}
    \arrow["{(\bullet)^{m/n}}", from=1-1, to=1-2]
\end{tikzcd} defined above.
\begin{definition}[Stacks of rational filtrations and rational graded points]\label{definition: rational filtrations and gradings}
The stacks $\Filt_\Q(\cX)$ of \emph{rational filtrations} and $\Grad_\Q(\cX)$ of \emph{rational graded points} are the colimits
\[\Filt_\Q(\cX)\coloneqq \colim D_{\Filt(\cX)} \quad \text{and}\quad \Grad_\Q(\cX)\coloneqq \colim D_{\Grad(\cX)}\]
in the cocomplete 2-category $\cat{St_\fppf}$.
\end{definition}

\begin{remark}
There are also maps $\ev_1\colon \Filt_\Q(\cX)\to \cX$, $\gr\colon \Filt_\Q(\cX)\to\Grad_\Q(\cX)$, etcetera, relating the stacks $\Filt_\Q(\cX)$, $\Grad_\Q(\cX)$ and $\cX$, just because the version of these maps for $\Filt$ and $\Grad$ are compatible with the colimits defining $\Filt_\Q$ and $\Grad_\Q$.
\end{remark}

\begin{proposition}\label{proposition: algebraicity of stacks of rational filtrations and rational graded points}
Let $\cX$ be an algebraic stack over an algebraic space $B$, satisfying Assumption \ref{assumption: basic assumptions}. Then $\Filt_\Q(\cX)$ and $\Grad_\Q(\cX)$ are algebraic and satisfy \Cref{assumption: basic assumptions}.
\end{proposition}
\begin{proof}
By \cite[Proposition 1.3.11]{_HalpernLeistner_Onthestructureofinstabilityinmodulitheory}, the “rising to the $n$th power” map $\Filt(\cX)\to \Filt(\cX)$ is a closed and open immersion for $n>0$. The same argument shows that the analogue map $\Grad(\cX)\to \Grad(\cX)$ is a closed and open immersion too. Thus the algebraicity result follows from Lemma \ref{lemma: filtered colimit open immersions is algebraic}. Since $\Filt_\Q(\cX)$ and $\Grad_\Q(\cX)$ are increasing unions of closed and open  substacks isomorphic to $\Filt(\cX)$ and $\Grad(\cX)$ respectively, they also satisfy \Cref{assumption: basic assumptions}.
\end{proof}

\begin{lemma}\label{lemma: filtered colimit open immersions is algebraic}
Let $I$ be a filtered poset, seen as a category, and let $D\colon I\to \cat{St}_\fppf$ be a pseudofunctor such that for all arrows $s\to t$ in $I$ the induced $D(s)\to D(t)$ is representable by open immersions. Let $\cY=\colim D$ in the 2-category $\cat{St}_\fppf$. If $D(s)$ is algebraic for all objects $s$ of $D$, then so is $\cY$.
\end{lemma}
\begin{proof}
In this proof, we consider the site of \emph{affine} schemes with the fppf topology. This does not change the 2-category $\cat{St}_\fppf$ of stacks, but it will be useful to consider only quasi-compact test schemes. 

The stack $\cY$ is the stackification of the colimit $\cY^\pre=\colim D$ in the 2-category of prestacks (meaning presheaves of groupoids). A morphism $T\to \cY^\pre$ is a pair $(s,f)$, with $s$ an object of $I$ and $f\colon T\to D(s)$ a map. A 2-morphism $(s,f)\to (s',f')$ is a pair $(t,r)$ with $t\in I$ such that there are arrows $t/s\colon s\to t$ and $t/s'\colon s'\to t$, and $r\colon D(t/s)\circ f\to D(t/s')\circ f'$ a 2-morphism. From this description, and using the facts that (1) $I$ is filtered and (2) every object in the site considered, i.e. every affine scheme, is quasi-compact, it follows that $\cY^\pre$ is already a stack, so $\cY=\cY^\pre$. Moreover, each of the maps $D(i)\to \cY$ is an open immersion. Indeed, if $(s,f)\colon T\to \cY$ is a map, with $T$ affine, then $D(i)\times_\cY T=D(i)\times_{D(s')} T$ if $i,s\leq s'$, which is open in $T$. Thus $\bigsqcup_{s\in I} D(s)\to \cY$ is a smooth representable surjection, so $\cY$ is algebraic.
\end{proof}

\begin{remark}[Functor of points of {$\Grad_\Q(\cX)$} and {$\Filt_\Q(\cX)$}]\label{remark: simple description points in Filt_Q}
From the proof of \Cref{lemma: filtered colimit open immersions is algebraic} we get a simple description of points in $\Filt_\Q(\cX)$ and $\Grad_\Q(\cX)$. Namely, if $T$ is a quasi-compact scheme, then a $T$-point of $\Filt_\Q(\cX)$ will be denoted as $\frac{1}{n}\lambda$, where $\lambda$ is a $T$-point of $\Filt(\cX)$ and $n$ is a positive integer. An isomorphism between $T$-points $\frac{1}{n}\lambda$ and $\frac{1}{n'}(\lambda')$ of $\Filt_\Q(\cX)$ is an isomorphism between $n'\lambda$ and $n(\lambda')$ in $\Filt(\cX)$. Here we are using additive notation for the “rising to the $n$th power” maps. A similar description applies to $\Grad_\Q(\cX)$.
\end{remark}

\begin{remark}
Since $(\N,\cdot, 1)$ also acts on $\Filt^n(\cX)$ and $\Grad^n(\cX)$ for all positive integers $n$, we can form the stacks $\Filt^n_\Q(\cX)$ and $\Grad^n_\Q(\cX)$ in a similar fashion. For the same reasons, they are algebraic and satisfy \Cref{assumption: basic assumptions}.
\end{remark}

\begin{example}\label{example: description filt and grad quotient stack}
Let $k$ be a field, let $X$ be a separated scheme of finite type over $k$, endowed with an action of a smooth affine algebraic group $G$ over $k$ that admits a $k$-split maximal torus $T$. We denote $W=N_G(T)/Z_G(T)$ the Weyl group. Then we have natural isomorphisms
\[\Grad^n(X/G)=\bigsqcup_{\lambda\in \Hom(\G_{m,k}^n,T)/W} X^{\lambda,0}/L(\lambda)\]
and
\[\Filt^n(X/G)=\bigsqcup_{\lambda\in \Hom(\G_{m,k}^n,T)/W} X^{\lambda,+}/P(\lambda)\]
by \cite[Theorem 1.4.8]{_HalpernLeistner_Onthestructureofinstabilityinmodulitheory}. The same description holds for $\Grad^n_\Q(X/G)$ and $\Filt^n_\Q(X/G)$ when replacing $\Hom(\G_{m,k}^n,T)/W$ by $\Q\otimes_\Z\Hom(\G_{m,k}^n,T)/W$. This can be seen by identifying the “rising to the $n$th power” maps with the “scaling by $n$” in $\Hom(\G_{m,k}^n,T)$. We will use this throughout.

Over a general base $B$, this description still holds for a quotient stack of the form $X/\GL_N$ with $X$ an algebraic space that is quasi-separated and locally finitely presented over $B$ \cite[Theorem 1.4.7]{_HalpernLeistner_Onthestructureofinstabilityinmodulitheory}.
\end{example}

The formation of $\Filt^n_\Q(\cX)$ and $\Grad^n_\Q(\cX)$ is well-behaved with respect to base change from a target algebraic space.

\begin{proposition}\label{proposition: Filt Grad and base change}
Let $\cX\to B$ and $\cX'\to B'$ satisfy \Cref{assumption: basic assumptions} and let
\[\begin{tikzcd}[ampersand replacement=\&]
    {\cX'} \& \cX \\
    {X'} \& X
    \arrow[from=2-1, to=2-2]
    \arrow[from=1-2, to=2-2]
    \arrow[from=1-1, to=2-1]
    \arrow[from=1-1, to=1-2]
    \arrow["\lrcorner"{anchor=center, pos=0.125}, draw=none, from=1-1, to=2-2]
\end{tikzcd}\]
be a cartesian square with $X$ and $X'$ algebraic spaces. Then 
\[\Grad_\Q^n(\cX')\cong \Grad_\Q^n(\cX)\times_\cX\cX'\cong \Grad_\Q^n(\cX)\times_X X'\]
and
\[\Filt_\Q^n(\cX')\cong \Filt_\Q^n(\cX)\times_{\ev_1,\cX}\cX'\cong \Filt_\Q^n(\cX)\times_X X'\]
for all $n$. The same holds for $\Filt^n$ and $\Grad^n$.
\end{proposition}
\begin{proof}
The case of $\Filt^n(\cX')$ and $\Grad^n(\cX')$ is \cite[Corollary 1.3.17]{_HalpernLeistner_Onthestructureofinstabilityinmodulitheory}. The result follows for $\Filt^n_\Q(\cX')$ and $\Grad^n_\Q(\cX')$ after covering these by copies of $\Filt^n(\cX')$ and $\Grad^n(\cX')$.
\end{proof}
\begin{proposition}\label{proposition: Filt Grad pullback closed immersion}
Let $\cX$ be an algebraic stack defined over an algebraic space $B$, satisfying \Cref{assumption: basic assumptions}, and let $\cX'\to \cX$ be a closed immersion. Then 
\[\Grad_\Q^n(\cX')\cong \Grad_\Q^n(\cX)\times_\cX\cX'\]
and
\[\Filt_\Q^n(\cX')\cong \Filt_\Q^n(\cX)\times_{\ev_1,\cX}\cX'\]
for all $n$. The same holds for $\Filt^n$ and $\Grad^n$.
\end{proposition}
\begin{proof}
As above, it is enough to see the fact for $\Filt^n$ and $\Grad^n$, which is \cite[Proposition 1.3.1]{_HalpernLeistner_Onthestructureofinstabilityinmodulitheory}.
\end{proof}

It will be useful for the sequel the fact that $\Grad$ preserves properness.

\begin{proposition}\label{proposition: f proper implies Grad f proper}
Let $f\colon \cX\to \cY$ be a representable proper finitely presented morphism of algebraic stacks over a base algebraic space $B$ satisfying \Cref{assumption: basic assumptions}. Then \[\Grad(f)\colon \Grad(\cX)\to \Grad(\cY)\] and \[\Grad_\Q(f)\colon \Grad_\Q(\cX)\to \Grad_\Q(\cY)\] are representable and proper.
\end{proposition}
\begin{proof}(Halpern-Leistner)
It is enough to prove the statement for $\Grad(f)$. Let $T$ be a scheme and $T\to \Grad(\cY)$ a map, corresponding to $B\G_{m,T}\to \cY$. Form a cartesian square

\[ \begin{tikzcd}
\cZ \arrow[r,""]\arrow[d,swap,""]\arrow[dr, phantom, "\ulcorner", very near start] & B\G_{m,T} \arrow[d,""] \\
\cX \arrow[r,""]& \cY.
\end{tikzcd}
\]
The 1-category of representable algebraic stacks over $B\G_{m,T}$ is equivalent to the category of algebraic spaces over $T$ endowed with a $\G_{m,T}$-action, and the equivalence is given by pullback along $T\to B\G_{m,T}$. Therefore $\cZ=Z/\G_{m,T}$ for a $T$-algebraic space $Z$ acted on by $\G_{m,T}$. Forming now the fibre product 
\[ \begin{tikzcd}
\cU \arrow[r,""]\arrow[d,swap,""]\arrow[dr, phantom, "\ulcorner", very near start] & T \arrow[d,""] \\
\Grad(\cX) \arrow[r,""]& \Grad(\cY)
\end{tikzcd}
\]
and given a $T$-scheme $S$, a map $S\to \cU$ over $T$ is a section of $\cZ\to B\G_{m,T}$ over $B\G_{m,S}\to B\G_{m,T}$, which is in turn a $\G_{m,T}$-equivariant map $S\to Z$. Therefore $\cU=Z^{\G_{m,T}}$, the fixed points of $Z$, as a stack over $T$. Since $Z\to T$ is finitely presented, we have by \cite[Proposition~1.4.1]{_HalpernLeistner_Onthestructureofinstabilityinmodulitheory} that the map $Z^{\G_{m,T}}\to Z$ is a closed immersion and also, by hypothesis, that $Z\to T$ is proper. Thus $\cU\to T$ is proper.
\end{proof}

We recall the definition of \emph{$\Z$-flag spaces} from \cite[Definition 1.1.15]{_HalpernLeistner_Onthestructureofinstabilityinmodulitheory} and introduce the natural counterpart of $\emph{$\Q$-flag spaces}$.

\begin{definition}[Flag spaces]\label{definition: flag spaces}
Let $\cX$ be an algebraic stack over an algebraic space $B$ satisfying \Cref{assumption: basic assumptions}, and let $x\colon T\to \cX$ be a scheme-valued point. We define the $\Z$-flag space $\Flag(x)$ to be the fibre product
\[\begin{tikzcd}[ampersand replacement=\&]
    {\Flag(x)} \& T \\
    {\Filt(\cX)} \& {\cX}
    \arrow["{\ev_1}", from=2-1, to=2-2]
    \arrow["x", from=1-2, to=2-2]
    \arrow[from=1-1, to=2-1]
    \arrow[from=1-1, to=1-2]
    \arrow["\lrcorner"{anchor=center, pos=0.125}, draw=none, from=1-1, to=2-2]
\end{tikzcd}\]
and the $\Q$-flag space $\Flag_\Q(x)$ as the fibre product
\[\begin{tikzcd}[ampersand replacement=\&]
    {\Flag_\Q(x)} \& T \\
    {\Filt_\Q(\cX)} \& {\cX.}
    \arrow["{\ev_1}", from=2-1, to=2-2]
    \arrow["x", from=1-2, to=2-2]
    \arrow[from=1-1, to=2-1]
    \arrow[from=1-1, to=1-2]
    \arrow["\lrcorner"{anchor=center, pos=0.125}, draw=none, from=1-1, to=2-2]
\end{tikzcd}\]
\end{definition}
By representability and separatedness of $\ev_1\colon \Filt(\cX)\to \cX$ and $\ev_1\Filt_\Q(\cX)\to \cX$,
the flag spaces $\Flag(x)$ and $\Flag_\Q(x)$ are separated algebraic spaces over $T$. Since $\Filt(\cX)\to \Filt_\Q(\cX)$ is an open and closed immersion, we have a natural open and closed immersion $\Flag(x)\to \Flag_\Q(x)$ of flag spaces.

In the case of a field-valued point, we can talk about a \emph{set of filtrations}.

\begin{definition}[Set of filtrations of a point]\label{definition: set of filtrations of a geometric point}
Let $\cX$ be an algebraic stack over an algebraic space $B$ satisfying \Cref{assumption: basic assumptions}. Let $k$ be a filed and let $x\colon \Spec(k)\to \cX$ be a $k$-point. The \emph{set of $\Z$-filtrations (or integral filtrations) of $x$} is defined to be
\[\Z\dash\Filt(\cX,x)\coloneqq \Flag(x)(k),\]
the set of $k$-points of the $\Z$-flag space of $x$.
Similarly, the \emph{set of $\Q$-filtrations (or rational filtrations) of $x$} is
\[\qfilt(\cX,x)\coloneqq \Flag_\Q(x)(k).\]
\end{definition}

The filtration of $x$ given by the composition of $\A^1_k/\G_{m,k}\to \Spec k$ and $x\colon \Spec k\to \cX$ is denoted $0$ and referred to as the \emph{trivial filtration}.

\begin{remark}[Filtrations of a quotient stack]\label{remark: filtrations of a quotient stack}
Let $k$ be a field and let $X$ be a separated scheme of finite type over $k$, endowed with an action $a\colon G\times X\to X$ by a linear algebraic group $G$. Form the quotient stack $\cX=X/G$ and let $x\in X(k)$ be a $k$-point. We abusively also denote by $x$ the composition
$\Spec k \xrightarrow{x} X\to \cX.$
If $\lambda\colon \G_{m,k}\to G$ is a cocharacter, we say that $\lim_{t\to 0}\lambda(t)x$ \emph{exists} in $X$ if the map \[\G_{m,k}\xrightarrow{\lambda} G\cong G\times \Spec k \xrightarrow{\id_G\times x} G\times X\xrightarrow{a} X\]
extends to a map $\overline \lambda_x\colon \A^1_k\to X$ (in which case it does so uniquely), where we regard $\G_{m,k}$ as the open subscheme $\A^1_k\setminus \{0\}$ of $\A^1_k$. If $\lim_{t\to 0}\lambda(t)x$ exists, we let $\lim_{t\to 0}\lambda(t)x$ denote the $k$-point $\overline \lambda_x(0)$ of $X$. For $n\in \Z_{>0}$, we have that $\lim_{t\to 0}\lambda(t)x$ exists if and only if $\lim_{t\to 0}\lambda^n(t)x$ exists, so it makes sense to define this notion for a rational one-parameter subgroup $\lambda \in \Gamma^\Q(G)$. The $k$-points of $X$ such that $\lim_{t\to 0}\lambda(t)x$ exists are in bijections with the $k$-points of the attractor $X^{\lambda,0}$. Thus it follows readily from \cite[Theorem 1.4.8 and Remark 1.4.9]{_HalpernLeistner_Onthestructureofinstabilityinmodulitheory} that we have an identification:
\begin{equation}\label{equation: explicit description of set of filtrations in a quotient stack}
\qfilt(\cX,x)=\{\lambda\in\Gamma^\Q(G)\st \lim_{t\to 0}\lambda(t)x \text{ exists}\}/\sim,
\end{equation}
where $\lambda\sim \lambda'$ if there is $g\in P(\lambda)(k)$ such that $\lambda^g=\lambda'$.
\end{remark}

We conclude this section with a couple of facts about maps induced on sets of filtrations.

\begin{proposition}\label{proposition: proper map induces bijection of set of filtrations}
Let $f\colon \cX\to \cY$ be a schematic proper morphism of algebraic stacks over an algebraic space $B$ satisfying \Cref{assumption: basic assumptions}. Let $k$ be a field, let $x\in \cX(k)$ and $y\in \cY(k)$ be $k$-points, and let $f(x)\to y$ be an isomorphism. Then the induced maps $\Z\dash\Filt(\cX,x)\to \Z\dash\Filt(\cY,y)$ and $\qfilt(\cX,x)\to \qfilt(\cY,y)$ of sets of filtrations are bijective.
\end{proposition} 
\begin{proof}
It is enough to deal with the case of integral filtrations.
An element of $\Z\dash\Filt(\cY,y)$ is a pair $(\lambda, \alpha)$ fitting in a commutative diagram
\[\begin{tikzcd}[column sep=large,row sep=2.25em]
    {\A^1_k/\G_{m,k}} & \cY \\
    {\Spec(k)}
    \arrow[""{name=0, anchor=center, inner sep=0}, "y"', from=2-1, to=1-2]
    \arrow["\lambda", from=1-1, to=1-2]
    \arrow["1", from=2-1, to=1-1]
    \arrow["\alpha"'{pos=0.4}, shift left=1, shorten <=2pt, Rightarrow, from=0, to=1-1]
\end{tikzcd}\]
and similarly for $\Z\dash\Filt(\cX,x)$. Fix such a $(\lambda,\alpha)\in \Z\dash\Filt(\cY,y)$. Now form the fibre product
\[ \begin{tikzcd}
X/\G_{m,k} \arrow[r,""]\arrow[d,swap,"r"]\arrow[dr, phantom, "\ulcorner", very near start] & \A^1_k/\G_{m,k} \arrow[d,"\lambda"] \\
\cX \arrow[r,"f"]& \cY.
\end{tikzcd}
\]
The base change is indeed of the form $X/\G_{m,k}$ for a scheme $X$ because $f$ is schematic, and $X/\G_{m,k}\to \A^1_k/\G_{m,k}$ is given by a $\G_{m,k}$-equivariant map $X\to \A^1_k$. There is a commutative solid square
\begin{equation}\label{diagram: lift proper 1}\begin{tikzcd}
    {(\A^1_k\setminus 0)/\G_{m,k}} & {\A^1_k/\G_{m,k}} \\
    {X/\G_{m,k}} & {\A^1_k/\G_{m,k}}
    \arrow[from=1-1, to=1-2]
    \arrow["{1_{\A^1/\G_m}}"{pos=0.4}, from=1-2, to=2-2]
    \arrow[from=2-1, to=2-2]
    \arrow["u"', from=1-1, to=2-1]
    \arrow["h"{description}, dotted, from=1-2, to=2-1]
\end{tikzcd}\end{equation}
and an isomorphism $r\circ u \to x$. An element of $\Z\dash\Filt(\cX,x)$ mapping to $(\lambda, u)$ is specified by a lift $h$ of the square \ref{diagram: lift proper 1} in the 2-categorical sense. Thus we want to prove that there is a unique such lift. There is a unique lift $g$ of
\[\begin{tikzcd}
    {\A^1_k\setminus 0} & {\A^1_k} \\
    X & {\A^1_k}
    \arrow[from=1-1, to=1-2]
    \arrow["{1_{\A^1}}"{pos=0.5}, from=1-2, to=2-2]
    \arrow[from=2-1, to=2-2]
    \arrow[from=1-1, to=2-1]
    \arrow["g"{description}, dotted, from=1-2, to=2-1]
\end{tikzcd}\]
by properness, see \cite[\href{https://stacks.math.columbia.edu/tag/0BX7}{Tag 0BX7}]{stacks-project}. We just need to prove that $g$ is $\G_{m,k}$-equivariant. Equivariance amounts to the commutativity of

\[ \begin{tikzcd}
\G_{m,k}\times \A^1_k \arrow[r,""]\arrow[d,swap,"1_{\G_m}\times g"] & \A^1_k \arrow[d,"g"] \\
\G_{m,k}\times X \arrow[r,""]& X.
\end{tikzcd}
\]
Both compositions agree when restricted to $\G_{m,k}\times (\A^1_k\setminus 0)$, which is schematically dense in $\G_{m,k}\times \A^1_k$. Thus, by separatedness of $X\to \A^1_k$, the square commutes.
\end{proof}

\begin{proposition}\label{proposition: representable and separated induced injection on set of filtrations}
Suppose $f\colon \cX\to \cY$ is a representable and separated morphism of algebraic stacks over an algebraic space $B$ satisfying Assumption \ref{assumption: basic assumptions}. Let $k$ be a field, let $x\in \cX(k)$ and $y\in \cY(k)$ be $k$-points, and let $f(x)\to y$ be an isomorphism. Then the induced maps $\Z\dash\Filt(\cX,x)\to \Z\dash\Filt(\cY,y)$ and $\qfilt(\cX,x)\to \qfilt(\cY,y)$ of sets of filtrations are injective
\end{proposition}
\begin{proof}
It is enough to prove the claim for integral filtrations. If $\lambda_1, \lambda_2$ are two filtrations of $x$ that give the same filtration of $y$, we can form a commutative diagram 
\[ \begin{tikzcd}
\Spec k \arrow[r,"x"]\arrow[d,swap,"1"] & \cX \arrow[d,"\Delta_f"] \\
\A^1_k/\G_{m,k} \arrow[ur, dashed]\arrow[r,"{(\lambda_1,\lambda_2)}"]& \cX\times_\cY \cX.
\end{tikzcd}
\]
Since $\Delta_f$ is a closed immersion and $1\colon \Spec k\to \A^1_k/\G_{m,k}$ is a schematically dense open immersion, there is a unique dashed arrow filling the diagram, which gives the isomorphism between $\lambda_1$ and $\lambda_2$.
\end{proof}

\subsection{Normed stacks}
We now recall from \cite{_HalpernLeistner_Onthestructureofinstabilityinmodulitheory} the notion of a norm on graded points of a stack. If $G$ is an algebraic group, norms on graded points of $BG$ are in bijection with norms on cocharacters of $G$ (\Cref{proposition: norm on BG is norm on cocharacters of G}). 

\begin{definition}[Nondegenerate graded point]\label{definition: nondegenerate graded point}
Let $\cX$ be a quasi-separated algebraic stack, let $k$ be a field and let $x\colon B\G^n_{m,k}\to \cX$ be a $\Z^n$-graded point. We say that the graded point $x$ is \emph{nondegenerate} if $\ker(\G_{m,k}^n\to \Aut(x\vert_{\Spec k}))$ is finite.

Suppose that $\cX$ is defined over some algebraic space $B$ and it satisfies Assumption \ref{assumption: basic assumptions}. We say that a connected component $\cZ$ of $\Grad^n(\cX)$ is \emph{nondegenerate} if there is a field $k$ and a point $x\in \cZ(k)$ with $x$ nondegenerate.
\end{definition}

\begin{remark}
If $\cZ$ is a nondegenerate component, then by \cite[Proposition 1.3.9]{_HalpernLeistner_Onthestructureofinstabilityinmodulitheory} we have that for \emph{all} fields $k$ and points $x\in \cZ(k)$, the point $x$ is nondegenerate.
\end{remark}

We recall the notion of norm on graded points of a stack from \cite[Definition 4.1.12]{_HalpernLeistner_Onthestructureofinstabilityinmodulitheory}.

\begin{definition}[Norm on graded points]\label{definition: norm on graded points}
Let $\cX$ be an algebraic stack over an algebraic space $B$, satisfying \Cref{assumption: basic assumptions}. A \emph{(rational quadratic) norm} $q$ on graded points of $\cX$ (or simply a \emph{norm} on $\cX$) is a locally constant function
\[q\colon \abs{\Grad(\cX)}\to \Q_{\geq 0}\]
such that for every field $k$ and every nondegenerate $\Z^n$-graded point $x\colon B\G^n_{m,k}\to \cX$, the induced map $q_x\colon \Gamma^\Z(\G^n_{m,k})\to \Q$ is the quadratic form of a rational inner product on the finite free $\Z$-module $\Gamma^\Z(\G^n_{m,k})$ of cocharacters of $\G^n_{m,k}$.

A \emph{normed algebraic stack} is an algebraic stack endowed with a norm on graded points. 
\end{definition}

Let us clarify what the map $q_x$ is. If $\lambda\colon \G_{m,k}\to \G_{m,k}^n$ is a cocharacter, the composition $B\G_{m,k}\xrightarrow{B\lambda} B\G_{m,k}^n\xrightarrow{x} \cX$ defines a point $p\in \abs{\Grad(\cX)}$, and we let $q_x(\lambda)=q(p)$.

\begin{remark}\label{remark: norm extends to rational graded points}
A norm on graded points $q\colon \abs{\Grad(\cX)}\to \Q$ extends canonically to a map \[q\colon \abs{\Grad_\Q(\cX)}\to \Q\] by setting $q(\frac{1}{n}\lambda)=\frac{1}{n^2}q(\lambda)$, for a rational graded point $\frac{1}{n}\lambda$.
\end{remark}

\begin{remark}[Notation for norms]\label{remark: notation for norms}
If $\cX$ is endowed with a norm on graded points $q$ and $\lambda\in \Grad_\Q(\cX)(k)$ is a graded point, with $k$ a field, then we will often denote $\norm{\lambda}\coloneqq \sqrt{q(\lambda)}$.
\end{remark}

In some circumstances we can pull back a norm under a morphism.

\begin{propositionanddefinition}[Pulling back norms]\label{proposition: pulling back norms}
Let $\cX$ and $\cY$ be algebraic stacks over an algebraic space $B$, satisfying \Cref{assumption: basic assumptions}. Let $f\colon \cX\to \cY$ be a morphism such that the relative inertia $\cI_f\to \cX$ has proper fibres (for example if $f$ is representable or separated). Let $q$ be a norm on $\cY$ and denote $f^*q$ the composition $\begin{tikzcd}[column sep=small]
\abs{\Grad(\cX)} \arrow[r,""] & \abs{\Grad(\cY)} \arrow[r,"q"] & \Q
\end{tikzcd}$. Then $f^*q$ is a norm on $\cX$, called the \emph{pulled back norm}. 

If $\cX$ is endowed with a norm $q'$, we say that the morphism $f$ is \emph{norm-preserving} if $f^*q=q'$.
\end{propositionanddefinition}
\begin{proof}
We need to see that if $u\colon B\G_{m,k}^n\to \cX$ is nondegenerate, then so is $f\circ u$. Let $x=u\vert_{\Spec k}\in \cX(k)$ be the point that $u$ grades. We have induced algebraic group homomorphisms
\[\begin{tikzcd}
    {\G_{m,k}^n} \\
    {\Aut(x)} & {\Aut(f(x)).}
    \arrow["l", from=1-1, to=2-2]
    \arrow["s"', from=1-1, to=2-1]
    \arrow["r"', from=2-1, to=2-2]
\end{tikzcd}\]
The kernel of $r$ is proper over $\Spec(k)$, since it is the fibre over $x$ of the relative inertia morphism, and $\ker s$ is finite by hypothesis. We have a sequence 
\[ \begin{tikzcd}
\ker l \arrow[r,"a"] & \ker l/\ker s \arrow[r,"b"] & \ker r
\end{tikzcd}
\]
where $a$ is finite and $b$ is a closed immersion. Thus $\ker l$ is proper over $k$, and therefore finite, because $\G^n_{m,k}$ is affine. This proves that $f\circ u$ is nondegenerate.
\end{proof}

The notion of norm on graded points of an algebraic stack is a generalisation of the classical notion of norm on cocharacters of a group.

\begin{definition}[Norm on cocharacters of a group]\label{definition: norm on cocharacters of a group}
Let $G$ be a smooth affine algebraic group, over a field $k$, that has a $k$-split maximal torus. A \emph{(rational quadratic) norm on cocharacters of $G$} is a map $q\colon \Gamma^\Z(G)\to \Q$ that is invariant under the action of $G(k)$ on $\Gamma^\Z(G)$ by conjugation and such that for every $k$-split torus $T$ of $G$, the restriction of $q$ to $\Gamma^\Z(T)$ is the quadratic form of a rational inner product on $\Gamma^\Z(T)$.
\end{definition}

Now fix a smooth connected linear algebraic group $G$ over a field $k$. Any two maximal $k$-split tori are conjugate by an element of $G(k)$ \cite[Theorem~C.2.3]{_Conrad_PseudoreductiveGroups}, and there exists one such torus for dimension reasons. The Weyl group $W(G,T)=N_G(T)/Z_G(T)$, for a $k$-split torus $T$, is a finite constant group scheme \cite[Proposition~C.2.10]{_Conrad_PseudoreductiveGroups}, so we identify it with a finite abstract group. By conjugacy of maximal $k$-split tori and \cite[Lemma~2.8]{_Mumford_GeometricInvariantTheory}, it follows that $\Gamma^\Z(T)/W(G,T)=\Gamma^\Z(G)/G(k)$ if $T$ is a maximal $k$-split torus of $G$. Thus we deduce:

\begin{proposition}
If $T$ is a $k$-split maximal torus of $G$, then the data of a norm on cocharacters of $G$ is equivalent to a rational quadratic inner product on $\Gamma^\Z(T)$, invariant under the action of $W(G,T)$. \hfill \qedsymbol{}
\end{proposition}

The link with the concept of a norm on graded points is given in the following well-known proposition.

\begin{proposition}\label{proposition: norm on BG is norm on cocharacters of G}
Suppose $G$ has a $k$-split maximal torus $T$. Then norms on $BG$ are in natural bijection with norms on cocharacters of $G$.
\end{proposition}
\begin{proof}
Let $W=W(G,T)$ be the Weyl group. Then, by \cite[Theorem~1.4.8]{_HalpernLeistner_Onthestructureofinstabilityinmodulitheory}, we can explicitly describe the stack of graded points as
\[\Grad(BG)=\bigsqcup_{\lambda \in\Gamma^\Z(T)/W} BL(\lambda),\]
where $L(\lambda)$ is the centraliser of a choice of representative of $\lambda$. Thus $\abs{\Grad(BG)}=\Gamma^\Z(T)/W$, and the identification is compatible in the sense that if $\lambda\colon \G_m\to G$ is a cocharacter, the point in $\abs{\Grad(BG)}$ defined by $B\lambda\colon B\G_m\to BG$ is the class of $\lambda$ in $\Gamma^\Z(T)/W=\Gamma^\Z(G)/G(k)$. To conclude, just note that if $T\to T'$ is a map of $k$-split tori with finite kernel, then $\Gamma^\Z(T)\to \Gamma^\Z(T')$ is injective, so $\Gamma^\Z(T)$ inherits an inner product if $\Gamma^\Z(T')$ has one.
\end{proof}

Norms on cocharacters are a source of norms on quotient stacks.

\begin{proposition}[Norms on quotient stacks]
Suppose $G$ has a $k$-split maximal torus and acts on an quasi-separated algebraic space $X$ of finite type over $k$. If $G$ is endowed with a norm on cocharacters $q$, then $X/G$ is naturally endowed with a norm on graded points.
\end{proposition}
\begin{proof}
Since $p\colon X/G\to BG$ is representable, the pullback $p^*q$ is a norm on $X/G$ by \Cref{proposition: pulling back norms}.
\end{proof}

\begin{definition}[Norm on {$\Grad_\Q(\cX)$}]\label{definition: norm on GradX from norm on X}
Let $\cX$ be an algebraic stack over an algebraic space $B$ satisfying Assumption \ref{assumption: basic assumptions} and endowed with a norm on graded points $q$. We denote $\Grad(q)$ (resp. $\Grad_\Q(q)$) the norm on $\Grad(\cX)$ (resp. $\Grad_\Q(\cX)$) given as the pullback of $q$ along the forgetful morphism $u \colon \Grad(\cX)\to \cX$ (resp. $\Grad_\Q(\cX)\to \cX$), which is representable.
\end{definition}
It follows from Proposition \ref{proposition: pulling back norms} that $\Grad(q)$ is a norm on $\Grad(\cX)$. 
If $\cX$ is a normed stack, we will always regard $\Grad(\cX)$ (resp. $\Grad_\Q(\cX)$) as a normed stack, endowed with the norm $\Grad(q)$ (resp. $\Grad_\Q(q)$).

\subsection{Linear forms on stacks}

We now recall the notion of linear form on graded points of a stack from \cite{_HalpernLeistner_Onthestructureofinstabilityinmodulitheory} and show how to get linear forms from line bundles.

\begin{definition}[Linear form]\label{definition: linear form on stack}
Let $\cX$ be an algebraic stack over an algebraic space $B$, satisfying \Cref{assumption: basic assumptions}. A \emph{rational linear form $\ell$ on graded points of $\cX$} (or simply a \emph{linear form on $\cX$}) is a locally constant function
\[\ell\colon \abs{\Grad(\cX)}\to \Q\]
such that, for every field $k$ and every $\Z^n$-graded point $B\G_{m,k}^n\to\cX$, the induced map $\Gamma^\Z(\G_{m,k}^n)\to \Q$ on cocharacters of the torus is $\Z$-linear. 
\end{definition}

If $\lambda\colon B\G_{m,k}\to \cX$ is a graded point, we denote by either $\langle \lambda,\ell \rangle$ or $\ell(\lambda)$ the value of $\ell$ at the point of $\Grad(\cX)$ defined by $\lambda$.

\begin{remark}\label{remark: linear form extends to rational graded points}
A linear form $\ell\colon \abs{\Grad(\cX)}\to \Q$ extends canonically to a map $\ell\colon \abs{\Grad_\Q(\cX)}\to \Q$ by setting $\langle \frac{1}{n}\lambda,\ell\rangle=\frac{1}{n}\langle \lambda, \cL\rangle$, for a rational graded point $\frac{1}{n}\lambda$.
\end{remark}

Line bundles are an important source of linear forms. The following definition essentially comes from \cite{_Heinloth_HilbertMumfordstabilityonalgebraicstacksandapplicationstoGbundlesoncurves} and \cite{_HalpernLeistner_Onthestructureofinstabilityinmodulitheory}.

\begin{definition}[Linear form associated to a line bundle]\label{definition: linear form from line bundle}
Let $\cX$ be an algebraic stack over an algebraic pace $B$, satisfying \Cref{assumption: basic assumptions}, and let $\cL$ be a line bundle on $\cX$. Define a map
\[\langle -,\cL\rangle\colon \abs{\Grad(\cX)}\to \Q\]
as follows. If $g\in \abs{\Grad(\cX)}$ is represented by a field $k$ and a map $\lambda\colon B\G_{m,k}\to \cX$, let $\langle g, \cL\rangle\coloneqq \langle \lambda, \cL\rangle\coloneqq -\wt(\lambda^*\cL)$, the opposite of the weight of the one-dimensional representation $\lambda^*\cL$ of $\G_{m,k}$.
\end{definition}
\begin{remark}[Sign conventions]\label{remark: sign convention on weights}
Our sign convention for weights is as follows. Let $k$ be a field. We denote by $\cO_{B\G_{m,k}}(n)$, for $n\in \Z$, the representation of $\G_{m,k}$ with underlying vector space $k$ and such that $t*1=t^{n}$ for $t\in k^\times$. By definition, $\wt \cO_{B\G_{m,k}}(n)=n$.

The total space of $\cO_{B\G_{m,k}}(n)$ is $\bA(\cO_{B\G_{m,k}}(n))=\Spec_{B\G_{m,k}}(\Sym (\cO_{B\G_{m,k}})^\vee)$. If we let $\G_{m,k}$ act on $\A^1_k$ by the formula $t\ast s=t^ns$, for any $k$-algebra $R$, $t\in \G_{m,k}(R)=R^\times$ and $s\in \A^1_k(R)=R$, then $\bA(\cO_{B\G_{m,k}}(n))=\A^1_k/\G_{m,k}$ for this action. If we write $\A^1_k=\Spec k[x]$, the standard coordinate $x$ has weight $-n$.

Let $G$ be a linear algebraic group over $k$ acting linearly on a finite dimensional vector space $V$, and let $p\colon\P(V)/G\to B\G_m$ be the map $\P(V)/G=(V\setminus \{0\})/G\times \G_m\to B(G\times \G_m)\to B\G_m$, where $\G_m$ acts by scaling on $V$. Then $p^*(\cO_{B\G_{m,k}}(1))=\cO_{\P(V)/G}(1)$ is the standard ample line bundle on $\P(V)/G$. Let $\lambda\colon \G_{m,k}\to G$ be a one-parameter subgroup, and let $x\in \P(V)(k)$ be a point fixed by $\lambda$, defining a graded point $\lambda_x\colon B\G_{m,k}\to \P(V)/G$. The point $x$ corresponds to a one-dimensional subspace $L\subset V$ invariant by $\lambda$. Thus there is $n\in \Z$ such that $\lambda(t)v=t^nv$ for all $k$-algebras $R$, $v\in R\otimes_k L$ and $t\in R^\times$. In other words, $L$ is regarded as a $\G_{m,k}$-representation and we let $n=\wt L$. We can identify $L=\lambda_x^*\left(\cO_{\P(V)/G}(-1)\right)$. Thus our sign conventions are such that $\langle \lambda_x,\cO_{\P(V)/G}(1)\rangle=-\wt (L^\vee)=n$.
\end{remark}

\begin{proposition}
Let $\cX$ and $\cL$ be as in \Cref{definition: linear form from line bundle}. Then $\langle -,\cL\rangle$ is well-defined and a linear form on $\cX$.
\end{proposition}
\begin{proof}
If $k'\vert k$ is a field extension, inducing a map $g\colon B\G_{m,k'}\to B\G_{m,k}$, and $U$ is a line bundle on $B\G_{m,k'}$, then $\wt U=\wt (g^*U)$, so  $\langle -,\cL\rangle$ is well-defined. 

To see that $\langle -,\cL\rangle$ is locally constant on $\abs{\Grad(\cX)}$, it is enough to prove that for any map $f\colon \Spec A \to \Grad(\cX)$ with $A$ any commutative ring, the composition $\abs{\Spec A}\to \abs{\Grad(\cX)}\to \Q$ of $\abs{f}$ and $\langle -,\cL\rangle$ is locally constant. Let $h\colon B\G_{m,A}\to \cX$ be the map corresponding to $f$. We may assume that $h^*\cL$ is trivial when restricted to $\Spec A$, so $h^*\cL$ is an $A$-module direct sum decomposition $A=\bigoplus_{n\in \Z} A_n$. Let $1_n$ be the degree $n$ part of $1\in A$. Each $A_n$ is generated by $1_n$ as an $A$-module. On the nonvanishing locus $D(1_n)$ of $1_n$, we have that the restriction $A_n\vert_{D(1_n)}=A\vert_{D(1_n)}$ and $A_m\vert_{D(1_n)}=0$ if $m\neq n$. Thus $\Spec A= \bigsqcup_{n\in \Z} D(1_n)$, and the composition $\abs{D(1_n)}\to \abs{\Grad(\cX)}\xrightarrow{\langle-,\cL\rangle} \Q$ is constant with value $-n$.

Now let $\alpha\colon B\G_{m,k}^n\to \cX$ be a $\Z^n$-graded point. The pullback $\alpha^*\cL$ corresponds to a character $\chi\in \Gamma_\Z(\G_{m,k}^n)$, and the map $\Gamma^\Z(\G_{m,k}^n)\to \Q$ that $\langle -,\cL\rangle$ induces is just the pairing $\lambda\mapsto -\langle \lambda,\chi\rangle$. It is thus linear.
\end{proof}

\begin{remark}
\Cref{definition: linear form from line bundle} naturally extends to \emph{rational} line bundles $\cL\in \Pic(\cX)\otimes_\Z \Q$.
\end{remark}

Now suppose that $\cX$ is an algebraic stack over an algebraic space $B$ satisfying Assumption \ref{assumption: basic assumptions}, and endowed with a norm on graded points $q$. 

\begin{definition}[Canonical linear form of a norm]\label{definition: canonical linear form on graded points}
The \emph{canonical linear form $\ell_q$ on $\Grad_\Q(\cX)$} induced by $q$ is the linear form on graded points of $\Grad_\Q(\cX)$ defined as follows. A graded point $\mu\in \Grad(\Grad_\Q(\cX))(k)$ lying over a point $\lambda/n\in \Grad_\Q(\cX)(k)$, with $\lambda\in \Grad(\cX)(k)$ and $n\in \Z_{>0}$, gives a map $\alpha=(\mu,\lambda)\colon B\G_{m,k}^2\to \cX$ and $q$ gives an inner product on the set $\Gamma^\Z(\G_{m,k}^2)$ of cocharacters of $\G_{m,k}^2$. Let $e_1,e_2$ be the standard basis of $\Gamma^\Z(\G_{m,k}^2)$. We define
\[\langle \lambda, \ell_q\rangle=\frac{1}{n}(e_1,e_2)_q.\]
\end{definition}

\begin{remark}
Note that $\ell_q$ determines $q$, since for a graded point $\lambda\colon B\G_{m,k}\to \cX$, if $o\colon B\G_{m,k}^2\to B\G_{m,k}$ is induced by $\G_{m,k}^2\to \G_{m,k}\colon (t,t')\to tt'$, then $q(\lambda)=\ell_q(\lambda\circ o)$.
\end{remark}

\begin{definition}[Algebraic norm]\label{definition: algebraic norm}
We say that the norm $q$ on $\cX$ is \emph{algebraic} if the canonical linear form $\ell_q$ on $\Grad_\Q(\cX)$ is induced, on each connected component $\cZ$ of $\Grad_\Q(\cX)$ by a rational line bundle on $\cZ$.
\end{definition}

\begin{remark}
If $q$ is an algebraic norm on $\cX$ and $f\colon \cY\to \cX$ is a morphism with proper relative automorphism groups, then the pulled back norm $f^*q$ (\Cref{proposition: pulling back norms}) is algebraic. Indeed, if $\ell_q=\langle -,\cM\rangle$ for a rational line bundle $\cM$ on $\Grad_\Q(\cX)$, then $\ell_{f^*q}=\langle -, \Grad_\Q(f)^*\cM\rangle$. We will tacitly use this fact in the sequel.
\end{remark}

\subsection{\texorpdfstring{$\Theta$}{Theta}-stratifications}

We now discuss the notion of \emph{$\Theta$-stratification} for algebraic stacks, which is a generalisation due to Halpern-Leistner of the Hesselink-Kempf-Kirwan-Ness stratification in GIT \cite[Chapter 12]{_Kirwan_Cohomologyofquotientsinsymplecticandalgebraicgeometry} and of the stratification by Harder-Narasimhan type for vector bundles on a curve \cite{_Atiyah_TheYangMillsEquationsoverRiemannSurfaces}. We follow \cite{_HalpernLeistner_Onthestructureofinstabilityinmodulitheory} with some modifications. We use the stack $\Filt_\Q(\cX)$ of rational filtrations instead of $\Filt(\cX)$ and we work with a reformulation of the original definition of $\Theta$-stratification that is closer to our definition of \emph{sequential stratification} (\Cref{definition: sequential stratification}). We will focus on $\Theta$-stratifications induced by a pair $(\ell,q)$, where $q$ is a norm and $\ell$ is a linear form on graded points of a stack $\cX$.

We fix an algebraic stack $\cX$ over an algebraic space $B$ satisfying Assumption \ref{assumption: basic assumptions}. The following is a variant of \cite[Definition 2.1.2]{_HalpernLeistner_Onthestructureofinstabilityinmodulitheory}.

\begin{definition}[$\Theta$-stratification]\label{definition: Theta-stratification my version}
Let $\Gamma$ be a partially ordered set. A \emph{(weak) $\Theta$-stratification} of $\cX$ indexed by $\Gamma$ is a family $(\cS_c)_{c\in \Gamma}$ of open substacks of $\Filt_\Q(\cX)$ satisfying:
\begin{enumerate}
\item For every $c\in \Gamma$, the composition $\cS_c\to \Filt_\Q(\cX)\xrightarrow{\ev_1} \cX$ is a locally closed immersion (resp. locally finite radicial\footnote{We say that a morphism $\cY\to \cZ$ is \emph{locally finite radicial} if it factors as $\cY\to \cU\to \cZ$ with $\cY\to \cU$ finite and radicial and $\cU\to \cZ$ an open immersion.}).
\item The $\ev_1(\abs{\cS_c})$ are pairwise disjoint and cover $\abs{\cX}$.
\item For every $c\in \Gamma$, the stratum $\cS_c$ is the preimage along $\gr\colon \Filt_\Q(\cX)\to \Grad_\Q(\cX)$ of an open substack $\cZ_c$ of $\Grad_\Q(\cX)$, called the \emph{centre} of $\cS_c$.
\item For every $c\in \Gamma$, the set $\abs{\cX_{\leq c}}\coloneqq\bigcup_{c'\leq c}\ev_1(\abs{\cS_{c'}})$ is open in $\abs{\cX}$, and it thus defines an open substack $\cX_{\leq c}$ of $\cX$.
\end{enumerate}
\end{definition}
\begin{remark}
If $\sigma\colon \Grad_\Q(\cX)\to \Filt_\Q(\cX)$ denotes the “split filtration” map, then for all $c\in\Gamma$ the centre of $\cS_c$ is $\cZ_c=\sigma^{-1}(\cS_c)$ and thus it is uniquely determined.
\end{remark}
The strata $\cS_c$ in a $\Theta$-stratification $(\cS_c)_{c\in \Gamma}$ are locally closed $\Theta$-strata in the following sense.
\begin{definition}\label{definition: locally closed Theta-stratum}
A \emph{locally closed $\Theta$-stratum} of $\cX$ is an open substack $\cS$ of $\Filt_\Q(\cX)$ that is is the preimage along $\gr$ of an open substack $\cZ$ of $\Grad_\Q(\cX)$ (its \emph{centre}) and such that the composition $\cS\to \Filt_\Q(\cX)\xrightarrow{\ev_1} \cX$ is a locally closed immersion.
\end{definition}

\begin{proposition}\label{proposition: Halpern-Leistners definition of theta-stratification is equivalent to mine}
Let $\Gamma$ be a totally ordered set. Let $(\cS_c)_{c\in \Gamma}$ be a (weak) $\Theta$-stratification of $\cX$. If each stratum $\cS_c$ is contained in the closed and open substack $\Filt(\cX)\subset \Filt_\Q(\cX)$ of integral filtrations, then $(\cS_c)_{c\in \Gamma}$ and $(\cX_{\leq c})_{c\in \Gamma}$ define a (weak) $\Theta$-stratification in the sense of \cite[Definition 2.1.2]{_HalpernLeistner_Onthestructureofinstabilityinmodulitheory}. Conversely, any (weak) $\Theta$-stratification in the sense of \cite[Definition 2.1.2]{_HalpernLeistner_Onthestructureofinstabilityinmodulitheory} defines a (weak) $\Theta$-stratification.
\end{proposition}
\begin{proof}
Let $(\cS_c)_{c\in \Gamma}$ be a (weak) $\Theta$-stratification of $\cX$ in the sense of \Cref{definition: Theta-stratification my version}. It is enough to show that each $\cS_c$ is a (weak) $\Theta$-stratum of $\cX_{\leq c}$ \cite[Definition 2.1.1]{_HalpernLeistner_Onthestructureofinstabilityinmodulitheory}, that is
\begin{enumerate}
\item Each $\cS_c\to \Filt_\Q(\cX)$ factors through $\Filt_\Q(\cX_{\leq c})$ and $\cS_c\to \Filt_\Q(\cX_{\leq c})$ is an open and closed immersion.

\item The composition $\cS_c\to \Filt_\Q(\cX_{\leq c})\to \cX_{\leq c}$ is a closed immersion (resp. finite and radicial).
\end{enumerate}
Then the conditions in \cite[Definition 2.1.2]{_HalpernLeistner_Onthestructureofinstabilityinmodulitheory} are trivially satisfied by construction. Note that Halpern-Leistner also demands the $\cS_c$ to be integral and $\Gamma$ to be a total order.

We have a diagram
\[\begin{tikzcd}
    {\cZ_c} & {\Grad_\Q(\cX)} & \cX \\
    {\cS_c} & {\Filt_\Q(\cX)} & \cX
    \arrow[Rightarrow, no head, from=1-3, to=2-3]
    \arrow["\sigma", from=1-2, to=2-2]
    \arrow[from=1-1, to=2-1]
    \arrow[from=2-1, to=2-2]
    \arrow[from=1-1, to=1-2]
    \arrow[from=1-2, to=1-3]
    \arrow["{\ev_1}", from=2-2, to=2-3]
    \arrow["\ulcorner"{anchor=center, pos=0.125}, draw=none, from=1-1, to=2-2]
\end{tikzcd}\]
and $\cS_c\to \cX$ factors through $\cX_{\leq c}$. Therefore $\cZ_c\to \cX$ also factors through $\cX_{\leq c}$ and thus $\cZ_c\subset \Grad_\Q(\cX_{\leq c})$ since the formation of $\Grad_\Q$ is compatible with immersions. Since the natural square
\[\begin{tikzcd}[ampersand replacement=\&]
    {\Filt_\Q(\cX_{\leq c})} \& {\Filt_\Q(\cX)} \\
    {\Grad_\Q(\cX_{\leq c})} \& {\Grad_\Q(\cX)}
    \arrow["\gr", from=1-2, to=2-2]
    \arrow["\gr", from=1-1, to=2-1]
    \arrow[from=2-1, to=2-2]
    \arrow[from=1-1, to=1-2]
    \arrow["\ulcorner"{anchor=center, pos=0.125}, draw=none, from=1-1, to=2-2]
\end{tikzcd}\]
is cartesian \cite[Proposition 1.3.1 (3)]{_HalpernLeistner_Onthestructureofinstabilityinmodulitheory} and $\cS_c=\gr^{-1}(\cZ_c)$, we have $\cS_c\subset \Filt_\Q(\cX_{\leq c})$.

Now the composition $\cS_c\to \Filt_\Q(\cX_{\leq c})\to \cX_{\leq c}$ is a locally closed immersion (resp. locally finite radicial) and has closed image. It is thus a closed immersion (resp. finite and radicial). Since $\Filt_\Q(\cX_{\leq c})\to \cX_{\leq c}$ is representable and separated, the map $\cS_c\to \Filt_\Q(\cX_{\leq c})$ is also a closed immersion (resp. finite and radicial). Since $\cS_c\to \Filt_\Q(\cX_{\leq c})$ is also an open immersion, it has to be an open and closed immersion.

For the converse, if $(\cS_c)_{c\in \Gamma}$ and $(\cX_{\leq c})_{c\in \Gamma}$ define a $\Theta$-stratification in the sense of \cite[Definition 2.1.1]{_HalpernLeistner_Onthestructureofinstabilityinmodulitheory}, just note that $\cS_c\to \Filt_\Q(\cX_{\leq c})$ being an open and closed immersion implies that $\cS_c=\gr^{-1}(\cZ_c)$ for some $\cZ_c\subset \Grad_\Q(\cX_{\leq c})$ open and closed, because $\Filt_\Q(\cX_{\leq c})$ and $\Grad_\Q(\cX_{\leq c})$ have the same connected components \cite[Lemma 1.3.8]{_HalpernLeistner_Onthestructureofinstabilityinmodulitheory}. The other conditions of \Cref{definition: Theta-stratification my version} are readily seen to be satisfied.
\end{proof}

\begin{remark}\label{remark: integrality of strata is not very important}
If each $\cS_c$ intersects only a finite number of connected components of $\Filt_\Q(\cX)$, then it becomes integral after scaling up, that is, after replacing $\cS_c$ by its image under the “rising to the $n$th power” map $\Filt_\Q(\cX)\to \Filt_\Q(\cX)$ for big enough $n$. After suitably subdividing each $\cS_c$, the condition is always satisfied.
\end{remark}

A $\Theta$-stratification of $\cX$ assigns, for each geometric point $x$ of $\cX$, a canonical rational filtration of $x$:

\begin{definition}[{HN filtrations, \cite[Lemma 2.1.4]{_HalpernLeistner_Onthestructureofinstabilityinmodulitheory}}]\label{definition: HN filtration}
Consider a $\Theta$-stratification $(\cS_c)_{c\in \Gamma}$ of $\cX$, and let $x\colon \Spec(k)\to \cX$ be a field-valued point. There is a unique $c\in \Gamma$ such that the image of $\cS_c\to \cX$ contains the point of $\abs{\cX}$ defined by $x$; and there is a unique, up to unique isomorphism, lift of $x$ to $\cS_c$, which gives a filtration $\lambda\in \qfilt(\cX,x)$ called the \emph{Harder-Narasimhan filtration} (or the \emph{HN filtration}) of $x$.
\end{definition}

\begin{remark}
In the case of a weak $\Theta$-stratification, the HN filtration of a $k$-point is defined over a finite purely inseparable extension of $k$. We will only use the case of $\Theta$-stratifications.
\end{remark}

The following is a reformulation of \cite[Definition 2.3.1]{_HalpernLeistner_Onthestructureofinstabilityinmodulitheory}, convenient for our purposes.

\begin{propositionanddefinition}[Induced {$\Theta$}-stratifications]\label{propositionanddefinition: induced Theta-stratifications}
Let $\cX'$ be an algebraic stack over an algebraic space $B$, satisfying \Cref{assumption: basic assumptions}, and let $h\colon \cX'\to \cX$ be either a closed immersion or a base change of a map between algebraic spaces like in \Cref{proposition: Filt Grad and base change}. Let $(\cS_c)_{c\in \Gamma}$ be a (weak) $\Theta$-stratification of $\cX$. For each $c\in \Gamma$ let $h^*\cS_c$ be the pullback
\[\begin{tikzcd}[ampersand replacement=\&]
    {h^*\cS_c} \& {\cS_c} \\
    {\Filt_\Q(\cX')} \& {\Filt_\Q(\cX).}
    \arrow[from=1-1, to=1-2]
    \arrow[""{name=0, anchor=center, inner sep=0}, from=2-1, to=2-2]
    \arrow[from=1-2, to=2-2]
    \arrow[from=1-1, to=2-1]
    \arrow["\lrcorner"{anchor=center, pos=0.125}, draw=none, from=1-1, to=0]
\end{tikzcd}\]
Then the family $(h^*\cS_c)_{c\in \Gamma}$ is a (weak) $\Theta$-stratification of $\cX'$ called the \emph{$\Theta$-stratification induced by} $h$ and $(\cS_c)_{c\in \Gamma}$.
\end{propositionanddefinition}
\begin{proof}
This is the content of \cite[Lemmas 2.3.2 and 2.3.3]{_HalpernLeistner_Onthestructureofinstabilityinmodulitheory} in the case of integral filtration, from which the case of rational filtrations follows easily.
\end{proof}

Now we fix a rational quadratic norm $q$ and a rational linear form $\ell$ on graded points of $\cX$. We regard $q$ and $\ell$ as maps $\abs{\Grad_\Q(\cX)}\to \Q$ by \Cref{remark: norm extends to rational graded points,remark: linear form extends to rational graded points}. Moreover, $q$ and $\ell$ induce functions on $\abs{\Filt_\Q(\cX)}$ by precomposing along $\abs{\gr}\colon \abs{\Filt_\Q(\cX)}\to \abs{\Grad_\Q(\cX)}$ that we will still denote $q$ and $\ell$.

The norm $q$ and the linear form $\ell$ give rise to two other interesting functions.

\begin{definition}[Associated numerical invariant {\cite[4.1.1 and 4.1.14]{_HalpernLeistner_Onthestructureofinstabilityinmodulitheory}}]
We define the \emph{numerical invariant} $\mu$ associated to $q$ and $\ell$ to be the function $\mu\colon \abs{\Grad_\Q(\cX)}\to \R_{\geq 0}$ such that
\begin{enumerate}
\item on the open and closed substack $\sigma\colon \cX\to \Grad_\Q(\cX)$ defined by the “trivial grading” map, $\mu$ takes the value $0$; and
\item on $\abs{\Grad_\Q(\cX)}\setminus \abs{\cX}$, we set $\mu=\dfrac{\ell}{\sqrt{q}}$.
\end{enumerate}
\end{definition}

We extend $\mu$ to a function on $\abs{\Filt_\Q(\cX)}$ by taking the composition 

\[ \begin{tikzcd}
\abs{\Filt_\Q(\cX)} \arrow[r,"\abs{\gr}"] & \abs{\Grad_\Q(\cX)} \arrow[r,"\mu"] & \R_{\geq 0},
\end{tikzcd}
\]
which we will still denote $\mu$.

\begin{definition}[Stability function {\cite[4.1.1]{_HalpernLeistner_Onthestructureofinstabilityinmodulitheory}}]
The \emph{stability function} $M^\mu\colon \abs{\cX}\to [0,\infty]$ associated to $\mu$ is defined by
\[M^\mu(x)\coloneqq \sup \{\mu(\lambda)\st \lambda\in \abs{\Filt_\Q(\cX)},\ \ev_1(\lambda)=x\}.\]
\end{definition}

\begin{remark}
In \cite{_HalpernLeistner_Onthestructureofinstabilityinmodulitheory}, the stack $\Grad(\cX)$ is used instead of $\Grad_\Q(\cX)$. This is not an important difference because the numerical invariant $\mu$ is scale-invariant.
\end{remark}
\begin{definition}[Semistable locus]\label{definition: semistable locus}
The \emph{semistable locus} $\abs{\cX^\ss}$ with respect to the linear form $\ell$ on $\cX$ is the subset
\[\{x\in \abs{\cX}\st \ell(\lambda)\leq 0,\ \text{for all } \lambda\in \abs{\Filt_\Q(\cX)} \text{ with} \ \ev_1(\lambda)=x\}\subset \abs{\cX}.\]
If the semistable locus $\abs{\cX^\ss}$ is open, then it defines an open substack of $\cX$ denoted $\cX^\ss$.
\end{definition}

\begin{remark}
If $k$ is a field and $x\colon \Spec k\to \cX$ is a point, to see whether $x$ is semistable it suffices to check that $\ell(\lambda)\leq 0$ for $\lambda\in \qfilt(\cX,x\vert_{\overline k})$. This follows from $\Filt_\Q(\cX)\to \cX$ being representable and locally of finite presentation and $\ell$ being locally constant on $\Filt_\Q(\cX)$.
\end{remark}

\begin{definition}[{$\Theta$}-stratification defined by a linear form and a norm]\label{definition: ell q define theta-stratification}
We say that the pair $(\ell,q)$ \emph{defines a (weak) $\Theta$-stratification} if the following holds:
\begin{enumerate}
\item The semistable locus $\abs{\cX^\ss}\subset \abs{\cX}$ is open. If $s\colon \cX\to \Filt_\Q(\cX)$ is the “trivial filtration”, which is an open and closed immersion, we denote $\cS_0=s(\cX^\ss)$, an open substack of $\Filt_\Q(\cX)$ isomorphic to $\cX^\ss$.
\item For all $c\in \Q_{>0}$, the subset $\abs{\cS_c}\coloneqq \{\lambda \in \abs{\Filt_\Q(\cX)}\st \langle \lambda,\ell \rangle=1 \text{ and } \mu(\lambda)=M^\mu(\ev_1(\lambda))=\sqrt{c}\}$ of $\abs{\Filt_\Q(\cX)}$ is open, and thus it defines an open substack $\cS_c$ of $\Filt_\Q(\cX)$.
\item The family $(\cS_c)_{c\in \Q_{\geq 0}}$ is a (weak) $\Theta$-stratification of $\cX$ indexed by $\Q_{\geq 0}$, referred to as the $\Theta$-stratification \emph{induced by} $(\ell,q)$.
\end{enumerate}
\end{definition}

\begin{remark}
If $(\ell,q)$ defines a $\Theta$-stratification, then for all $c\in \Q_{\geq 0}$ we have the equality  $\abs{\cX_{\leq c}}=\{x\in \abs{\cX}\st M^\mu(x)\leq \sqrt{c}\}$ of $\abs{\cX}$, using the notation of \Cref{definition: Theta-stratification my version}. In particular $\cX^\ss=\cX_{\leq 0}$, and note as well that the centre $\cZ_0$ of the minimal stratum $\cS_0=\cX^\ss$ is canonically isomorphic to $\cX^\ss$.
\end{remark}

\begin{remark}
In Halpern-Leistner's definition, the stack $\Filt(\cX)$ is used instead of $\Filt_\Q(\cX)$. This is not an important difference, as we now explain. In \cite[Definition 4.1.3]{_HalpernLeistner_Onthestructureofinstabilityinmodulitheory}, the $\Theta$-stratification depends, in principle, on the choice of a complete set of representatives for $\pi_0(\Filt(\cX))/\N^*$, although two different choices give rise to isomorphic locally closed $\Theta$-strata. Since $\pi_0(\Filt_\Q(\cX))/\N^*=\pi_0(\Filt(\cX))/\N^*$, and since two connected components $a$ and $b$ of $\Filt_\Q(\cX)$ and $\Filt(\cX)$ respectively that represent the same class are isomorphic via the action of $\N^*$, we can instead take a complete set of representatives in $\pi_0(\Filt_\Q(\cX))$. We are using canonical set of representatives for the unstable strata, namely the components of $\Filt_\Q(\cX)$ on which $\ell$ takes the value $1$. This observation, together with \Cref{proposition: Halpern-Leistners definition of theta-stratification is equivalent to mine}, implies that the pair $(\ell,q)$ defines a $\Theta$-stratification in the sense of \cite[Definition 4.1.3]{_HalpernLeistner_Onthestructureofinstabilityinmodulitheory} if and only if it does so in the sense above, and that in that case the locally closed $\Theta$-strata that we get are isomorphic to Halpern-Leistner's. However, Definition \ref{definition: ell q define theta-stratification} has the advantage that it does not depend on noncanonical choices. Moreover, for \Cref{conjecture: on asymptotics of gradient flows affine GIT} it is important to take the canonical choice of HN filtration and not just consider it up to scaling.
\end{remark}

\begin{remark}[Characterisation of the HN filtration]\label{remark: characterisation of the HN filtration}
It follows from the definitions that, if $(\ell,q)$ defines a $\Theta$-stratification on $\cX$, then the HN filtration of a field-valued unstable point $x\in \cX(k)$ is the unique $\lambda\in \qfilt(\cX,x)$ such that $\langle \lambda,\ell\rangle\geq 1$ (equivalently, $\langle \lambda,\ell\rangle= 1$) and for all $\gamma\in \qfilt(\cX,x\vert_{\overline k})$ such that $\langle \gamma,\ell\rangle\geq 1$ we have $q(\lambda)\leq q(\gamma)$.
\end{remark}
\begin{proposition}[Compatibility with pullback]\label{proposition: compatibility of Theta-stratification with pullback}
Let $\cX'$ be an algebraic stack over an algebraic space $B$, satisfying \Cref{assumption: basic assumptions}, and let $h\colon \cX'\to \cX$ be either a closed immersion or a base change of a map between algebraic spaces like in \Cref{proposition: Filt Grad and base change}. Suppose that the pair $(\ell,q)$ defines a $\Theta$-stratification $(\cS_c)_{c\in \Q_{\geq 0}}$. Then $(h^*\ell,h^*q)$ defines a $\Theta$-stratification of $\cX'$, equal to the induced stratification $(h^*\cS_c)_{c\in \Q_{\geq 0}}$ of \Cref{propositionanddefinition: induced Theta-stratifications}.
\end{proposition}
\begin{proof}
The claim follows at once from the observation that, for every field $k$ and point $x\in \cX'(k)$, the induced map $\qfilt(\cX',x)\to \qfilt(\cX,h(x))$ is a bijection, compatible with the values of $(\ell,q)$ and $(h^*\ell,h^*q)$.
\end{proof}

\subsection{\texorpdfstring{$\Theta$}{Theta}-stratifications for stacks proper over a normed good moduli stack}
We now get to the main result (\Cref{theorem: theta stratification proper over gms}) about existence and properties of $\Theta$-stratifications that we will use. It is an extension of \cite[Theorem 5.6.1]{_HalpernLeistner_Onthestructureofinstabilityinmodulitheory} where we also establish existence of good moduli spaces for the centres of the strata, along with other improvements. We first introduce the notion of \emph{positive linear form on graded points}, a slight variant of \cite[Definition 5.3.1]{_HalpernLeistner_Onthestructureofinstabilityinmodulitheory}.

\begin{definition}[Positive linear form]\label{definition: positive linear form}
Let $\cX$ and $\cY$ be algebraic stacks over an algebraic space $B$, satisfying Assumption \ref{assumption: basic assumptions}. Let $f\colon \cY\to \cX$ be a proper representable morphism and let $\ell$ be a linear form on graded points of $\cY$. We say that $\ell$ is \emph{$f$-positive} provided that for all fields $k$ and all diagrams
\[\begin{tikzcd}[ampersand replacement=\&]
    {\P^1_k/\G_{m,k}} \& \cY \\
    {B\G_{m,k}} \& \cX
    \arrow[from=2-1, to=2-2]
    \arrow["f", from=1-2, to=2-2]
    \arrow[from=1-1, to=2-1]
    \arrow["\phi", from=1-1, to=1-2]
\end{tikzcd}\]
such that the induced map $\P^1_k/\G_{m,k}\to B\G_{m,k}\times_\cX \cY$ is finite, where the action of $\G_{m,k}$ on $\P^1_k$ is given by $t[a,b]=[ta,b]$ in projective coordinates and we denote $0=[0,1]$ and $\infty=[1,0]$, we have
\[\ell(\phi\vert_{\{0\}/\G_{m,k}})<\ell(\phi\vert_{\{\infty\}/\G_{m,k}}).\]
\end{definition}

\begin{example}
If $\cL$ is a rational line bundle on $\cY$ that is $f$-ample and $\ell=\langle -,\cL\rangle$, then $\ell$ is $f$-positive. Indeed, for all commutative diagrams as in \Cref{definition: positive linear form}, the pullback $\phi^*\cL$ is ample relative to $\P^1_k/\G_{m,k}\to B\G_{m,k}$, and the claim follows after embedding $\P^1_k$ in a bigger projective space  and looking at the weights of the corresponding $\G_{m,k}$-representation.
\end{example}
\begin{theorem}\label{theorem: theta stratification proper over gms}
Let $\cX$ be a noetherian algebraic stack with affine diagonal and a good moduli space $\pi\colon \cX\to X$, endowed with a norm on graded points $q$. Let $f\colon \cY\to \cX$ be a proper representable morphism and let $\ell$ be an $f$-positive linear form on graded points of $\cY$. Then
\begin{enumerate}
\item The pair $(\ell,f^*q)$ defines a $\Theta$-stratification $(\cS_c)_{c\in \Q_{\geq 0}}$ of $\cY$.
\item For every $c\in \Q_{\geq 0}$, the centre $\cZ_c$ has a good moduli space $\cZ_c\to Z_c$. We denote $\cY^\ss=\cZ_0$ and $Y^\ss=Z_0$.
\item For $c\in \Q_{\geq 0}$, let $\cX_c$ be the union of connected components of $\Grad_\Q(\cX)$ intersecting the image of $\cZ_c\to \Grad_\Q(\cY)\to \Grad_\Q(\cX)$. Then $\cX_c$ is quasi-compact and has a good moduli space $\cX_c\to X_c$.
\item For every $c\in \Q_{\geq 0}$, the induced map $Z_c\to X_c$ is proper.
\item If $\ell=\langle -,\cL\rangle$ for an $f$-ample rational line bundle $\cL$ on $\cY$, then $Y^\ss\to X$ is projective. If in addition the norm on graded points $q$ is algebraic, then for all $c\in \Q_{\geq c}$ the map $Z_c\to X_c$ is projective.
\end{enumerate}
\end{theorem}
\begin{remark}
The fact that we get a $\Theta$-stratification instead of a weak $\Theta$-stratification implies in particular rationality of HN filtrations. This comes at the expense of demanding $\cX$ to have a good moduli space instead of an adequate moduli space \cite{_Alper_Adequatemodulispacesandgeometricallyreductivegroupschemes}, which is the less restrictive notion in positive characteristic.
\end{remark}

\begin{example}\label{example: stratification on a blow-up}
Let $\cX$ be as in the statement of the theorem, and let $f\colon \cY\to \cX$ be the blow-up of $\cX$ along a closed substack $\cZ$. Let $\cE=f^{-1}(\cZ)$ be the exceptional divisor. The ideal sheaf $\cO_\cY(-\cE)$ of the exceptional divisor is $f$-ample. Therefore we may apply the theorem with $\ell=\langle - , \cO_\cY(-\cE)\rangle$ to get a $\Theta$-stratification of $\cY$. In this case, the semistable locus $\cY^\ss$ is the \emph{saturated blow-up} of $\cX$ along $\cZ$ \cite[Definition 3.2]{_Edidin_CanonicalreductionofstabilizersforArtinstackswithgoodmodulispaces}, by \cite[Proposition 3.17]{_Edidin_CanonicalreductionofstabilizersforArtinstackswithgoodmodulispaces} and \cite[Theorem 5.6.1, (2)]{_HalpernLeistner_Onthestructureofinstabilityinmodulitheory}. 

More generally, if $f$ is projective and $\ell=\langle -,\cL\rangle$, with $\cL$ an $f$-ample line bundle. Then $\cY^\ss$ is the \emph{saturated $\Proj$} \cite[Definition 3.1]{_Edidin_CanonicalreductionofstabilizersforArtinstackswithgoodmodulispaces}, $\cY^\ss=\Proj^\pi_\cX\left(\bigoplus_{n\in \N}f_*(\cL^{\otimes n})\right)$.
\end{example}
\begin{proof}
Note that $\cX$ satisfies \Cref{assumption: basic assumptions} with $B=X$ (\Cref{example: noetherian good moduli stack with affine diagonal}). Let $\mu$ denote the numerical invariant defined by $(\ell, f^*q)$.

\medskip
\noindent \emph{Step 1. The numerical invariant $\mu$ is strictly $\Theta$-monotone over $X$ \cite[Definition 5.2.1]{_HalpernLeistner_Onthestructureofinstabilityinmodulitheory} and strictly $S$-monotone \cite[Definition 5.5.7]{_HalpernLeistner_Onthestructureofinstabilityinmodulitheory} over $X$.}
\medskip

We first note that the map $\pi\colon \cX\to X$ is $\Theta$-reductive \cite[Definition 3.10]{_Alper_Existenceofmodulispacesforalgebraicstacks} and $S$-complete \cite[Definition 3.38]{_Alper_Existenceofmodulispacesforalgebraicstacks} by \cite[Theorem 5.4]{_Alper_Existenceofmodulispacesforalgebraicstacks}. We follow the argument in the proof of \cite[Proposition~5.3.3]{_HalpernLeistner_Onthestructureofinstabilityinmodulitheory}.
Let $R$ be a discrete valuation ring over $X$ with uniformiser $\pi$ and residue field $k$. Let $\overline{\operatorname{ST}}_R=\Spec(R[s,t]/(st-\pi))/\G_m$, where $s$ has weight $1$ and $t$ has weight $-1$ \cite[Section 3.5.1]{_Alper_Existenceofmodulispacesforalgebraicstacks}. Let $\cV$ be the stack $\Theta_R$ (resp. the stack $\overline{\operatorname{ST}}_R$), and let $\cV'=\cV\setminus \{(0,0)\}$. Suppose given a morphism $v\colon \cV'\to \cY$. Because $\cX$ is $\Theta$-reductive (resp. $S$-complete), the morphism $f\circ v\colon \cV'\to \cX$ extends to a map $u\colon \cV\to \cX$. Let $\cW$ be the schematic image of $\cV'$ inside the pullback $\cV\times_{u,\cX,f} \cY$. We have a diagram
\[\begin{tikzcd}[ampersand replacement=\&]
    \& \cW \& \cY \\
    {\cV'} \& \cV \& \cX
    \arrow["f", from=1-3, to=2-3]
    \arrow["u"', from=2-2, to=2-3]
    \arrow["p", from=1-2, to=2-2]
    \arrow["{\overline v}"', from=1-2, to=1-3]
    \arrow[hook, from=2-1, to=2-2]
    \arrow[hook, from=2-1, to=1-2]
    \arrow["v", shift left, curve={height=-30pt}, from=2-1, to=1-3]
\end{tikzcd}\]
and we want to show that $\cW$ and $\overline v$ satisfy the conditions of \cite[Definitions 5.2.1 and 5.5.7]{_HalpernLeistner_Onthestructureofinstabilityinmodulitheory}, of which the only nontrivial one is $(3)$. For this, suppose given a commutative square
\[\begin{tikzcd}[ampersand replacement=\&]
    {\P^1_k/\G_{m,k}} \& \cW \\
    {B\G_{m,k}} \& \cV
    \arrow["g", from=2-1, to=2-2]
    \arrow["p", from=1-2, to=2-2]
    \arrow[from=1-1, to=2-1]
    \arrow["\phi", from=1-1, to=1-2]
\end{tikzcd}\]
where $g$ is a 
multiple of the canonical graded point $\{(0,0)\}/\G_{m,k}\to\cV$ and the induced $\P^1_k/\G_{m,k}\to B\G_{m,k}\times_\cV \cW$ is finite. We want to show that
\[\ell_c(\overline v\circ \phi\vert_{\{0\}/\G_{m,k}})<\ell_c(\overline v\circ \phi\vert_{\{\infty\}/\G_{m,k}}).\]
The map $\P^1_k/\G_{m,k}\to B\G_{m,k}\times_\cX \cY$ is also finite, so the inequality follows from $\ell$ being $f$-positive.

We remark that properness of $p$ is one of the conditions of \cite[Definitions 5.2.1 and 5.5.7]{_HalpernLeistner_Onthestructureofinstabilityinmodulitheory} and it is guaranteed by properness of $f$.

\medskip
\noindent \emph{Step 2. The pair $(\ell, f^*q)$ defines a weak $\Theta$-stratification.}
\medskip

Since $\cY$ is quasi-compact, the numerical invariant $\mu$ satisfies HN-boundedness \cite[Proposition 4.4.2]{_HalpernLeistner_Onthestructureofinstabilityinmodulitheory}. Therefore $(\ell,f^*q)$ defines a weak $\Theta$-stratification if $\mu$ is strictly $\Theta$-monotone \cite[Theorem 5.2.3]{_HalpernLeistner_Onthestructureofinstabilityinmodulitheory}. Thus the claim follows from Step 1.

\medskip
\noindent \emph{Step 3. The weak $\Theta$-stratification defined by $(\ell, f^*q)$ is a $\Theta$-stratification.}
\medskip

Checking that a finite morphism $\cU\to \cV$ is a closed immersion can be done after base change along all geometric points $\Spec k\to \cV$. Since the stratification is preserved by base change along any map $X'\to X$ \Cref{proposition: compatibility of Theta-stratification with pullback}, we may assume that $X=\Spec k$ is the spectrum of an algebraically closed field. If $\Spec k$ is of characteristic $0$, then the weak $\Theta$-stratification is a $\Theta$-stratification by \cite[Theorem 5.6.1]{_HalpernLeistner_Onthestructureofinstabilityinmodulitheory}. Suppose $k$ is of positive characteristic $p$. The stack $\cX$ is of the form $\cX=\Spec A/G$, where $G$ is a linearly reductive group over $k$, by \ref{corollary: fibres of good moduli space are quotient stacks}.
Therefore we have $\cY=Y/G$, where $Y\to \Spec A$ is $G$-equivariant and projective.
By Nagata's Theorem \cite[Chapter IV, Section 3, Theorem 3.6]{_Demazure_Gabriel_Groupesalgebriques}, the identity component $G^\circ$ is a group of multiplicative type and $p$ does not divide the order of $G/G^\circ$. By \cite[Lemma 2.1.7, (2)]{_HalpernLeistner_Onthestructureofinstabilityinmodulitheory} it is enough to prove that for any $c\in \Q_{\geq 0}$ and for any $k$-point $\lambda$ of $\cS_c$, the induced map $\phi\colon\Lie(\Aut_{\cS_c}(\lambda))\to \Lie(\Aut_\cY(\ev_1(\lambda))$ on Lie algebras is surjective. The filtration $\lambda$ is contained in some open substack of $\Filt_\Q(\cY)$ of the form $Y^{\gamma,+}/P(\gamma)$ for some one-parameter subgroup $\gamma$ of $G$, by \cite[Theorem 1.4.8]{_HalpernLeistner_Onthestructureofinstabilityinmodulitheory}, and it corresponds to a point $z\in Y^{\gamma,+}(k)$. Thus $\Aut_{\cS_c}(\lambda)=\Stab_{P(\gamma)}(z)$, while $\Aut_{\cX}(\ev_1(\lambda))=\Stab_G(z)$. Both groups have the same identity component, equal to $\Stab_{G^\circ}(z)^{\circ}$, since $G^\circ$ is contained in $P(\lambda)$. Thus the map $\phi$ on Lie algebras is actually an isomorphism. 

\medskip
\noindent \emph{Step 4. The semistable locus $\cY^\ss$ has a good moduli space $\cY^\ss \to Y^\ss$ and $Y^\ss\to X$ is separated.}
\medskip

By \cite[Theorem 5.4]{_Alper_Existenceofmodulispacesforalgebraicstacks}, $\cX$ is $\Theta$-reductive and $S$-complete. Therefore $\mu$ is strictly $\Theta$-monotone and strictly $S$-monotone over $X$ by Step 1. By \cite[Theorem 5.5.8]{_HalpernLeistner_Onthestructureofinstabilityinmodulitheory}, the semistable locus $\cY^\ss$ is $\Theta$-reductive and $S$-complete over $X$. Therefore it is enough to show, by \cite[Theorem 5.4]{_Alper_Existenceofmodulispacesforalgebraicstacks}, that the stabiliser of every closed point of $\cY^\ss$ is linearly reductive. Let $k$ be an algebraically closed field and let $x\in \cY^\ss(k)$ be a closed point. By \cite[Proposition 3.47]{_Alper_Existenceofmodulispacesforalgebraicstacks}, the automorphism group $\Aut(x)$ is geometrically reductive. If the characteristic of $k$ is $0$, then $\Aut(x)$ is also linearly reductive. If $k$ is of positive characteristic $p$, then $f(x)$ specialises to a $k$-point $y$ closed in the fibre of $\pi\colon \cX\to X$, whose stabiliser $\Aut(y)$ is thus linearly reductive. By Nagata's Theorem \cite[Chapter IV, Section 3, Theorem 3.6]{_Demazure_Gabriel_Groupesalgebriques}, the identity component $\Aut(y)^\circ$ is a group of multiplicative type and $p$ does not divide the order of $\Aut(y)/\Aut(y)^\circ$. Since $\Aut(x)$ is a subgroup of $\Aut(y)$, the same holds for $\Aut(x)$ and it is thus linearly reductive.

\medskip
\noindent \emph{Step 5. The map $Y^\ss\to X$ is proper.}
\medskip

 Since $\cY\to \cX$ is proper and $\cX\to X$ is universally closed \cite[Theorem 4.16]{_Alper_GoodmodulispacesforArtinstacks}, we have that $\cY\to X$ is universally closed. Now the Semistable Reduction Theorem \cite[Corollary 6.12]{_Alper_Existenceofmodulispacesforalgebraicstacks} implies that $\cY^\ss\to X$ is also universally closed. Therefore $Y^\ss\to X$ is universally closed. Since it is also separated, by Step 4, and of finite type, because $f$ and $\pi$ are, it is proper.

\medskip
\noindent \emph{Step 6. The map $Y^\ss\to X$ is projective if $\ell=\langle -,\cL\rangle$} for a rational line bundle $\cL$.
\medskip

We may scale $\cL$ up to assume that it is a line bundle, since this does not change the semistable locus. Then we have, by \cite{_HalpernLeistner_Onthestructureofinstabilityinmodulitheory}, that $Y^\ss=\Proj_X\left( \pi_*f_*(\oplus_{n\in\N}\cL^{\otimes n}) \right)$.

\medskip
\noindent \emph{Step 7. For all $c\in \Q_{\geq 0}$, the stack $\cX_c$ is quasi-compact and has a good moduli space $\cX_c\to X_c$.}
\medskip

Since $\cZ_c$ is quasi-compact, and because of the existence of a norm on graded points on $\cX$, we have by \cite[Proposition 3.8.2]{_HalpernLeistner_Onthestructureofinstabilityinmodulitheory} that $\cX_c$ is quasi-compact. The claim follows from \Cref{lemma: quasi-compact components of GradcX are affine over cX}.

\medskip
\noindent \emph{Step 8. The centres $\cZ_c$ have good moduli spaces $\cZ_c\to Z_c$ for all $c\in \Q_{c\geq 0}$ and $Z_c\to X_c$ is proper.}
\medskip

For $c=0$, the claim is the content of Step 4 and Step 5. If $c\in \Q_{>0}$, let $\overline \cZ_c$ be the union of those connected components of $\Grad_\Q(\cY)$ intersecting $\cZ_c$. Again, because of the existence of a norm on graded points on $\cX$ and $\cY$, we have by \cite[Proposition 3.8.2]{_HalpernLeistner_Onthestructureofinstabilityinmodulitheory} that $\overline \cZ_c$ is quasi-compact. Let $f_c\colon \overline \cZ_c\to \cX_c$ be the restriction of $\Grad_\Q(f)$ to $\overline \cZ_c$ and $\cX_c$. By \Cref{proposition: f proper implies Grad f proper}, the map $f_c$ is representable and proper.

Let us denote $\ell\vert_{\overline\cZ_c}$ the pullback of the linear form $\ell$ on $\cY$ along $\overline\cZ_c\to \Grad_\Q(\cY)\to \cY$. We denote $\ell_{f^*q}$ the linear form on $\Grad_\Q(\cY)$ induced by the norm $f^*q$ on $\cY$ (\Cref{definition: canonical linear form on graded points}), and $\ell_{f^*q}\vert_{\overline \cZ_c}$ its restriction to $\overline \cZ_c$. Note that $\ell_{f^*q}\vert_{\overline \cZ_c}=(f_c)^*(\ell_{q}\vert_{\cX_c})$ is the pullback of the linear form $\ell_q\vert_{\cX_c}$ on $\cX_c$. We will consider the \emph{shifted linear form} 
\begin{equation}\label{equation: shifted linear form}
\ell_c\coloneqq \ell\vert_{\overline \cZ_c}-c\ell_{f^*q}\vert_{\overline\cZ_c}
\end{equation}
on $\overline\cZ_c$. We will use the following.
\begin{theorem}[Linear Recognition Theorem \cite{tame-theta-stratifications}]\label{theorem: linear recognition} 
The centre $\cZ_c$ is the semistable locus inside $\overline\cZ_c$ with respect to the shifted linear form $\ell_c$ on $\overline \cZ_c$.
\end{theorem}

Therefore, we have a representable proper morphism $f_c\colon \overline\cZ_c\to \cX_c$ and a linear form $\ell_c$ and a norm on graded points $q\vert_{\cX_c}$ on $\cX_c$. By steps 4 and 5 applied to $f_c$ in place of $f$, it is enough to show that $\ell_c$ is $f_c$-positive. For this, suppose given a commutative square
\[\begin{tikzcd}[ampersand replacement=\&]
    {\P^1_k/\G_{m,k}} \& {\overline \cZ_c} \\
    {B\G_{m,k}} \& {\cX_c}
    \arrow["g", from=2-1, to=2-2]
    \arrow[from=1-2, to=2-2]
    \arrow[from=1-1, to=2-1]
    \arrow["\phi", from=1-1, to=1-2]
\end{tikzcd}\]
where the induced $\P^1_k/\G_{m,k}\to B\G_{m,k}\times_{\cX_c} \overline \cZ_c$ is finite. We want to show that
\[\ell_c(\phi\vert_{\{0\}/\G_{m,k}})<\ell_c(\phi\vert_{\{\infty\}/\G_{m,k}}).\]
First, note that 
\[\left(\ell_{f^*q}\vert_{\overline \cZ_c}\right)(\phi\vert_{\{0\}/\G_{m,k}})=\ell_q(g)=\left(\ell_{f^*q}\vert_{\overline \cZ_c}\right)(\phi\vert_{\{\infty\}/\G_{m,k}}).\]
On the other hand, the map $\overline \cZ_c\to \cX_c\times_\cX \cY$ is proper, since $\overline \cZ_c\to \cX_c$ and $\cY\to \cX$ are proper; and also affine, since $\overline \cZ_c\to \cY$ and $\cX_c\to \cX$ are affine. Thus the induced map $\P^1_k/\G_{m,k}\to B\G_{m,k}\times_\cX \cY$ is finite. Therefore we have
\begin{align*}
&\ell_c(\phi\vert_{\{\infty\}/\G_{m,k}})-\ell_c(\phi\vert_{\{0\}/\G_{m,k}})=\\
&\ell\vert_{\overline\cZ_c}(\phi\vert_{\{\infty\}/\G_{m,k}})-\ell\vert_{\overline\cZ_c}(\phi\vert_{\{0\}/\G_{m,k}})>0,
\end{align*}
because $\ell$ is $f$-positive.

\medskip
\noindent \emph{Step 9. If $\ell=\langle -,\cL\rangle$ for an $f$-ample line bundle $\cL$ and $q$ is algebraic, then for all $c\in \Q_{\geq 0}$ the map $Z_c\to X_c$ is projective.}
\medskip

We begin by showing that the pullback $\cM\coloneqq (\Grad_\Q(\cY)\to \cY)^*\cL$ is relatively ample with respect to $g=\Grad_\Q(f)\colon \Grad_\Q(\cY)\to \Grad_\Q(\cX)$.
For this, we may assume that $\cL$ is a line bundle and we want to show that the canonical map $\cY\to \Proj\left(\bigoplus_{n\in \N} g_*(\cM^{\otimes n})\right)$ is everywhere defined and an open immersion (although it is actually an isomorphism by properness of $g$). This can be checked étale locally on $\Grad_\Q(\cX)$. Chose a surjective, affine and strongly étale morphism $\rho\colon \Spec A/\GL_n\to \cX$ (\Cref{theorem: local structure}).
Pulling back along $\rho$ we get a cube
\[\begin{tikzcd}[ampersand replacement=\&,column sep=tiny,row sep=small]
    \& {\bigsqcup_{\lambda\in C}Y^{\lambda,0}/L(\lambda)} \&\& {\bigsqcup_{\lambda\in C}(\Spec A)^{\lambda,0}/L(\lambda)} \\
    {\Grad_\Q(\cY)} \&\& {\Grad_\Q(\cX)} \\
    \& {Y/\GL_n} \&\& {\Spec A/\GL_n} \\
    \cY \&\& \cX
    \arrow[from=2-1, to=4-1]
    \arrow[from=4-1, to=4-3]
    \arrow[from=3-2, to=3-4]
    \arrow[from=1-4, to=3-4]
    \arrow[from=2-1, to=2-3]
    \arrow[from=2-3, to=4-3]
    \arrow[from=3-4, to=4-3]
    \arrow["c"', from=3-2, to=4-1]
    \arrow["d"{description}, from=1-4, to=2-3]
    \arrow["a"', from=1-2, to=2-1]
    \arrow["b", from=1-2, to=1-4]
    \arrow[from=1-2, to=3-2]
\end{tikzcd}\]
where we are using \cite[Theorem 1.4.7]{_HalpernLeistner_Onthestructureofinstabilityinmodulitheory} for the description of the stack of graded points of a quotient stack. Here, $Y$ is a scheme acted on by $\GL_n$ with an equivariant map into $\Spec A$. We are denoting $C=\Gamma^\Z(T)/W$, where $T$ is the standard maximal torus of $\GL_n$ and  $W$ is the Weyl group. For $\lambda\in \Gamma^\Z(T)$, $Y^{\lambda,0}$ and $(\Spec A)^{\lambda,0}$ denote the fixed point loci by the cocharacter $\lambda$, and $L(\lambda)$ is the centraliser of $\lambda$ in $\GL_n$.

The arrow $d$ is an étale cover of $\Grad_\Q(\cX)$. Thus we want to show that $a^*\cM$ is ample relative to $b$. In fact, since the $(\Spec A)^{\lambda,0}$ are affine, it is enough to show that $a^*\cM\vert_{Y^{\lambda,0}}$ is ample on the scheme $Y^{\lambda,0}$. Now, $a^*\cM\vert_{Y^{\lambda,0}}=\left(\cL\vert_Y\right)\vert_{Y^{\lambda,0}}$, the line bundle $\cL\vert_Y$ is ample and $Y^{\lambda,0}\to Y$ is affine, so $a^*\cM\vert_{Y^{\lambda,0}}$ is ample too. This shows that $\cM$ is relatively $g$-ample. In particular, if we let $\cM_c=\cM\vert_{\overline \cZ_c}$, then $\cM_c$ is relatively ample with respect to $f_c$.

Since $q$ is algebraic, the linear form $\ell_{f^*q}\vert_{\overline\cZ_c}$ is induced by a rational line bundle of the form $f_c^*\cN$, where $\cN$ is a rational line bundle on $\cX_c$. Thus $\ell_c$ is induced by the rational line bundle $\cM_c-cf_c^*\cN$, which is relatively $f_c$-ample because $\cM_c$ is. Therefore the claim follows from Step 6 applied to $f_c$ in place of $f$.
\end{proof}

In the proof we used the following lemmas.

\begin{lemma}\label{lemma: quasi-compact components of GradcX are affine over cX}
Let $\cX$ be a noetherian algebraic stack with affine diagonal and a good moduli space $\pi\colon \cX\to X$. For every quasi-compact open and closed substack $\cZ$ of $\Grad_\Q(\cX)$, the composition
\[\cZ\to \Grad_\Q(\cX)\xrightarrow{u} \cX\]
is affine. In particular, $\cZ$ has a good moduli space $Z$ which is affine over $X$.
\end{lemma}
\begin{proof}
Let $\Spec A/\GL_n\to \cX$ be affine, surjective and strongly étale (\Cref{theorem: local structure}).
We have a cartesian diagram
\[\begin{tikzcd}[ampersand replacement=\&]
    {\Grad_\Q(\Spec A/\GL_n)} \& {\Grad_\Q(\cX)} \\
    {\Spec A/\GL_n} \& \cX
    \arrow[from=1-2, to=2-2]
    \arrow[""{name=0, anchor=center, inner sep=0}, from=2-1, to=2-2]
    \arrow[from=1-1, to=2-1]
    \arrow[from=1-1, to=1-2]
    \arrow["\ulcorner"{anchor=center, pos=0.125}, draw=none, from=1-1, to=0]
\end{tikzcd}\]
and $\Grad_\Q(\Spec A/\GL_n)$ is a disjoint union of schemes of the form $\Spec(A)^{\lambda,0}/L(\lambda)$ with $\lambda$ a rational cocharacter of $\GL_n$. We conclude by \Cref{lemma: Grad has gms affine case} and descent.
\end{proof}
\begin{lemma}\label{lemma: Grad has gms affine case}
Let $A$ be a commutative ring, and consider an action of $\GL_N$ on $X=\Spec A$ (over $\Z$) such that $X/\GL_N$ has a good moduli space. Let $\lambda\colon \G_m\to \GL_N$ be a cocharacter. Then the natural map $X^{\lambda, 0}/L(\lambda)\to X/\GL_N$ is affine, where $X^{\lambda, 0}$ is the fixed point locus of the $\G_m$-action on $X$ induced by $\lambda$ and $L(\lambda)$ is the centraliser of $\lambda$.
\end{lemma}
\begin{proof}
There is a cartesian square 
\[ \begin{tikzcd}
\GL_N\times^{L(\lambda)} X^{\lambda, 0} \arrow[r,""]\arrow[d,swap,""]\arrow[dr, phantom, "\ulcorner", very near start] & X \arrow[d,""] \\
X^{\lambda, 0}/L(\lambda) \arrow[r,""]& X/\GL_N
\end{tikzcd}
\]
where $\GL_N\times^{L(\lambda)} X^{\lambda, 0}$ is the stack quotient of $\GL_N\times X^{\lambda, 0}$ by the diagonal action of $L(\lambda)$. Since the action is free, $\GL_N\times^{L(\lambda)} X^{\lambda, 0}$ is an algebraic space. Now, $L(\lambda)$ is isomorphic to a product of $\GL_{N_i}$'s and it is thus geometrically reductive \cite[Definition 9.1.1]{_Alper_Adequatemodulispacesandgeometricallyreductivegroupschemes}. Since $\GL_N\times X^{\lambda,0}=\Spec B$ is affine, the $L(\lambda)$-invariants give an adequate moduli space $\GL_N\times^{L(\lambda)} X^{\lambda, 0}\to \Spec \left(B^{L(\lambda)}\right)$  \cite[Theorem 9.1.4]{_Alper_Adequatemodulispacesandgeometricallyreductivegroupschemes}. By universality for adequate moduli spaces \cite[Theorem 3.12]{_Alper_Theetalelocalstructureofalgebraicstacks}, we get an isomorphism $\GL_N\times^{L(\lambda)} X^{\lambda,0}=\Spec B^{L(\lambda)}$. Therefore $\GL_N\times^{L(\lambda)} X^{\lambda, 0}$ is affine and we are done by descent.
\end{proof}

\begin{remark}
Since $L(\lambda)$ is geometrically reductive and $X^{\lambda,0}$ is affine, taking $L(\lambda)$-invariants gives an adequate moduli space for $X^{\lambda, 0}/L(\lambda)$. However, unless $A$ is of characteristic $0$, an extra argument is needed to show that the adequate moduli space is indeed a good moduli space.
\end{remark}

\section{Sequential stratifications and the iterated balanced filtration}
The goal of this section is the construction of the \emph{balancing stratification} for every noetherian normed good moduli stack $\cX$ with affine diagonal (\Cref{theorem: sequential stratifications for good moduli stacks} and \Cref{definition: balancing stratification}) and derive from it the definition of the \emph{iterated balanced filtration} of a point of $\cX$ (\Cref{definition: iterated balanced filtration}).

We start by defining a first approximation of the iterated balanced filtration, simply called the \emph{balanced filtration} (\Cref{definition: balanced filtration}). We then construct a stack of sequential filtrations $\Filt_{\Q^\infty}(\cX)$ for a stack $\cX$, and define a notion of \emph{sequential stratification} (\Cref{definition: sequential stratification}), roughly a weak analogue of $\Theta$-stratification with the stack $\Filt_{\Q^\infty}(\cX)$ used instead of $\Filt_\Q(\cX)$. For the purpose of using induction in the construction of the balancing stratification, we introduce the concept of \emph{central rank} of a stack (\Cref{definition: central rank}) and show that it increases after taking the centre of an unstable stratum in a blow-up (\Cref{lemma: centres of unstable strata have bigger central rank}). After our main construction (\Cref{theorem: sequential stratifications for good moduli stacks}), we prove some functorial properties of the balancing stratification (\Cref{proposition: compatibility of balancing stratification with pullback}). We finish the section with a collection of natural examples of normed good moduli stacks in moduli theory, including GIT quotients, moduli of Bridgeland semistable objects and moduli of K-semistable Fano varieties.

\subsection{The balanced filtration}

Let $\cX$ be an algebraic stack over an algebraic space $B$, satisfying Assumption \ref{assumption: basic assumptions}. Let $\cZ\to \cX$ be a closed substack, let $k$ be a field, let $x$ be a $k$-point of $\cX$ and let $\lambda\in \qfilt(\cX,x)$ be a rational filtration of $x$. We would like to formalise the idea of velocity at which $\lambda(t)$ tends to $\cZ$ when $t$ tends to $0$. For that, we first write $\lambda=\gamma/m$ with $\gamma\in\Z\dash\Filt(\cX,x)$ an integral filtration and $m\in \Z_{>0}$. We then form the pullback 
\[ \begin{tikzcd}
R/\G_{m,k} \arrow[r,""]\arrow[d,swap,""]\arrow[dr, phantom, "\ulcorner", very near start] & \cZ \arrow[d,""] \\
\A^1_k/\G_{m,k} \arrow[r,"\gamma"]& \cX
\end{tikzcd}
\]
of $\cZ\to \cX$ along $\gamma$, which is given by a $\G_{m,k}$-equivariant closed subscheme $R$ of $\A^1_k$, thus necessarily of the form $R=\Spec(k[x]/(x^n))$ for some $n\in \N\cup \{\infty\}$. The following concept is used by Kempf in \cite{_Kempf_InstabilityinInvariantTheory}.

\begin{definition}[Kempf's intersection number]\label{definition: Kempfs intersection number}
We define \emph{Kempf's intersection number} (or simply \emph{Kempf's number}) to be $\langle \lambda, \cZ \rangle \coloneqq n/m\in \Q_{\geq 0}\cup \{\infty\}$.
\end{definition}

The definition does not depend on the presentation $\lambda=\gamma/m$ because if $l\in \Z_{>0}$, then $\langle l\gamma, \cZ \rangle =ln$. More generally, the linearity property $\langle c\lambda,\cZ\rangle=c\langle \lambda,\cZ\rangle$ holds for all $c\in \Q_{> 0}$. We have $\langle \lambda, \cZ\rangle=\infty$ precisely when $x$ is in $\cZ$, and $\langle \lambda,\cZ \rangle = 0$ if and only if $\ev_0(\lambda)$ is not in $\cZ$. If $x$ is not in $\cZ$ but $\ev_0(\lambda)$ is in $\cZ$, then $\langle \lambda,\cZ\rangle$ is a positive rational number that can be thought of as the velocity at which $\lambda$ approaches $\cZ$.

\begin{proposition}\label{proposition: Kempfs intersection number and blowups}
Suppose $x$ is not in $\cZ$. Let $p\colon \Bl_\cZ\cX\to \cX$ be the blow-up of $\cX$ along $\cZ$, let $\cL$ be the standard $p$-ample line bundle on $\Bl_\cZ \cX$, and let $x'\in \Bl_\cZ\cX(k)$ be a lift of $x$ to $\Bl_\cZ\cX$. Let $\lambda'\in \qfilt(\Bl_\cZ\cX,x')$ be the unique lift of $\lambda$ to a rational filtration of $x'$, which exists by \Cref{proposition: proper map induces bijection of set of filtrations}. Then $\langle \lambda', \cL\rangle =\langle \lambda, \cZ \rangle$.
\end{proposition}
\begin{proof}

By linearity of $\langle -,\cZ\rangle$, we may assume that $\lambda$ is integral. By definition, $\cL$ is the ideal sheaf of the exceptional divisor $\cE=p^{-1}(\cZ)$. We have a diagram
\[\begin{tikzcd}[ampersand replacement=\&]
    {R/\G_{m,k}} \& \cE \& \cZ \\
    {\Theta_k} \& {\Bl_\cZ\cX} \& \cX,
    \arrow[from=1-2, to=1-3]
    \arrow[from=1-3, to=2-3]
    \arrow["p", from=2-2, to=2-3]
    \arrow[from=1-2, to=2-2]
    \arrow["\lrcorner"{anchor=center, pos=0.125}, draw=none, from=1-2, to=2-3]
    \arrow["{\lambda'}", from=2-1, to=2-2]
    \arrow[from=1-1, to=1-2]
    \arrow[from=1-1, to=2-1]
    \arrow["\lrcorner"{anchor=center, pos=0.125}, draw=none, from=1-1, to=2-2]
\end{tikzcd}\]
where $R=\Spec k[x]/(x^n)$ for some $n\in \N$, the variable $x$ having weight $-1$. There is a natural injection $\cL\to \cO_{\Bl_\cZ \cX}$. Pulling back along $\lambda'$ we get a map $(\lambda')^*\cL\to \cO_{\Theta_k}$ whose image $\cI=(x^n)$ is the ideal sheaf of $R/\G_{m,k}$. Since there is a surjective map $(\lambda')^*\cL \to \cI$ and both source and target are line bundles, the map should be an isomorphism. Thus $(\lambda')^*\cL\cong \cI\cong \cO_{\Theta_k}(-n)$, where $\cO_{\Theta_k}(-n)$ is the pullback of $\cO_{B\G_{m,k}}(-n)$ along the structure map $\Theta_k\to B\G_{m,k}$, because $x^n$ has weight $-n$. Therefore  $\langle \lambda',\cL\rangle=-\wt \left(\left((\lambda')^*\cL\right)\vert_{B\G_{m,k}}\right)=n=\langle \lambda, \cZ \rangle$, as desired.
\end{proof}
The following result is a generalisation of Kempf's Theorem \cite[Theorem 3.4]{_Kempf_InstabilityinInvariantTheory} to stacks with good moduli space. See also \cite[Example 5.3.7]{_HalpernLeistner_Onthestructureofinstabilityinmodulitheory}.

\begin{theorem}[Kempf]\label{theorem: Kempf}
Let $\cX$ be a noetherian algebraic stack with affine diagonal and a good moduli space $\pi\colon \cX\to X$. Let $k$ be a field, let $\cZ\to \cX$ be a closed substack and let $x$ be a $k$-point of $\cX\setminus \cZ$. Then:
\begin{enumerate}
\item\label{item 1 kempf theorem} The intersection  $\abs{\cZ}\cap \abs{\pi^{-1}\pi(x)}$ is nonempty if and only if there is a filtration $\lambda\in \Z\dash\Filt(\cX,x)$ with $\ev_0(\lambda)$ in $\cZ$.

\item\label{item 2 kempf theorem} Suppose that $\cX$ is endowed with a norm on graded points. If $\abs{\cZ}\cap \abs{\pi^{-1}\pi(x)}\neq \varnothing$, then there is a unique $\lambda\in\qfilt(\cX,x)$ with $\langle\lambda,\cZ\rangle\geq 1$ and such that for all $\gamma\in \qfilt(\cX,x\vert_{\overline k})$ with $\langle \gamma, \cZ \rangle \geq 1$ we have $\norm{\lambda}\leq \norm{\gamma}$.
\end{enumerate}
\end{theorem}
\begin{remark}
Note that we do not require $k$ to be perfect as in the original Kempf's Theorem. This is possible because we are working with good moduli spaces instead of the adequate moduli spaces, that are more general in positive characteristic.
\end{remark}
\begin{proof}
By replacing $\cX$ by $\pi^{-1}\pi(x)$ we may assume that $X=\Spec k$. By \Cref{corollary: fibres of good moduli space are quotient stacks}, $\cX=\Spec A/ \GL_{n,k}$, where $A$ is a $k$-algebra of finite type.

We first show \Cref{item 1 kempf theorem} assuming $k$ is algebraically closed. Suppose $\cZ\neq \varnothing$. Since $\cX$ has a unique closed point \cite[Proposition 9.1]{_Alper_GoodmodulispacesforArtinstacks}, we also have that $\overline{\{x\}}\cap \abs{\cZ}$ is nonempty. By \cite[Theorem 1.4]{_Kempf_InstabilityinInvariantTheory}, there exists $\lambda\in \qfilt(\cX,x)$ such that $\ev_1(\lambda)=x$ and $\ev_0(\lambda)$ is in $\cZ$.

Now we show \Cref{item 2 kempf theorem} for any $k$. We are given a norm on graded points on $\cX$. Let $x'$ be a lift of $x$ to the blow-up $\Bl_\cZ \cX$. By Theorem \ref{theorem: theta stratification proper over gms} and \Cref{example: stratification on a blow-up} there is a $\Theta$-stratification on $\Bl_\cZ \cX$. By \Cref{proposition: proper map induces bijection of set of filtrations} we have $\qfilt(\cX,x\vert_{\overline k})=\qfilt(\Bl_\cZ \cX,x'\vert_{\overline k})$, and by \Cref{item 1 kempf theorem} in the algebraically closed case and \Cref{proposition: Kempfs intersection number and blowups} we have that $x'$ is semistable if and only if $\abs{\cZ}=\varnothing$. If $x'$ is unstable, then its HN filtration is the unique $\lambda\in \qfilt(\cX,x)$ in \Cref{item 2 kempf theorem}, also by \Cref{proposition: Kempfs intersection number and blowups}.

By choosing any norm on cocharacters of $\GL_{n,k}$ and pulling it back to $\cX$, we see that \Cref{item 2 kempf theorem} readily implies \Cref{item 1 kempf theorem} for any $k$.
\end{proof}

We now shift attention to the case where $\cZ$ is the locus $\cX^{\max}$ of maximal dimension of stabiliser groups of $\cX$, whose definition we recall below. Suppose $\cX$ is noetherian with affine diagonal. For each $d\in \N$, the set $\{x\in \abs{\cX}\st \dim \Aut(x)\geq d\}$ is closed \cite[Exposé VIb, Proposition 4.1]{SGA3}. Therefore, by quasi-compactness of $\cX$, it makes sense to define:
\begin{definition}[Maximal stabiliser dimension]\label{definition: maximal stabiliser dimension}
Let $\cX$ be a noetherian algebraic stack with affine diagonal. The \emph{maximal stabiliser dimension} $d(\cX)\in \N$ of $\cX$ is the maximal dimension of the stabiliser group of a point of $\cX$. For the empty stack we set $d(\varnothing)=-\infty$.
\end{definition}

As a topological space, $\abs{\cX^{\max}}=\{x\in \abs{\cX}\st \dim \Aut(x)=d(\cX)\}$, and it is a closed subset of $\abs{\cX}$. It is a nontrivial result \cite[Proposition C.5]{_Edidin_CanonicalreductionofstabilizersforArtinstackswithgoodmodulispaces} that if $\cX$ is a noetherian good moduli stack with affine diagonal, then $\abs{\cX^{\max}}$ can be given a natural structure of closed substack of $\cX$, denoted $\cX^{\max}$ and called the \emph{maximal dimension stabiliser locus} of $\cX$. It can be characterised étale locally by the property that, if $\cX=X/G$, where $X$ is an algebraic space and $G\to \Spec \Z$ is an affine flat group scheme of finite type that is either diagonalisable or a Chevalley group, with fibres of dimension $d(\cX)$, then $\cX^{\max}=X^{G_\circ}/G$, where $G_\circ$ is the reduced identity component of $G$ \cite[Section C.2]{_Edidin_CanonicalreductionofstabilizersforArtinstackswithgoodmodulispaces}. In general, the functorial definition of $\cX^{\max}$ is as follows. A map $f\colon T\to \cX$ factors through $\cX^{\max}$ if and only if the pullback $f^*\cI_\cX$ of the inertia of $\cX$ has a smooth closed subgroup all whose fibres are connected and of dimension $d(\cX)$.

One reason why Edidin-Rydh's stack structure on $\cX^{\max}$ is better behaved than the reduced structure is that it behaves well with respect to base change. If $\cY\to \cX$ is a closed immersion and $d(\cY)=d(\cX)$, then $\cY^{\max}=\cY\times_\cX \cX^{\max}$. Similarly, if $\cY$ has a good moduli space, $\cY\to \cX$ is representable, étale and separated, and $d(\cY)=d(\cX)$, then $\cY^{\max}=\cY\times_\cX \cX^{\max}$. See \cite[Proposition C.5]{_Edidin_CanonicalreductionofstabilizersforArtinstackswithgoodmodulispaces}. We get to the main definition of this section.

\begin{definition}[The balanced filtration]\label{definition: balanced filtration}
Let $\cX$ be a normed noetherian algebraic stack with affine diagonal and a good moduli space $\pi\colon \cX \to X$. Let $x\colon \Spec k\to \cX$ be a field-valued point and denote $\cF=\pi^{-1}\pi(x)$. We define the \emph{balanced filtration} $\lambda_{\text{b}}(x)$ of $x$ to be the unique element $\lambda$ of $\qfilt(\cX,x)$ satisfying $\langle \lambda,\cF^{\max}\rangle\geq 1$ and such that for all filtrations $\gamma\in \qfilt(\cX,x\vert_{\overline k})$ with $\langle \gamma,\cF^{\max}\rangle\geq 1$ we have $\norm{\lambda}\leq \norm{\gamma}$.
\end{definition}
Note that we have an identification $\qfilt(\cX,x)=\qfilt(\cF,x)$ and that existence and uniqueness of the balanced filtration is guaranteed by \Cref{theorem: Kempf}.  The balanced filtration of $x$ is $0$ precisely when $x$ is closed in $\cF$.

In the case where $d(\cF)=d(\cX)$, from the fact that $\cX^{\max}\cap \cF=\cF^{\max}$ it follows that $\langle \lambda,\cF^{\max}\rangle=\langle \lambda, \cX^{\max} \rangle$. We have a $\Theta$-stratification $(\cS_c)_{c\in\Q_{\geq 0}}$ of $\Bl_{\cX^{\max}}\cX$, and thus a stratification of $\cX$ where the strata are $\cX^{\max}$ and the $\cS_c\setminus \cE$ with $c\geq 0$, where $\cE$ is the exceptional divisor. This can be seen as a stratification of $\cX$ by type of balanced filtration. Our goal is to refine this stratification by iterating the blowing-up process from the centres $\cZ_c$ of the strata $\cS_c$ (\Cref{theorem: sequential stratifications for good moduli stacks}). The iterated strata will naturally live inside a stack of sequential filtrations, defined in the following section.

\subsection{Stacks of sequential filtrations}
The main goal of this section is to define a stack $\Filt_{\Q^\infty}(\cX)$ of \emph{$\Q^\infty$-filtrations} (or \emph{sequential filtrations}) for a suitable algebraic stack $\cX$. Here, the symbol $\Q^\infty$ denotes the direct sum $\Q^{\oplus \N}$ of countably many copies of $\Q$ with lexicographic order. We start with a simpler version of the stack $\Filt_{\Q^\infty}(\cX)$.

\begin{definition}[Stack of {$\Q^n_{\lex}$}-filtrations]
Let $\cX$ be an algebraic stack over an algebraic space $B$, satisfying \Cref{assumption: basic assumptions}. We define the stack $\Filt_{\Q^n_{\lex}}(\cX)$ of \emph{$\Q^n$-filtrations of $\cX$ with lexicographic order} to be the limit of the diagram
\[\begin{tikzcd}
    {\Filt_\Q(\cX)} & {\Filt_\Q\Grad_\Q(\cX)} & \cdots & {\Filt_\Q\Grad_\Q^{n-1}(\cX)} \\
    & {\Grad_\Q(\cX)} & {\Grad^2_\Q(\cX)} & {\Grad^{n-1}_\Q(\cX).}
    \arrow["{\ev_1}", shift right=5, from=1-4, to=2-4]
    \arrow[shorten <=10pt, from=1-3, to=2-4]
    \arrow[shorten <=2pt, from=1-3, to=2-3]
    \arrow["\gr", from=1-2, to=2-3]
    \arrow["{\ev_1}"', from=1-2, to=2-2]
    \arrow["\gr", from=1-1, to=2-2]
\end{tikzcd}\]
\end{definition}
Thus $\Filt_{\Q^n_\lex}(\cX)$ is just the fibre product
\[\Filt_\Q(\cX)\underset{\Grad_\Q(\cX)}{\times} \Filt_\Q\Grad_\Q(\cX)\underset{\Grad^2_\Q(\cX)}{\times} \cdots\underset{\Grad^{n-1}_\Q(\cX)}{\times}  \Filt_\Q\Grad^{n-1}_\Q(\cX).\]

We define, for each $n\in \Z_{>0}$, a map $\Filt_{\Q^n_\lex}(\cX)\to \Filt_{\Q^{n+1}_\lex}(\cX)$ by

\[\begin{tikzcd}[column sep=2.25em,row sep=small]
    {\Filt_{\Q^n_\lex}(\cX)} & {\Filt_{\Q^{n+1}_\lex}(\cX)} \\
    {\Filt_{\Q^n_\lex}(\cX)\times_{\Grad^n_\Q(\cX)}\Grad^n_\Q(\cX)} & {\Filt_{\Q^n_\lex}(\cX)\times_{\Grad^n_\Q(\cX)}\Filt_\Q\Grad^n_\Q(\cX),}
    \arrow["{1\times o}", from=2-1, to=2-2]
    \arrow[Rightarrow, no head, from=1-2, to=2-2]
    \arrow[Rightarrow, no head, from=1-1, to=2-1]
\end{tikzcd}\]
where $o \colon \Grad_\Q^n(\cX)\to \Filt_\Q\Grad^n_\Q(\cX)$ is the “trivial filtration” map. If $\cX$ satisfies Assumption \ref{assumption: basic assumptions}, then so does $\Grad^n_\Q(\cX)$, so by \cite[Proposition~1.3.9]{_HalpernLeistner_Onthestructureofinstabilityinmodulitheory} and the argument in \cite[Proposition 1.3.11]{_HalpernLeistner_Onthestructureofinstabilityinmodulitheory}, the morphism $o$ is an open and closed immersion. Thus each of the maps $\Filt_{\Q^n_\lex}(\cX)\to \Filt_{\Q^{n+1}_\lex}(\cX)$ is an open immersion. We also have morphisms $\Grad^n_\Q(\cX)\to \Grad_\Q\Grad^n_\Q(\cX)=\Grad_\Q^{n+1}(\cX)$ given by the “trivial grading” maps, that are also open and closed immersions.

\begin{definition}[Stack of sequential filtrations]\label{definition: Q^infty filtrations and gradings}
Let $\cX$ be an algebraic stack over an algebraic space $B$, satisfying \Cref{assumption: basic assumptions}. We define the stack $\Filt_{\Q^\infty}(\cX)$ of \emph{$\Q^\infty$-filtrations} (or \emph{sequential filtrations}) of $\cX$ as the colimit
\[\Filt_{\Q^\infty}(\cX)=\colim_{n\in \Z_{>0}}\Filt_{\Q^n_\lex}(\cX)\]
in the 2-category $\cat{St_\fppf}$ of stacks for the fppf site of schemes.
Similarly, we define the stack $\Grad_{\Q^\infty}(\cX)$ of $\Q^\infty$-graded points of $\cX$ as 
\[\Grad_{\Q^\infty}(\cX)=\colim_{n\in\Z_{>0}} \Grad^n_\Q(\cX).\]
\end{definition}
As a direct application of \Cref{lemma: filtered colimit open immersions is algebraic}, we get:
\begin{proposition}\label{proposition: algebraicity of Q^infty filtrations}
Let $\cX$ be an algebraic stack over an algebraic space $B$, satisfying Assumption \ref{assumption: basic assumptions}. Then the stacks $\Filt_{\Q^\infty}(\cX)$ and $\Grad_{\Q^\infty}(\cX)$ are algebraic, naturally defined over $B$, and also satisfy Assumption \ref{assumption: basic assumptions}. \hfill \qedsymbol{}
\end{proposition}

\begin{remark}
Here we are regarding $\Q^\infty$ as having lexicographic order. It would be more precise to use the notation $\Filt_{\Q^\infty_\lex}(\cX)$ for what we called $\Filt_{\Q^\infty}(\cX)$, to distinguish it from a stack of $\Q^\infty$-filtrations with product order, which can also be defined. Since we will not use such a stack, our notation will not be problematic.
\end{remark}

For each $n$, we have an associated graded map $\gr\colon \Filt_{\Q^n_\lex}(\cX)\to \Grad^n_\Q(\cX)$, defined as the composition of the projection $\Filt_{\Q^n_\lex}(\cX)\to \Filt_\Q\Grad_\Q^{n-1}(\cX)$ and the associated graded map $\Filt_\Q\Grad_\Q^{n-1}(\cX)\to \Grad_\Q^n(\cX)$. These maps glue to a morphism $\gr\colon \Filt_{\Q^\infty}(\cX)\to \Grad_{\Q^\infty}(\cX)$. Similarly, we get an “evaluation at 1” map $\ev_1\colon \Filt_{\Q^\infty}(\cX)\to \cX$, a “forgetful” map $u\colon \Grad_{\Q^\infty}(\cX)\to \cX$, and a “split filtration” map $\sigma\colon \Grad_{\Q^\infty}(\cX)\to \Filt_{\Q^\infty}(\cX)$ as in \Cref{subsection: stacks of rational filtrations and graded points}. We will also consider the “trivial filtration” map $\cX\to \Filt_{\Q^\infty}(\cX)$, defined as the composition of the usual trivial filtration map $\cX\to \Filt_\Q(\cX)$ and the map $\Filt_\Q(\cX)=\Filt_{\Q^1_{\lex}}(\cX)\to \Filt_{\Q^\infty}(\cX)$ given by the colimit. It is an open and closed immersion. Similarly, we have a “trivial grading” map $\cX\to \Grad_{\Q^\infty}(\cX)$, and it is also an open and closed immersion. 

\begin{proposition}\label{proposition: representability of ev1}
Let $\cX$ be an algebraic stack over an algebraic space $B$, satisfying Assumption \ref{assumption: basic assumptions}. Then the “evaluation at 1” map $\ev_1\colon \Filt_{\Q^\infty}(\cX)\to \cX$ is representable and separated.
\end{proposition}
\begin{proof}
It is enough to prove that $\Filt_{\Q^n_\lex}(\cX)\to \cX$ is representable and separated for each $n$. There is a cartesian square
\[ \begin{tikzcd}
\Filt_{\Q^{n+1}_\lex}(\cX) \arrow[r,"a_n"]\arrow[d,swap,""]\arrow[dr, phantom, "\ulcorner", very near start] & \Filt_{\Q^n_\lex}(\cX) \arrow[d,""] \\
\Filt_\Q\Grad^n_\Q(\cX) \arrow[r,"b_n"]& \Grad^n_\Q(\cX)
\end{tikzcd}
\]
for each $n$, where $b_n$ is the “evaluation at 1” map. Thus by \cite[Proposition~1.1.13]{_HalpernLeistner_Onthestructureofinstabilityinmodulitheory}, the map $a_n$ is representable and separated, being a base change of $b_n$. Expressing $\ev_1\colon \Filt_{\Q^n_\lex}(\cX)\to \cX$ as a composition of the $a_n$, we get the result.
\end{proof}

\begin{remark}
One can define stacks $\Filt_{\Z^\infty}(\cX)$ and $\Grad_{\Z^\infty}(\cX)$ in a similar vein. The monoid $(\N,\cdot,1)$ acts on these stacks, and $\Filt_{\Q^\infty}(\cX)$ and $\Grad_{\Q^\infty}(\cX)$ are obtained from these by localising the action as in \Cref{definition: rational filtrations and gradings}.
\end{remark}
\begin{remark}
The formation of $\Filt_{\Q^\infty}(\cX)$ is functorial in $\cX$. If $f\colon \cX\to \cY$ is a morphism, then there is an obvious induced map $\Filt_{\Q^\infty}(f)\colon \Filt_{\Q^\infty}(\cX)\to \Filt_{\Q^\infty}(\cY)$.
\end{remark}

\begin{remark}[Stack of polynomial filtrations]
Our definition of the stack $\Filt_{\Q^\infty}(\cX)$ as a colimit of the stack $\Filt_{\Q^n_\text{lex}}(\cX)$ is justified by the fact that the poset $\Q^\infty$ can be written as the colimit $\Q^\infty=\colim_n \Q^n_{\text{lex}}$
where the maps are $\Q^n_\text{lex}\to \Q^{n+1}_\text{lex}\colon (a_0,\ldots,a_{n-1})\mapsto (a_0,\ldots,a_{n-1},0)$. Alternatively, we can consider the diagram where the maps are $\Q^n_\text{lex}\to \Q^{n+1}_\text{lex}\colon (a_0,\ldots,a_{n-1})\mapsto (0,a_0,\ldots,a_{n-1})$, and the colimit is now $\Q[t]$, the set of polynomials in one variable with rational coefficients, where $p\leq q$ if $p(n)\leq q(n)$ for $n>>0$. We define the stack $\Filt_{\Q[t]}(\cX)$ of \emph{polynomial filtrations} of $\cX$ to be the colimit of the corresponding diagram of open and closed immersions $\Filt_{\Q^n_\text{lex}}(\cX)\to \Filt_{\Q^{n+1}_\text{lex}}(\cX)$, which are given as the composition
\begin{align*}
\Filt_{\Q^n_{\lex}}(\cX) \to \Filt_{\Q^n_\lex}\left(\Grad_\Q(\cX)\right)= \Grad_\Q(\cX)\times_{\Grad_\Q(\cX)}\Filt_{\Q^n_\lex}\left(\Grad_\Q(\cX)\right) \\
 \xrightarrow{\sigma\times \id} \Filt_\Q(\cX)\times_{\Grad_\Q(\cX)}\Filt_{\Q^n_\lex}\left(\Grad_\Q(\cX)\right)=\Filt_{\Q^{n+1}_\lex}(\cX).
\end{align*}
 Everything we have said about $\Filt_{\Q^\infty}(\cX)$ applies also to $\Filt_{\Q[t]}(\cX)$ with similar arguments.
\end{remark}

The stack of sequential filtrations behaves well with respect to pullback with from an algebraic space.

\begin{proposition}\label{proposition: qinfinity Filt and base change}
Let $\cX\to B$ and $\cX'\to B'$ satisfy \Cref{assumption: basic assumptions}, and let
\[\begin{tikzcd}[ampersand replacement=\&]
    {\cX'} \& \cX \\
    {X'} \& X
    \arrow[from=2-1, to=2-2]
    \arrow[from=1-2, to=2-2]
    \arrow[from=1-1, to=2-1]
    \arrow[from=1-1, to=1-2]
    \arrow["\lrcorner"{anchor=center, pos=0.125}, draw=none, from=1-1, to=2-2]
\end{tikzcd}\]
be a cartesian square, with $X, X'$ algebraic spaces. Then 
\[\Filt_{\Q^\infty}(\cX')\cong \Filt_{\Q^\infty}(\cX)\times_{\ev_1,\cX} \cX'\cong \Filt_{\Q^\infty}(\cX)\times_X X'.\]
\end{proposition}
\begin{proof}
Since $\Filt_{\Q^\infty}(\cX')$ is an increasing union of the stacks $\Filt_{\Q^n_{\text{lex}}}(\cX')$, it is enough to show the analogue claim for these stacks. For $n=1$, this is \Cref{proposition: Filt Grad and base change}. For $n>1$, we have
\begin{align*}
\cX'\times_{\cX,\ev_1}\Filt_{\Q^n_\lex}(\cX)=\cX'\times_{\cX,\ev_1}\Filt_{\Q^{n-1}_\lex}(\cX)\times_{\gr,\Grad^{n-1}_\Q(\cX),\ev_1} \Filt_\Q(\Grad^{n-1}_\Q(\cX))\\
=\Filt_{\Q^{n-1}_\lex}(\cX')\times_{\gr,\Grad_\Q^{n-1}(\cX')}\Grad_\Q^{n-1}(\cX')\times_{\Grad_\Q^{n-1}(\cX),\ev_1}\Filt_\Q\left(\Grad_\Q^{n-1}(\cX)\right)\\
=\Filt_{\Q^{n-1}_\lex}(\cX')\times_{\gr,\Grad_\Q^{n-1}(\cX'),\ev_1}\Filt_\Q\left(\Grad_\Q^{n-1}(\cX')\right)=\Filt_{\Q^{n-1}_\lex}(\cX')
\end{align*}
again by \Cref{proposition: Filt Grad and base change}.
\end{proof}
\begin{proposition}\label{proposition: infFiltclosed immersion}
Let $\cX$ be an algebraic stack defined over an algebraic space $B$, satisfying \Cref{assumption: basic assumptions}, and let $\cX'\to \cX$ be a closed immersion. Then
\[\Filt_{\Q^\infty}(\cX')\cong \Filt_{\Q^\infty}(\cX)\times_{\ev_1,\cX}\cX'.\]
\end{proposition}
\begin{proof}
The statement follows in the same way as \Cref{proposition: qinfinity Filt and base change}, but using \Cref{proposition: Filt Grad pullback closed immersion} instead of \Cref{proposition: Filt Grad and base change}.
\end{proof}

\begin{proposition}\label{proposition: iterative cartesian diagram sequential filtrations}
Let $\cX$ be an algebraic stack over an algebraic space $B$, satisfying \Cref{assumption: basic assumptions}. Then there is a cartesian diagram
\[\begin{tikzcd}
    {\Filt_{\Q^\infty}(\cX)} & {\Filt_{\Q^\infty}\Grad_\Q(\cX)} \\
    {\Filt_\Q(\cX)} & {\Grad_\Q(\cX)}
    \arrow[""{name=0, anchor=center, inner sep=0}, "\gr", from=2-1, to=2-2]
    \arrow["{\ev_1}", from=1-2, to=2-2]
    \arrow[from=1-1, to=2-1]
    \arrow[from=1-1, to=1-2]
    \arrow["\lrcorner"{anchor=center, pos=0.125}, draw=none, from=1-1, to=0]
\end{tikzcd}\]
\end{proposition}
\begin{proof}
The claim follows from cartesianity of the diagram
\[\begin{tikzcd}
    {\Filt_{\Q^{n+1}_{\text{lex}}}(\cX)} & {\Filt_{\Q^n_{\text{lex}}}\Grad_\Q(\cX)} \\
    {\Filt_\Q(\cX)} & {\Grad_\Q(\cX)}
    \arrow[""{name=0, anchor=center, inner sep=0}, "\gr", from=2-1, to=2-2]
    \arrow["{\ev_1}", from=1-2, to=2-2]
    \arrow[from=1-1, to=2-1]
    \arrow[from=1-1, to=1-2]
    \arrow["\lrcorner"{anchor=center, pos=0.125}, draw=none, from=1-1, to=0]
\end{tikzcd}\]
by taking the colimit when $n$ tends to $\infty$.
\end{proof}

We define $\Q^\infty$-flag spaces in analogy with \Cref{definition: flag spaces}.

\begin{definition}[{$\Q^\infty$-Flag spaces}]\label{definition: Q^infty flag space}
Let $\cX$ be an algebraic stack over an algebraic space $B$, satisfying Assumption \ref{assumption: basic assumptions}, and let $x\colon T\to \cX$ be a scheme-valued point. We define the \emph{space of $\Q^\infty$-flags $\Flag_{\Q^\infty}(x)$  of $x$} as the fibre product 
\[ \begin{tikzcd}
\Flag_{\Q^\infty}(x) \arrow[r,""]\arrow[d,swap,""]\arrow[dr, phantom, "\ulcorner", very near start] & T \arrow[d,"x"] \\
\Filt_{\Q^\infty}(\cX) \arrow[r,"\ev_1"]& \cX,
\end{tikzcd}
\]
which is, by \Cref{proposition: representability of ev1}, a separated algebraic space over $T$ locally of finite presentation.
\end{definition}

In the case of a field-valued point $x$, it is reasonable to talk about $\Q^\infty$-filtrations.

\begin{definition}[{$\Q^\infty$}-filtrations of a point]\label{definition:Q^infty filtrations of a point}
Let $\cX$ be an algebraic stack over an algebraic space $B$, satisfying Assumption \ref{assumption: basic assumptions}. Let $k$ be a field and let $x\colon \Spec k\to \cX$ point. We define the set $\qinfilt(\cX,x)$ of \emph{$\Q^\infty$-filtrations} (or \emph{sequential filtrations}) of $x$ to be $\Q^\infty\text{-}\Filt(\cX,x)\coloneqq \Flag_{\Q^\infty}(x)(k)$, the set of $k$-points of the space of $\Q^\infty$-flags of $x$.
\end{definition}

\begin{remark}\label{remark: description set of infty filtrations}
The set $\qinfilt(\cX,x)$ can be described as follows. An element $\lambda$ of $\qinfilt(\cX,x)$ is uniquely determined by a sequence $(\lambda_n)_{n\in\N}$ where 
\begin{enumerate}
\item $\lambda_0\in \qfilt(\cX,x)$, 
\item $\lambda_n\in\qfilt(\Grad^n_\Q(\cX),\gr\lambda_{n-1})$ for $n\geq 1$, and 
\item $\lambda_n= 0$ for $n>>0$.
\end{enumerate}
The trivial $\Q^\infty$-filtration corresponds to the sequence where $\lambda_n=0$ for all $n$. In general, 
\[N\coloneqq\min\{n\st \lambda_i=0,\ \forall i\geq n\}\]
is the minimal natural number such that $\lambda$ is in $\Filt_{\Q^N_{\lex}}(\cX)\subset \Filt_{\Q^\infty}(\cX)$. The description follows at once from the definition of $\Filt_{\Q^N_{\lex}}(\cX)$ as a fibre product.
\end{remark}

\begin{remark}\label{remark: bijection between sets of sequential filtrations}
If $f\colon \cX'\to \cX$ is either a closed immersion or a base change from a map of algebraic spaces, and if $x\in \cX'(k)$ is a field-valued point, then, by \Cref{proposition: infFiltclosed immersion,proposition: qinfinity Filt and base change}, we have a canonical bijection $\qfilt(\cX',x)\cong\qfilt(\cX,f(x))$. We will use this fact throughout.
\end{remark}

\begin{remark}[Sequential filtrations on quotient stacks]\label{remark: description sequential filtrations quotient stacks}
Let $k$ be a field, and consider a quotient stack $\cX=X/G$ where $X$ is a separated scheme of finite type over $k$ and $G$ is a linear algebraic group over $k$. Let $x\in X(k)$ be a $k$-point and call also $x\in \cX(k)$ its image in $\cX$. From \Cref{remark: filtrations of a quotient stack,remark: description set of infty filtrations}, it follows that we have an identification of the set $\qinfilt(\cX,x)$ of sequential filtrations of $x$ with the set of sequences $(\lambda_n)_{n\in \N}$ where
\begin{enumerate}
\item each $\lambda_n\in \Gamma^\Q(G)$; 
\item for all $n,m\in \N$, $\lambda_n$ and $\lambda_m$ commute; 
\item $\lambda_n=0$ for $n>>0$; 
\item for every $n\in \Z_{>0}$, the iterated limit 
\[\lim_{t_n\to 0}\lambda_n(t_n)\left(\lim_{t_{n-1}\to 0} \lambda_{n-1}(t_{n-1})\left(\cdots \lim_{t_0\to 0}\lambda_0(t_0)x\right)\cdots\right)\] 
exists in $X$;
\end{enumerate}
subject to the equivalence relation that identifies $(\lambda_n)_{n\in \N}\sim (\lambda'_n)_{n\in \N}$ if there are $g_n\in P_{L(\lambda_0,\cdots,\lambda_n)}(\lambda_n)$ such that $(\lambda_n)^{g_0g_1\cdots g_n}=\lambda'_n$ for all $n\in \N$. To see this, note that we can explicitly describe components of $\Grad_\Q^n(\cX)$ as $X^{\lambda_0,\cdots,\lambda_n, 0}/L(\lambda_0,\cdots,\lambda_n)$, where $\lambda_0,\cdots,\lambda_n$ are commuting rational cocharacters of $G$, $X^{\lambda_0,\cdots,\lambda_n, 0}$ denotes the fixed-point locus by $\lambda_0,\cdots,\lambda_n$ and $L(\lambda_0,\cdots,\lambda_n)$ is the centraliser of $\lambda_0,\cdots,\lambda_n$ inside $G$, and we can apply \Cref{remark: filtrations of a quotient stack} to these quotient stacks.
\end{remark}

\subsection{Sequential stratifications}
In this section, we give the definition of \emph{sequential stratification} and we study some pullback and pushforward operations for sequential strata. All algebraic stacks are defined over an algebraic space $B$ and are assumed to satisfy Assumption \ref{assumption: basic assumptions}.

\begin{definition}[Sequential stratification]\label{definition: sequential stratification}
Let $\cX$ be an algebraic stack and let $\Gamma$ be a partially ordered set. A \emph{sequential stratification} (or \emph{$\Q^\infty$-stratification}) of $\cX$ indexed by $\Gamma$ is a family $(\cS_i)_{i\in \Gamma}$ of locally closed substacks of $\Filt_{\Q^\infty}(\cX)$ such that:
\begin{enumerate}
\item Each composition $r_i\colon \cS_i\to \Filt_{\Q^\infty}(\cX)\to \cX$ is a locally closed immersion.
\item The $r_i(\abs{\cS_i})$ are pairwise disjoint and cover $\abs{\cX}$.
\item For each $i\in \Gamma$, the union $\bigcup_{j\leq i}r_j(\abs{\cS_j})$ is open in $\abs{\cX}$.
\end{enumerate}
\end{definition}

\begin{remark}[Sequential and polynomial {$\Theta$}-stratifications]
This definition is not quite the analogue of the notion of $\Theta$-stratification for $\Q^\infty$-filtrations. That would require in addition that each $\cS_i$ is open in $\Filt_{\Q^\infty}(\cX)$ and is the preimage of an open substack of $\Grad_{\Q^\infty}(\cX)$, but these conditions do not hold for the balancing stratification, that we will construct later (\Cref{definition: balancing stratification}). We may call this stronger notion \emph{sequential $\Theta$-stratifications}. Similarly, if we use the stack $\Filt_{\Q[t]}(\cX)$ of polynomial filtrations instead of $\Filt_{\Q^\infty}(\cX)$, we get a notion of \emph{polynomial $\Theta$-stratification}. The stratification of the stack of pure coherent sheaves on a polarised projective scheme over a noetherian base by polynomial Harder-Narasimhan filtration \cite{_Nitsure_SchematicHarderNarasimhanstratification} should be a polynomial $\Theta$-stratification. Another example should be given by the stratifications for moduli spaces of principal $\rho$-sheaves considered in \cite{https://doi.org/10.48550/arxiv.2107.03918}.
\end{remark}

\begin{definition}[Pulling back sequential stratifications]\label{definition: pullback of sequential stratifications}
Let $f\colon \cX'\to \cX$ be a morphism of algebraic stacks such that
\[\begin{tikzcd}[ampersand replacement=\&]
    {\Filt_{\Q^\infty}(\cX')} \& {\Filt_{\Q^\infty}(\cX)} \\
    {\cX'} \& \cX
    \arrow["f", from=2-1, to=2-2]
    \arrow["{\ev_1}", from=1-2, to=2-2]
    \arrow["{\ev_1}"', from=1-1, to=2-1]
    \arrow[from=1-1, to=1-2]
    \arrow["\ulcorner"{anchor=center, pos=0.125}, draw=none, from=1-1, to=2-2]
\end{tikzcd}\]
is cartesian (for example a closed immersion or a base change from a map of algebraic spaces, see \Cref{proposition: infFiltclosed immersion,proposition: qinfinity Filt and base change}). Let $(\cS_i)_{i\in \Gamma}$ be a sequential stratification of $\cX$. For each $i\in \Gamma$, set $f^*\cS_i\coloneqq\Filt_{\Q^\infty}(\cX')\times_{\Filt_{\Q^\infty}(\cX)} \cS_i$, which is a locally closed substack of $\Filt_{\Q^\infty}(\cX')$. Then $(f^*\cS_i)_{i\in \Gamma}$ is a sequential stratification of $\cX'$, called the \emph{pulled back} sequential stratification.
\end{definition}

To see that $(f^*\cS_i)_{i\in \Gamma}$ is indeed a sequential stratification, just note that we have $f^*\cS_i=\cX'\times_\cX \cS_i$ for all $i\in \Gamma$. We refer to the $\cS_i$ as \emph{sequential strata}. More generally, we define:

\begin{definition}\label{definition: sequential stratum}
A \emph{sequential stratum} $\cS$ of $\cX$ is a locally closed substack of $\Filt_{\Q^\infty}(\cX)$ such that the composition
\[\cS\to \Filt_{\Q^\infty}(\cX)\to \cX\]
is a locally closed immersion. We refer to the morphism $\cS\to \cX$ as the \emph{structure map}.
\end{definition}

\begin{remark}\label{remark: enough for a sequential stratum that the map to X is a locally closed immersion}
If $a\colon \cS\to \Filt_{\Q^\infty}(\cX)$ is a morphism such that the composition 
\[\cS\to \Filt_{\Q^\infty}(\cX)\to \cX\]
 is a locally closed immersion, then $a$ is a locally closed immersion as well, since $\ev_1\colon \Filt_{\Q^\infty}(\cX)\to \cX$ is representable and separated, so its diagonal is a closed immersion.
\end{remark}

It will be useful to have a few ways of constructing sequential strata from given ones.

\begin{definition}[Pushforward along a locally closed immersion]\label{definition: pushforward of sequential stratum along locally closed immersion}
If $\iota\colon \cX\to \cY$ is a locally closed immersion and $\cS$ is a sequential stratum of $\cX$, then we define a sequential stratum $\iota_*\cS$ as follows. As a stack, $\iota_*\cS=\cS$, and the structure map is the composition
\[\cS\to \Filt_{\Q^\infty}(\cX)\to \Filt_{\Q^\infty}(\cY),\]
which is a locally closed immersion because $\Filt_{\Q^\infty}(\iota)$ is\footnote{The fact that $\Filt_\Q$ and $\Grad_\Q$ preserve closed immersions \cite[Proposition 1.3.1]{_HalpernLeistner_Onthestructureofinstabilityinmodulitheory} easily implies that $\Filt_{\Q^\infty}(\iota)$ is a closed immersion.}.
The map $\cS\to \cY$ factors as $\cS\to \cX\to \cY$, so it is a locally closed immersion.
\end{definition}

\begin{definition}[Induction of a sequential stratum from centre of a {$\Theta$}-stratum]\label{definition: Induction of a sequential stratum from centre of a Theta-stratum}
Suppose $\cS$ is a locally closed $\Theta$-stratum of an algebraic stack $\cX$, with centre $\cZ$ (\Cref{definition: locally closed Theta-stratum}), and let $\cS'$ be a sequential stratum of $\cZ$. We define a sequential stratum $\ind_{\cZ}^\cX(\cS')$ of $\cX$, as follows. As a stack, $\ind_\cZ^\cX(\cS')$ is the pullback
\[\begin{tikzcd}
    {\ind_\cZ^\cX(\cS')} & {\cS'} \\
    \cS & {\cZ.}
    \arrow["\gr", from=2-1, to=2-2]
    \arrow[from=1-2, to=2-2]
    \arrow[from=1-1, to=1-2]
    \arrow[from=1-1, to=2-1]
    \arrow["\ulcorner"{anchor=center, pos=0.125}, draw=none, from=1-1, to=2-2]
\end{tikzcd}\]
The structure morphism $a$ is obtained by pulling back the square
\[\begin{tikzcd}
    {\cS'} & {\Filt_{\Q^\infty}\Grad_\Q(\cX)} \\
    \cZ & {\Grad_\Q(\cX)}
    \arrow[from=2-1, to=2-2]
    \arrow[from=1-2, to=2-2]
    \arrow[from=1-1, to=2-1]
    \arrow[from=1-1, to=1-2]
\end{tikzcd}\]
along $\Filt_\Q(\cX)\xrightarrow{\gr}\Grad_\Q(\cX)$, obtaining a square of the form
\[\begin{tikzcd}
    {\ind_\cZ^\cX(\cS')} & {\Filt_{\Q^\infty}(\cX)} \\
    \cS & {\Filt_\Q(\cX)}
    \arrow["a", from=1-1, to=1-2]
    \arrow[from=1-2, to=2-2]
    \arrow[from=2-1, to=2-2]
    \arrow[from=1-1, to=2-1]
\end{tikzcd}\]
by \Cref{proposition: iterative cartesian diagram sequential filtrations}.
The induced morphism $\ind_{\cZ}^\cX(\cS')\to \cX$ is the composition of the locally closed immersions $\ind_{\cZ}^\cX(\cS')\to \cS$ and $\cS\to\cX$ and it is thus also a locally closed immersion.
\end{definition}

\begin{definition}[Pushforward along a blow-up]\label{definition: pushforward sequential stratum along a blow-up}
Let $p\colon \cY\to \cX$ be a blow-up with exceptional divisor $\cE\subset \cY$. Let $\cS$ be a sequential stratum of $\cY$. We define a sequential stratum $p_*(\cS\setminus \cE)$ of $\cX$ as follows. If $r\colon \cS\to \cY$ is the locally closed immersion, then we set $p_*(\cS\setminus \cE)=\cS\setminus r^{-1}(\cE)$ as a stack. The structure map is the composition
\[\cS\setminus r^{-1}(\cE)\to \cS\to \Filt_{\Q^\infty}(\cY)\to \Filt_{\Q^\infty}(\cX).\]
We have a diagram
\[\begin{tikzcd}
    {\cS\setminus r^{-1}(\cE)} & {\Filt_{\Q^\infty}(\cY)} & {\Filt_{\Q^\infty}(\cX)} \\
    \cY\setminus\cE & \cY & \cX.
    \arrow[from=1-3, to=2-3]
    \arrow["p"', from=2-2, to=2-3]
    \arrow[from=1-2, to=1-3]
    \arrow[from=1-1, to=1-2]
    \arrow["b"', from=2-1, to=2-2]
    \arrow[from=1-2, to=2-2]
    \arrow["a"', from=1-1, to=2-1]
\end{tikzcd}\]
The map $a$ is a locally closed immersion, and $p\circ b$ is an open immersion. Thus $\cS\setminus r^{-1}(\cE)\to \cX$ is a locally closed immersion. By Remark \ref{remark: enough for a sequential stratum that the map to X is a locally closed immersion}, the structure map is also a locally closed immersion.
\end{definition}

\subsection{Central rank\texorpdfstring{ and $B\G_m^n$-actions}{}}
Let $\cX$ be an algebraic stack over an algebraic space $B$, satisfying Assumption \ref{assumption: basic assumptions}, and further that $\cX$ is noetherian and has affine diagonal.

\begin{definition}[Central rank]\label{definition: central rank}
We define the \emph{central rank} $z(\cX)\in \N$ of $\cX$ to be the biggest $n\in \N$ such that there is a union $\cZ$ of nondegenerate components of $\Grad^n(\cX)$ (Definition \ref{definition: nondegenerate graded point}) such that the composition $\cZ\to \Grad^n(\cX)\xrightarrow{u} \cX$ is an isomorphism. For the empty stack we define $z(\varnothing)=\infty$.
\end{definition}

The stack $B\G_m^n$ is a group stack. It turns out that the data of a $B\G_m^n$ action on a stack $\cX$ satisfying Assumption \ref{assumption: basic assumptions} is equivalent to the data of a section $s\colon \cX \to \Grad^n(\cX)$ of $u\colon \Grad^n(\cX)\to \cX$ such that $s$ is a closed and open immersion. This is \cite[Corollary 1.4.2.1]{_HalpernLeistner_DerivedThetastratificationsandtheDequivalenceconjecture} in the case of the group stack $\A^1/\G_m$, but the same proof works for $B\G_m^n$.
 We say that the $B\G_m^n$-action is \emph{nondegenerate} if all components of $s(\cX)$ are nondegenerate. If $m\colon B\G_m^n\times \cX\to \cX$ is the action map, then the action is nondegenerate if, for all $x\in \cX(k)$, the homomorphism $\G_{m,k}^n\to \Aut(x)$ induced by $m$ has finite kernel.

 Using the stack $\Grad_\Q^n(\cX)$ instead of $\Grad^n(\cX)$ we can talk about rational $B\G_{m}^n$-actions.

 \begin{definition}\label{definition: rational BG_m^n-action}
A \emph{rational $B\G_{m}^n$-action} on $\cX$ is a section $s\colon \cX\to \Grad_\Q^n(\cX)$ of \[u\colon \Grad^n_\Q(\cX)\to \cX\] that is a closed and open immersion.
 \end{definition}

\begin{remark}
It is tacitly understood that an isomorphism $u\circ s \sim \id_\cX$ is part of the data of a section $s$.
\end{remark}

\begin{lemma}\label{lemma: lifting BG_m^n actions along blow-ups}
Suppose that $\cX$ has a good moduli space $\pi\colon \cX\to X$, and let $p\colon \cY\to \cX$ be a blow-up. Then any (rational) $B\G_m^n$-action on $\cX$ lifts canonically to $\cY$.
\end{lemma}
\begin{proof}
We prove the lemma for integral actions, the proof for rational actions being identical after replacing $\Grad$ with $\Grad_\Q$. The $B\G_m^n$-action corresponds to a section $s\colon \cX\to \Grad^n(\cX)$ of $\Grad^n(\cX)\to \cX$ that is an open and closed immersion. With this language, what we want to show is that the preimage $\cY'$ of $s(\cX)$ along $\Grad^n(\cY)\to \Grad^n(\cX)$ is a closed and open substack of $\Grad^n(\cY)$ such that the composition $\cY'\to \Grad^n(\cY)\to \cY$ is an isomorphism.

If $\cX$ is of the form $\cX=\Spec A/\GL_l$, then 
\[\Grad^n(\cX)=\bigsqcup_{\lambda\in \Hom(\G_m^n,T)/W} (\Spec A)^{\lambda,0}/L(\lambda),\]
 where $T$ is the standard maximal torus of $\GL_l$ and $W$ is the symmetric group of degree $l$ \cite[Theorem 1.4.7]{_HalpernLeistner_Onthestructureofinstabilityinmodulitheory}. Thus $\cX$ is isomorphic to a union of connected components of one of the stacks $(\Spec A)^{\lambda,0}/L(\lambda)$, with $\lambda\colon \G_m^n\to \GL_l$ having finite kernel. Any blow-up $p\colon \cY\to \cX$ is then of the form $Y/L(\lambda)$ with $\lambda(\G_m^l)$ acting trivially on $Y$. This proves the lemma in the case $\cX=\Spec A/\GL_l$.

For general $\cX$, the claim can be checked étale locally on $\cX$, since for $\cX'\to \cX$ representable and étale, the square
\[\begin{tikzcd}[ampersand replacement=\&]
    {\Grad^n(\cX')} \& {\Grad^n(\cX)} \\
    {\cX'} \& \cX
    \arrow[from=1-2, to=2-2]
    \arrow[from=2-1, to=2-2]
    \arrow[from=1-1, to=2-1]
    \arrow[from=1-1, to=1-2]
    \arrow["\ulcorner"{anchor=center, pos=0.125}, draw=none, from=1-1, to=2-2]
\end{tikzcd}\]
is cartesian \cite[Corollary 1.1.7]{_HalpernLeistner_Onthestructureofinstabilityinmodulitheory}, and since blow-ups commute with flat base change. We can cover $\cX$ by representable étale neighbourhoods of the form $\Spec A/\GL_l$ (\Cref{theorem: local structure}), which proves the lemma for general $\cX$.
\end{proof}
\begin{remark}
It should not be essential that $\cX$ has a good moduli space for \Cref{lemma: lifting BG_m^n actions along blow-ups}, but this is all we will need.
\end{remark}

\begin{lemma}\label{lemma: centres of unstable strata have bigger central rank}
Let $\cX$ be a normed noetherian good moduli stack with affine diagonal. Let $p\colon \cY\to \cX$ be a blow-up of $\cX$ at some closed substack, and let $\cE$ be the exceptional divisor. The stack $\cY$ is endowed with the $\Theta$-stratification induced by the norm on $\cX$ and the $p$-ample line bundle $\cO_\cY(-\cE)$ (Theorem \ref{theorem: theta stratification proper over gms} and \Cref{example: stratification on a blow-up}). Let $\cZ_c$ be the centre of some unstable stratum of $\cY$. Then $z(\cX)<z(\cZ_c)$.
\end{lemma}
\begin{proof}
Let $n=z(\cX)$. There is a nondegenerate $B\G_m^n$ action on $\cX$ that, since $p$ is a blow-up, lifts canonically to $\cY$ by \Cref{lemma: lifting BG_m^n actions along blow-ups}. The $B\G_m^n$ action gives a closed and open immersion $\cY\to \Grad^n(\cY)$.

Since $\cZ_c$ is the centre of a stratum, it comes equipped with a rational $B\G_m$-action inherited from a closed and open immersion $\cZ_c\to \Grad_\Q(\cY)$. Scaling up, we get a closed and open immersion $\cZ_c\to \Grad(\cY)$ and thus an integral $B\G_m$-action on $\cZ_c$. Applying $\Grad$ to the closed and open immersion $\cY\to \Grad^n(\cY)$ and composing with $\cZ_c\to \Grad(\cY)$, we get a closed and open immersion $\cZ_c\to \Grad^{n+1}(\cY)$, which gives a $B\G_m\times B\G_m^n$-action on $\cZ_c$.

Let $x\in \cZ_c(k)$ be a $k$-point for some field $k$. The $B\G_m\times B\G_m^n$-action provides cocharacters $\lambda,\beta_1,\ldots,\beta_n$ of the centre $Z(\Aut(x))$. Since the $\beta_1,\ldots,\beta_n$ come from $\cX$, we have $\langle \beta_i,\cO_\cY(-\cE)\rangle=0$, where $\cE$ is the exceptional divisor, while since $\cZ_c$ is the centre of an \emph{unstable} stratum, we have $\langle \lambda,\cO_\cY(-\cE)\rangle >0$. Therefore $\lambda(\G_m)$ is not contained in the image of $(\beta_1,\ldots,\beta_n)\colon \G_m^n\to \Aut(x)$, and thus $(\lambda,\beta_1,\ldots,\beta_n)\colon \G_m^{n+1}\to \Aut(x)$ has finite kernel. Therefore $z(\cZ_c)\geq n+1$, as desired.
\end{proof}

\subsection{The balancing stratification and the iterated balanced filtration}

We now get to the main construction of the paper, namely the \emph{balancing stratification} for normed good moduli stacks (\Cref{theorem: sequential stratifications for good moduli stacks}), and the canonical sequential filtration it defines for every point (the \emph{iterated balanced filtration}, \Cref{definition: iterated balanced filtration}). We also show that the balancing stratification is preserved under certain pullbacks (\Cref{proposition: compatibility of balancing stratification with pullback}). Let us first introduce the indexing poset labelling the stratification.

\begin{definition}[Indexing poset for the balancing stratification]\label{definition: totally ordered poset bGamma}
We define a totally ordered set $\bGamma$ as follows. As a set, $\bGamma$ consists of the sequences
\[((d_0,c_0),(d_1,c_1),\ldots,(d_n,c_n))\]
with
\begin{enumerate}
\item $n\in \N$,
\item $d_0\geq d_1 \geq \cdots \geq d_n$ in $\N$,
\item $d_i\geq i$ for all $0\leq i \leq n$,
\item $c_0,\ldots,c_{n-1}\in \Q_{>0}$,
\item $c_n=\infty$.
\end{enumerate}
As a poset, we write
\[((d_0,c_0),\ldots,(d_n,c_n))<((d_0',c_0'),\ldots,(d_m',c_m'))\]
if there is $0\leq i\leq \min(n,m)$ such that $(d_j,c_j)=(d'_j,c'_j)$ for $j<i$ and either $d_i<d'_i$ or $d_i=d_i'$ and $c_i<c_i'$. This makes $\bGamma$ a totally ordered set.
\end{definition}
If $\alpha\in \bGamma$, we use the notation
\[\alpha=((d^\alpha_0,c^\alpha_0),\ldots,(d^\alpha_{n(\alpha)},c^\alpha_{n(\alpha)})).\]

\begin{theorem}[Existence and characterisation of the balancing stratification]\label{theorem: sequential stratifications for good moduli stacks}
There is a unique way of assigning, to every normed noetherian good moduli stack $\cX$ with affine diagonal, a sequential stratification $(\cS^\cX_\alpha)_{\alpha\in \bGamma}$ of $\cX$ indexed by $\bGamma$ in such a way that the following properties are satisfied for every such $\cX$:
\begin{enumerate}
\item \label{item 1 stratification theorem}The stratum indexed by $(d(\cX),\infty)$ is $\cS^\cX_{(d(\cX),\infty)}=\cX^{\max}$, with structure map
\[\cX^{\max}\to \cX\xrightarrow{\text{}} \Filt_{\Q^\infty}(\cX),\]
where the second arrow is the “trivial filtration” morphism.
\item \label{item 2 stratification theorem}Let $\pi\colon \cX\to X$ be the good moduli space and let $\cU=\cX\setminus \pi^{-1}\pi(\cX^{\max})$. Denote $j\colon \cU\to \cX$ the open immersion. Then for all $\alpha\in\bGamma$ with $d_0^\alpha<d(\cX)$ we have $\cS_\alpha^\cX=j_*(\cS_\alpha^\cU)$.
\item \label{item 3 stratification theorem} Let $p\colon \cY=\Bl_{\cX^{\max}} \cX \to \cX$ be the blow-up of $\cX$ along $\cX^{\max}$ and let $\cE$ be the exceptional divisor. Let $\cZ_c$ be the centres of the $\Theta$-stratification on $\cY$ induced by the norm on $\cX$ and the line bundle $\cO_\cY(-\cE)$ (\Cref{example: stratification on a blow-up}). The $\cZ_c$ are endowed with the induced norm and are thus also normed noetherian good moduli stacks with affine diagonal by Theorem \ref{theorem: theta stratification proper over gms}. Then, for all $c\in \Q_{>0}$ and for all $\alpha\in \bGamma$ such that the concatenation $((d(\cX),c),\alpha)$ of $(d(\cX),c)$ and $\alpha$ also belongs to $\bGamma$, we have the equality
\[\cS_{((d(\cX),c),\alpha)}^\cX=p_*\left(\ind_{\cZ_c}^\cY(\cS^{\cZ_c}_\alpha)\setminus \cE\right).\]
\end{enumerate}
Moreover, for every $\alpha\in \bGamma$ such that $\cS_\alpha^\cX$ is nonempty, we have:
\begin{enumerate}
\setcounter{enumi}{3}
\item \label{item 4 stratification theorem}$d_0^\alpha\leq d(\cX)$; and
\item \label{item 5 stratification theorem}$d_i^\alpha \geq i+z(\cX)$, for all $0\leq i \leq n(\alpha)$.
\end{enumerate}
\end{theorem}
\begin{proof}
First we note that if $\cX$ is empty, then the only stratification is given by $\cS_\alpha^\cX=\varnothing$ for all $\alpha$, and it satisfies the required properties.

For a stack $\cX$ as in the statement of the theorem, define the number $N(\cX)=d(\cX)-z(\cX)$. We will use induction on $N(\cX)$.

Assume $\cX\neq \varnothing$ and let $\cU$ be as in \ref{item 2 stratification theorem}. Then clearly $d(\cU)<d(\cX)$ and $z(\cU)\geq z(\cX)$, so $N(\cU)<N(\cX)$. If $\cZ_c$ is as in \ref{item 3 stratification theorem}, with $c>0$, then $d(\cZ_c)\leq d(\cX)$, because $\cZ_c\to \cX$ is representable, and $z(\cZ_c)>z(\cX)$ by Lemma \ref{lemma: centres of unstable strata have bigger central rank}. Thus $N(\cZ_c)<N(\cX)$. Therefore the statement of the theorem makes sense if we fix an $N\in \N$ and we restrict to the class of normed noetherian good moduli stacks $\cX$ with affine diagonal and with $N(\cX)\leq N$. We prove the theorem for this class of stacks by induction on $N$.

If $N=0$, and $N(\cX)=0$, then the identity component of every stabiliser group of a point of $\cX$ is a split torus of dimension $d(\cX)$. Therefore $\abs{\cX^{\max}}=\abs{\cX}$ and $S^\cX_{(d(\cX),\infty)}=\cX^{\max}$ is the only nonempty stratum. This gives the desired sequential stratification.

Fix $N>0$ and assume the theorem is true for $N-1$. For $\cX$ with $N(\cX)=N$, define $\cS^\cX_\alpha$ as in the statement of the theorem when $d_0^\alpha\leq d(\cX)$, which makes sense because $N(\cU)<N$ and $N(\cZ_c)<N$. Define $\cS^\cX_\alpha=\varnothing$ otherwise. We show that $(\cS_\alpha^\cX)_{\alpha\in \bGamma}$ is a sequential stratification of $\cX$.

Denote $r_\alpha\colon \cS_\alpha^\cX\to \cX$ the induced locally closed immersions. First we show that the $r_\alpha\left(\abs{\cS_\alpha^\cX}\right)$ are pairwise disjoint and cover $\cX$. For $p\in \abs{\cX}$, there is a unique $p'\in \overline{\{x\}}$ that is closed in $\pi^{-1}\pi(p)$ \cite[Proposition 9.1]{_Alper_GoodmodulispacesforArtinstacks}. One and only one of the following situations takes place:
\begin{enumerate}
\item The dimension $\dim \Aut(p')<d(\cX)$. In this case $p\in \abs{\cU}$ and it is contained in exactly one of the $j_*\cS_\alpha^\cU=\cS_\alpha^\cX$ with $d_0^\cX<d(\cX)$, by induction hypothesis.
\item We have $\dim \Aut(p')=d(\cX)$ and $p=p'$. Then $p\in \abs{\cX^{\max}}=\abs{\cS^\cX_{d(\cX),\infty)}}$ and this is the only stratum containing $p$.
\item Again, $\dim \Aut(p')=d(\cX)$, but this time $p\neq p'$. In this case, $p\in \abs{\cX}\setminus \left( \abs{\cU}\cup \abs{\cX^{\max}}\right)$ and there is a unique point $q\in \abs{\cY}\setminus \abs{\cE}$ mapping to $p$. The point $q$ lies in a unique stratum $\cS_c$ ($c\in \Q_{\geq 0}$) of the $\Theta$-stratification of $\cY$. There exists a map $\lambda\colon \Theta_k\to \cX$ with $\lambda(0)=p'$ and $\lambda(1)=p$ for some field $k$ \cite[Lemma 3.24]{_Alper_Existenceofmodulispacesforalgebraicstacks}, and Kempf's number $\langle \lambda,\cX^{\max}\rangle>0$ is positive because $\lambda(0)\in \abs{\cX^{\max}}$ and $\lambda(1)\notin \abs{\cX^{\max}}$ (Definition \ref{definition: Kempfs intersection number}). The filtration $\lambda$ lifts uniquely to $\lambda\colon \Theta_k\to \cY$ (Proposition \ref{proposition: proper map induces bijection of set of filtrations}), and $\langle \lambda, \cO_\cY(-\cE)\rangle=\langle \lambda, \cX^{\max}\rangle >0$ (Proposition \ref{proposition: Kempfs intersection number and blowups}). Therefore $q$ is unstable in $\cY$ and thus $c>0$. Since the $\cS_\alpha^{\cZ_c}$ stratify $\cZ_c$, the $\ind_{\cZ_c}^\cY(\cS_\alpha^{\cZ_c})$ stratify $\cS_c$, so $q$ is contained in $\ind_{\cZ_c}^\cY(\cS_\alpha^{\cZ_c})$ for a unique $\alpha\in \bGamma$ and thus $p$ is contained in a unique $p_*\left(\ind_{\cZ_c}^\cY(\cS_\alpha^{\cZ_c})\setminus \cE\right)$. It is left to check that $((d(\cX),c),\alpha)\in \bGamma$. This follows because $d_i^\alpha\geq i+z(\cZ_c)\geq i+1$ by \ref{item 5 stratification theorem} and Lemma \ref{lemma: centres of unstable strata have bigger central rank}.
\end{enumerate}

Now we check that $\abs{\cS^\cX_{> \alpha}}\coloneqq \bigcup_{\beta > \alpha} r_\beta\left(\abs{\cS_\beta^\cX}\right)$ is closed for all $\alpha\in \bGamma$. If $\alpha=(d(\cX),\infty)$, then $\abs{\cS^\cX_{> \alpha}}=\varnothing$, which is closed. If $\alpha=((d(\cX),c),\alpha')$ with $c\in \Q_{>0}$, then 
\[\abs{\cS^\cX_{> \alpha}}=p\left( \gr^{-1}\left(\abs{\cS_{> \alpha'}^{\cZ_c}}\right)\cup \bigcup_{c'>c}\abs{\cS_{c'}}\right)\cup \abs{\cX^{\max}},\]
which is closed by induction and because $p$ is proper. If $d_0^\alpha<d(\cX)$, then $\abs{\cS^\cX_{> \alpha}}=\abs{\cS^\cU_{> \alpha}}\cup \abs{\pi^{-1}\pi(\cX^{\max})}$, which is also closed by induction and because $\pi$ is universally closed.

It is left to show properties \ref{item 4 stratification theorem} and \ref{item 5 stratification theorem}. The former is true by construction. The latter is true for $i=0$ because $d(\cX)\geq z(\cX)$, and it is true if $d_0^\alpha<d(\cX)$ by induction, using the result for $\cU$. The remaining case is when $\alpha$ is of the form $\alpha=((d(\cX),c),\alpha')$. Then $\cS^{\cZ_c}_{\alpha'}\neq \varnothing$ and thus $d^\alpha_{i+1}=d^{\alpha'}_i\geq i+z(\cZ_c)\geq i+1+z(\cX)$ by induction and by Lemma \ref{lemma: centres of unstable strata have bigger central rank}.
\end{proof}

\begin{definition}[The balancing stratification]\label{definition: balancing stratification}
Let $\cX$ be a normed noetherian good moduli stack with affine diagonal. The \emph{balancing stratification} of $\cX$ is the sequential stratification $(\cS^\cX_\alpha)_{\alpha\in \bGamma}$ of $\cX$ from \Cref{theorem: sequential stratifications for good moduli stacks}.
\end{definition}

\begin{remark}\label{remark: no reduncancy in the definition of Gamma}
It is not hard to construct, using the convex-geometric picture introduced later (\Cref{corollary: iterated balanced filtration for states equals that for stacks}), that for every $\alpha\in \Gamma$ there is a normed good moduli stack $\cX$ with $\cS^\cX_\alpha\neq \varnothing$. We can in fact take $\cX$ to be of the form $\cX=\A^n_k/\G_{m,k}^l$, with $k$ any field.
\end{remark}

\begin{remark}[Refining {$\Theta$}-stratifications]\label{remark: refining theta-stratifications}
Suppose that $\cX$ is an algebraic stack over an algebraic space $B$, satisfying \Cref{assumption: basic assumptions}, and endowed with a linear form $\ell$ and a norm $q$ that define a $\Theta$-stratification $(\cS_c)_{c\in \Q_{\geq 0}}$ such that all the centres $\cZ_c$ are quasi-compact with affine diagonal and have good moduli spaces. Then each $\cZ_c$ has a well-defined balancing stratification, and thus $\cX$ is endowed with the sequential stratification $(\ind_{\cZ_c}^\cX(\cS_{\alpha}^{\cZ_c}))_{(c,\alpha)\in \Q_{\geq 0}\times \bGamma}$, indexed by the poset $\Q_{\geq 0}\times \bGamma$ with lexicographic order. This gives a natural sequential refinement of the $\Theta$-stratification $(\cS_c)_{c\in \Q_{\geq 0}}$. The sequential filtration that this stratification associates to every point can be interpreted as a refined Harder-Narasimhan filtration.
\end{remark}

The balancing stratification is well-behaved under certain pullbacks.

\begin{proposition}[Compatibility with pullback]\label{proposition: compatibility of balancing stratification with pullback}
Let $f\colon \cX'\to \cX$ be a norm-preserving morphism between normed noetherian good moduli stacks with affine diagonal. Let $\pi\colon \cX\to X$ and $\pi'\colon \cX'\to X'$ be the good moduli spaces. Assume further that $f$ is either
\begin{enumerate}
\item a closed immersion, or
\item it fits in a cartesian diagram
\[\begin{tikzcd}
    {\cX'} & \cX \\
    {X'} & {X.}
    \arrow["h", from=2-1, to=2-2]
    \arrow[from=1-2, to=2-2]
    \arrow[from=1-1, to=2-1]
    \arrow["f", from=1-1, to=1-2]
    \arrow["\ulcorner"{anchor=center, pos=0.125}, draw=none, from=1-1, to=2-2]
\end{tikzcd}\]
\end{enumerate}
Then the balancing stratification on $\cX'$ is the pullback along $f$ (\Cref{definition: pullback of sequential stratifications}) of the balancing stratification on $\cX$, that is, for all $\alpha\in \bGamma$ we have $\cS^{\cX'}_\alpha=f^*\cS^\cX_\alpha$.
\end{proposition}
\begin{proof}
We prove the statement by induction on $N(\cX)=d(\cX)-z(\cX)$. The base case is $N(\cX)=-\infty$, corresponding to the empty stack, for which the statement is obvious.

We have an upper semicontinuous function $r$ on $\abs{X}$ given by
\[\begin{tikzcd}[ampersand replacement=\&,column sep=scriptsize,row sep=tiny]
    {r\colon \abs{X}} \& \N \\
    p \& {d(\pi^{-1}(p)),}
    \arrow[from=1-1, to=1-2]
    \arrow[maps to, from=2-1, to=2-2]
\end{tikzcd}\]
where $d(-)$ denotes “maximal stabiliser dimension” (\Cref{definition: maximal stabiliser dimension}). For $d\in \N$, let $X_{\leq d}$ be the open subspace of $X$ with $\abs{X_{\leq d}}=\{p\in \abs{X}\st r(p)\leq d\}$ and let $\cX_{\leq d}=\pi^{-1}(X_{\leq d})$. In the case where $\cX'=\cX_{\leq d}$ and $f$ is the inclusion $\cX_{\leq d}\to \cX$, the result follows from \Cref{item 2 stratification theorem} in \Cref{theorem: sequential stratifications for good moduli stacks}.

Note that in both cases we have $\cX'_{\leq d}= f^{-1}(\cX_{\leq d})$ for all $d\in \N$. To prove the statement for given $\alpha\in \bGamma$, we may assume $d(\cX)=d_0^\alpha$ by replacing $\cX$ by $\cX_{\leq d_0^\alpha}$ and $\cX'$ by $\cX'_{\leq d_0^\alpha}=f^{-1}(\cX_{\leq d_0^\alpha})$. If $d(\cX')<d(\cX)$, then $\cS_\alpha^{\cX'}=\varnothing=f^*(\cS_\alpha^\cX)$. Thus we may assume $d\coloneqq d(\cX)=d(\cX')$. We prove that $\cS_\alpha^{\cX'}=f^*\cS_\alpha^\cX$ in this case.

First, we have $(\cX')^{\max}=f^{-1}(\cX^{\max})$ by \cite[Proposition C.5]{_Edidin_CanonicalreductionofstabilizersforArtinstackswithgoodmodulispaces}, since in both cases $f$ is stabiliser-preserving (that is, the map $\cI_{\cX'}\to f^*\cI_\cX$ between inertia stacks is an isomorphism). Therefore $\cS^{\cX'}_{(d,\infty)}=f^*\cS^\cX_{(d,\infty)}$. We may thus assume $c\coloneqq c^\alpha_0\in \Q_{>c}$ and write $\alpha=((d,c),\beta)$ with $\beta\in \Gamma$. In both cases we have a diagram
\[\begin{tikzcd}
    {\cY'} & {X'\times_X\cY} & \cY \\
    & {X'} & X
    \arrow[from=1-3, to=2-3]
    \arrow[from=2-2, to=2-3]
    \arrow[from=1-2, to=2-2]
    \arrow[from=1-2, to=1-3]
    \arrow["\iota", hook, from=1-1, to=1-2]
    \arrow[from=1-1, to=2-2]
\end{tikzcd}\]
with $\iota$ a closed immersion. The $\Theta$-stratification $(\cS'_c)_{c\in \Q_\geq 0}$ on $\cY'$ is the pullback of the $\Theta$-stratification $(\cS_c)_{c\in \Q_{\geq 0}}$ on $\cY$ by \Cref{proposition: compatibility of Theta-stratification with pullback}. Let $c>0$, let $\cZ_c$ be the centre of $\cS_c$, let $\cZ_c'$ be the centre of $\cS_c'$, and let $g\colon \cZ_c'\to \cZ_c$ be the natural map. By cartesianity of the square
\[\begin{tikzcd}[ampersand replacement=\&,cramped]
    {\cZ'_c} \& {\cZ_c} \\
    {\cY'} \& \cY
    \arrow["f", from=2-1, to=2-2]
    \arrow[from=1-2, to=2-2]
    \arrow["g", from=1-1, to=1-2]
    \arrow[from=1-1, to=2-1]
\end{tikzcd}\]
we see that $g$ is a composition of a closed immersion and a base change of a map between the good moduli spaces of $\cZ'_c$ and $\cZ_c$.
By induction, since $N(\cZ_c)<N(\cX)$, we have $\cS_\beta^{\cZ'_c}=g^*(\cS^{\cZ_c}_\beta)$. It follows from \Cref{item 3 stratification theorem} in \Cref{theorem: sequential stratifications for good moduli stacks} that $\cS^{\cX'}_\alpha=f^*(\cS^\cX_\alpha)$.
\end{proof}

\begin{remark}
A crucial ingredient that makes \Cref{proposition: compatibility of balancing stratification with pullback}  possible is that the closed substack structure of the maximal dimension stabiliser locus is compatible with pullbacks. It would not hold if we had used the reduced substack structure on $\abs{\cX^{\max}}$ instead.
\end{remark}

We will often use the following form of \Cref{proposition: compatibility of balancing stratification with pullback}.

\begin{corollary}[Compatibility with fibres]\label{corollary: compatibility of balancing stratification with fibres}
Let $\cX$ be a normed noetherian algebraic stack with affine diagonal and a good moduli space $\pi\colon \cX\to X$, let $k$ be a field and let $x\colon \Spec k\to \cX$ be a point. Let $\cF=\pi^{-1}\pi(x)$. Then the balancing stratification of $\cF$ is the pullback along $\cF\to \cX$ of the balancing stratification of $\cX$. \hfill \qedsymbol{}
\end{corollary}

The balancing stratification defines, for every point $x$ of the stack $\cX$, a canonical sequential filtration of $x$.
\begin{definition}[Iterated balanced filtration]\label{definition: iterated balanced filtration}
Let $\cX$ be a normed noetherian good moduli stack with affine diagonal, let $k$ be a field and let $x\in \cX(k)$ be a $k$-point. The \emph{iterated balanced filtration} $\lambda_{\text{ib}}(x)\in \qinfilt(\cX,x)$ of $x$ is defined as follows. There is a unique $\alpha\in \bGamma$ such that the point $x$ factors through $\cS_\alpha^\cX\to \cX$. There is then a unique (up to unique isomorphism) $k$-point $y$ of $\cS_\alpha^\cX$ mapping to $x$. The point $y$ naturally belongs to $\qinfilt(\cX,x)$, and we set $\lambda_{\text{ib}}(x)\coloneqq y$.
\end{definition}

\begin{remark}[Derived stacks]
The concept of good moduli space for derived Artin stacks has recently been introduced in \cite{derivedgoodmodulispaces}, and blow-ups along the substack of points with maximal stabiliser dimension in this context have been studied in \cite{_Hekking_Rydh_Savvas_Derivedstabilizerreduction}. On the other hand, $\Theta$-stratifications are also available in derived geometry \cite{_HalpernLeistner_DerivedThetastratificationsandtheDequivalenceconjecture}. Therefore we expect that the theory of balancing stratifications extends to normed derived Artin stacks $\cX$ with good moduli space. However, in view of \Cref{proposition: compatibility of balancing stratification with pullback} we anticipate that the iterated balanced filtration a point $x$ of $\cX$ that one gets from this theory only depends on the classical truncation of $\cX$.
\end{remark}

\subsection{Examples}\label{section: examples}

In moduli theory there are many natural instances of normed good moduli stacks, for which the balancing stratification is defined. We collect here a few examples.

\subsubsection{Stacks proper over a good moduli stack}
Assume the framework of \Cref{theorem: theta stratification proper over gms}, that is, $\cX$ is a normed noetherian good moduli stack with affine diagonal, $f\colon \cY\to \cX$ is representable and proper and $\ell$ is an $f$-positive linear form on $\cY$. Then we have a $\Theta$-stratification $(\cS_c)_{c\in \Q_{\geq 0}}$ of $\cY$ and every centre $\cZ_c$ is again a normed noetherian good moduli stack with affine diagonal. Therefore each $\cZ_c$ has canonically a balancing stratification $(\cS^{\cZ_c}_\alpha)_{\alpha\in \bGamma}$. By \Cref{remark: refining theta-stratifications}, the $\Theta$-stratification $(\cS_c)_{c\in \Q_{\geq 0}}$ is canonically refined by the sequential stratification 
\[\left(\ind_{\cZ_c}^\cY\left(\cS_\alpha^{\cZ_c}\right)\right)_{(c,\alpha)\in \Q_{\geq 0}\times \bGamma}.\]

\subsubsection{Geometric Invariant Theory}

Let $k$ be a field, let $G$ be a linearly reductive affine algebraic group over $k$ admitting a split maximal torus $T$, with Weyl group $W$, and endowed with a norm on cocharacters (\Cref{definition: norm on cocharacters of a group}). Let $A$ be a finite type $k$-algebra and consider an action of $G$ on $\Spec A$. The quotient stack $\cX=\Spec(A)/G$ has a good moduli space $\cX\to \Spec(A^G)$. Given any $G$-equivariant projective morphism $f\colon Y\to \Spec A$ and an ample linearisation on $Y$ (that is, a line bundle on $Y/G$ ample with respect to $h=f/G\colon Y/G\to \Spec(A)/G$), the previous example applied to $h$ gives a sequential stratification of $Y/G$ indexed by $\Q_{\geq 0}\times \bGamma$. For every $k$-point $y\in Y^\ss(k)$, the iterated balanced filtration $\lambda_{\text{ib}}(y)$ of $y$ can be seen as a sequence $\lambda_0,\ldots,\lambda_n\in \Gamma^\Q(G)$ of commuting rational cocharacters of $G$, considered up to certain equivalence relation, by \Cref{remark: description sequential filtrations quotient stacks}.

This stratification was first defined by Kirwan \cite{_Kirwan_RefinementsoftheMorsestratificationofthenormsquareofthemomentmap} in the case where $k=\C$, $A=\C$ and $Y\to \Spec \C$ is smooth, building on the ideas introduced in \cite{_Kirwan_PartialDesingularisationsofQuotientsofNonsingularVarietiesandtheirBettiNumbers}. The indexing set used in \cite{_Kirwan_RefinementsoftheMorsestratificationofthenormsquareofthemomentmap} is different from the one used here, and it depends on the quotient presentation of $Y/G$. The strata obtained in \cite{_Kirwan_RefinementsoftheMorsestratificationofthenormsquareofthemomentmap} are open and closed substacks of the strata defined here, which does not make a substantial difference. This further partition of the strata arises in two different ways. First, instead of the $\Theta$-stratification $(\cS_c)_{c\in \Q_{\geq 0}}$ of $Y/G$ considered here, indexed by $\Q_{\geq 0}$, the stratification considered by Kirwan is indexed by $\Gamma^\Q(T)/W$ (see also \cite{_Kirwan_Cohomologyofquotientsinsymplecticandalgebraicgeometry}). The set $\Gamma^\Q(T)/W$ can be seen as parametrising certain unions of connected components on $\Grad_\Q(Y/G)$ by \cite[Theorem 1.4.8]{_HalpernLeistner_Onthestructureofinstabilityinmodulitheory}, and the strata obtained in this way are closed and open substacks of the $\cS_c$. The same kind of difference on indexing sets appears each time a blow-up is performed in the construction. Also, each time a locus of maximal dimension stabilisers is considered, for example $\left(Y^\ss/G\right)^{\max}$, Kirwan writes it as a disjoint union of loci of the form $G\left((Y^\ss)^R\right)/G$, where $R$ is a reductive subgroup of $G$ of maximal dimension such that $(Y^\ss)^R$ is nonempty. This further refines the indexing set and the stratification, but again only breaking down the strata into pieces that are closed and open.

Even in the case considered in \cite{_Kirwan_RefinementsoftheMorsestratificationofthenormsquareofthemomentmap}, the fact that the balancing stratification has the structure of a sequential stratification is a novelty of our approach. The definition of the iterated balanced filtration (\Cref{definition: iterated balanced filtration}) is new also in this case.

\subsubsection{Quiver representations}
Let $Q$ be a quiver with set of vertices $Q_0$, set of arrows $Q_1$ and source and target maps $s,t\colon Q_1\to Q_0$. Let $d$ be a dimension vector for $Q$ and consider the moduli stack $\cRep(Q,d)$ of representations of $Q$ with dimension vector $d$ over an algebraically closed field $k$ of characteristic $0$. A \emph{central charge} $Z$ for $Q$ is defined by a family $(a_i)_{i\in Q_0}$ with $a_i\in \Q\oplus i\Q_{>0}\subset \C$. For a finite dimensional representation $E$ of $Q$ we set $Z(E)=\sum_{i\in \Q}a_i\dim E_i$. This defines a linear form $\ell$ and a norm $q$ on $\cRep(Q,d)$ as follows. A graded point $g\colon B\G_{m,k}\to \cRep(Q,d)$ corresponds to a representation $E$ with dimension vector $d$ and a direct sum decomposition $E=\bigoplus_{c\in \Z} E_c$. We then set
\begin{align}\label{align: ell and q from Z}
\ell(g)&=\sum_{c\in \Z} -c\Re Z(E_c) \\ 
\label{equation: norm on graded points from central charge}q(g)&=\sum_{c\in \Z} c^2\Im Z(E_c).
\end{align}
Therefore the central charge $Z$ determines a normed good moduli stack $\cRep(Q,d)^\ss$ of semistable quiver representations. In fact, semistability with respect to $\ell$ is a King stability \cite{_King_Moduliofrepresentationsoffinitedimensionalalgebras} condition. We can also see $\ell$-semistable representations as Bridgeland semistable representations of slope $\pi/2$. The linear form $\ell$ comes from the rational line bundle on $\cRep(Q,d)$ given by the rational character $\prod_{i\in Q_0} \det_{\GL_{d_i}}^{\Re a_i}$ of $\prod_{i \in Q_0}\GL_{d_i}$.

The iterated balanced filtration of a semistable representation $E$ coincides with the iterated weight filtration (or HKKP filtration) of $E$ defined by Haiden-Katzarkov-Kontsevich-Pandit in \cite{_Haiden_Semistabilitymodularlatticesanditeratedlogarithms}. This fact is proven in \cite{IbanezNunez_Blowupsandlattices}. Hence the balancing stratification of $\cRep(Q,d)^\ss$ can be seen as a stratification by type of HKKP filtration.

\subsubsection{Vector bundles on a curve}\label{example: vector bundles curve}
Let $C$ be a smooth projective curve over $\C$ and consider the stack $\Bun(C)_{r,d}$ of vector bundles on $C$ of rank $r$ and degree $d$. The open substack $\Bun(C)_{r,d}^\ss\subset \Bun(C)_{r,d}$ of semistable vector bundles admits a good moduli space. This can be proven either using GIT \cite{_Huybrechts_Thegeometryofmodulispacesofsheaves} or by intrinsic methods \cite{alper_tajakka_projectivityvectorbundles}.
The stack $\Bun(C)_{r,d}$ has a norm on graded points $q$ given by the rank as follows. A graded point $g\colon B\G_{m,\C}\to \Bun(C)_{r,d}$ corresponds to a vector bundle $E$ of degree $d$ and rank $r$ and a direct sum decomposition $E=\bigoplus_{c\in \Z}E_c$ as sum of subbundles. We set
\[q(g)=\sum_{c\in \Z}c^2\rk(E_c).\]
Hence $\Bun(C)_{r,d}^\ss$ is a normed good moduli stack (and it is also noetherian with affine diagonal), so the balancing stratification and the iterated balanced filtration of every semistable vector bundle are defined. In fact, the centres of the $\Theta$-stratification of $\Bun(C)_{r,d}$ by Harder-Narasimhan type have good moduli spaces and inherit a norm from $\Bun(C)_{r,d}$. Therefore, the balancing stratification of the centres gives a sequential stratification of $\Bun(C)_{r,d}$ that refines the stratification by Harder-Narasimhan type (\Cref{remark: refining theta-stratifications}). This produces, for every vector bundle $E$, a sequential filtration of $E$ that refines its Harder-Narasimhan stratification.

In this setting, a coarser version of the balancing stratification of $\Bun(C)_{r,d}^\ss$ was defined and studied by Kirwan \cite{_Kirwan_ModulispacesofbundlesoverRiemannsurfacesandtheYangMillsstratificationrevisited}. Each stratum in Kirwan's stratification is a connected component of a locally closed substack of the form $\bigcup_{\alpha\in \bGamma}\cS^\cX_{((n,c),\alpha)}$, for $n\in \N$ and $c\in \Q_{> 0}\cup\{\infty\}$, where $\cX=\Bun(C)_{r,d}^\ss$. Therefore Kirwan's stratification can be thought of as the stratification by type of balanced filtration. Kirwan calls the filtrations associated to this stratification \emph{balanced $\delta$-filtrations of maximal triviality}.

In \cite{IbanezNunez_Blowupsandlattices}, we show that the iterated balanced filtration for a semistable vector bundle $E$ coincides with the iterated HKKP filtration of the lattice of semistable subbundles of $E$ of the same slope. Hence the (first step of) the HKKP filtration of $E$ should correspond to Kirwan's balanced $\delta$-filtration of maximal triviality of $E$.

\subsubsection{Bridgeland semistable objects}\label{example: Bridgeland}
Let $X$ be a projective scheme over an algebraically closed field $k$ of characteristic $0$, and consider a Bridgeland stability condition \cite{_Bridgeland_Stabilityconditionsontriangulatedcategories} given by the heart $\cC\subset \Db(X)$ of a t-structure on the derived category of $X$ and a central charge $Z\colon K_0^{\text{num}}(\Db(X))\to \C$, where $K_0^{\text{num}}(\Db(X))$ is the numerical Grothendieck group. Let $\cA=\Ind(\cC)$ be the ind-completion of $\cC$. Under certain natural assumptions on $(\cC,Z)$, good moduli stacks of Bridgeland semistable objects can be constructed as open substacks of $\cM_\cA$, the moduli stack of objects in $\cA$ defined in \cite[Section 7]{_Alper_Existenceofmodulispacesforalgebraicstacks} following \cite{_Artin_Abstracthilbertschemes}. We assume that $Z$ is algebraic in the sense that $Z\left(\Db(X)\right)\subset \Q\oplus i\Q$, that $Z$ factors through a finite free quotient of $K_0^{\text{num}}(\Db(X))$, that $\cC$ satisfies the generic flatness condition and that certain boundedness conditions also hold (see \cite[Theorem 6.5.3]{_HalpernLeistner_Onthestructureofinstabilityinmodulitheory} and \cite[Example 7.29]{_Alper_Existenceofmodulispacesforalgebraicstacks} for details). Under these assumptions, $\cM_\cA$ is an algebraic stack with affine diagonal locally of finite type over $k$. 

For a numerical class $v\in K_0^{\text{num}}(\Db(X))$, there is an open and closed substack $\cM_v\subset \cM_\cA$ of objects in class $v$ and, by our boundedness hypothesis, there is a quasi-compact open substack $\cM_v^\ss$ of Bridgeland semistable objects. From the general results in \cite{_Alper_Existenceofmodulispacesforalgebraicstacks}, it follows that $\cM_v^\ss$ admits a (proper) good moduli space \cite[Example 7.29]{_Alper_Existenceofmodulispacesforalgebraicstacks}. The imaginary part of the central charge defines a norm $q$ on graded points of $\cM_\cA$ as in the previous examples. If $g\colon B\G_{m,k}\to \cM_\cA$ corresponds to $E=\bigoplus_{c\in \Z}E_c$ in $\cC$, then we define $q(g)$ by the formula
\[q(g)=\sum_{c\in \Z} c^2\Im Z(E_c)\]
as above (this is the norm used in \cite[Chapter 6]{_HalpernLeistner_Onthestructureofinstabilityinmodulitheory} to define the numerical invariant on $\cM_v$, by \cite[Lemma 6.4.8]{_HalpernLeistner_Onthestructureofinstabilityinmodulitheory}). Therefore $\cM_v^\ss$ is naturally a normed good moduli stack, noetherian and with affine diagonal, and thus the balancing stratification and the iterated balanced filtration of every point are defined. 

For $E\in \cM_v^\ss(k)$, the iterated balanced filtration of $E$ coincides with the HKKP filtration of the lattice of semistable subobjects of $E$ of the same slope. This follows from \Cref{theorem: correspondence HKKP filtration and iterated balanced filtration}.

\subsubsection{K-semistable Fano varieties}
Let $k$ be a field of characteristic $0$. It has recently been established that there is an algebraic stack $\cX^{\text{K}}_{n,V}$ of finite type over $k$ with affine diagonal parametrising families of $K$-semistable Fano varieties over $k$ of dimension $n$ and volume $V$, and that $\cX^{\text{K}}_{n,V}$
admits a proper good moduli space \cite{_Alper_ReductivityofautgroupofKpolystableFano,openness-K-semistability,_Blum_OnpropernessofKmodulispacesandoptimaldegenerationsofFanovarieties} (see \cite{Xu-book-K-stability} for a book account of these results). The $L^2$-norm of a test configuration \cite[Section 2.3]{_Blum_OnpropernessofKmodulispacesandoptimaldegenerationsofFanovarieties} defines a norm on graded points of $\cX^{\text{K}}_{n,V}$ \cite[Lemmas 2.36 and 8.45]{Xu-book-K-stability}, and thus the balancing stratification of $\cX^{\text{K}}_{n,V}$ is defined.

It is plausible that the iterated balanced filtration of a smooth $K$-semistable Fano variety $X$ over $\C$ is related to the asymptotics of the Calabi flow on $X$. This would provide a refinement of the results in \cite{ChenSun-Calabiflow,ChenSunWang-KahlerRicciflow}.

\subsubsection{$G$-bundles and gauged maps}
Let $C$ be a smooth projective over $\C$ and let $G$ be a connected reductive group. The notion of semistability for principal $G$-bundles over $C$ was defined in \cite{ramanathan-stable}. A moduli space of semistable $G$-bundles on $C$ was constructed in \cite{ramanathan-thesis-i,ramanathan-thesis-ii} using GIT. The stack of semistable $G$-bundles $\Bun_G(C)^\ss$ is thus an open substack of the stack $\Bun_G(C)$ of all $G$-bundles that admits a good moduli space. A choice of norm on cocharacters of $G$ induces a norm on $\Bun_G(C)$ as follows. Let $T$ be a maximal torus of $G$ with Weyl group $W$, then we have identifications
\begin{align*}\Grad(\Bun_G(C))&=\uHom(B\G_{m,\C},\uHom(C,BG))=\uHom(C,\uHom(B\G_{m,\C},BG))\\
&=\bigsqcup_{\lambda\in \Gamma^\Z(T)/W}\uHom(L(\lambda),C)=\bigsqcup_{\lambda\in \Gamma^\Z(T)/W}\Bun_{L(\lambda)}(C)\end{align*}
by \cite[Theorem 1.4.8]{_HalpernLeistner_Onthestructureofinstabilityinmodulitheory}. Therefore there is a natural map 
\[\pi_0(\Grad(\Bun_G(C)))\to \pi_0(\Grad(BG))=\Gamma^\Z(T)/W,\] 
along which we can pullback the norm on $BG$ to get a norm on cocharacters of $\Bun_G(C)$. This gives $\Bun_G(C)^\ss$ a natural structure of normed good moduli stack, and therefore defines the iterated balanced filtration for any semistable $G$-bundle. We expect the Yang-Mills flow for $G$-bundles \cite{_Atiyah_TheYangMillsEquationsoverRiemannSurfaces} to be related to the iterated balanced filtration, in analogy with \cite{_Haiden_Iteratedlogarithmsandgradientflows}, which deals with the case $G=\GL_{n,\C}$.

The moduli stack of $G$-bundles is a special case of the general framework of gauged maps studied in \cite{gauged-maps}. For a projective-over-affine scheme $X$ over $\C$ endowed with an action of $G$ and $n\geq 0$, Halpern-Leistner and Fernandez Herrero define a stack $\cM^G_n(X)$ parametrising families of Kontsevich stable maps from $C$ to the quotient stack $X/G$. Taking $X=\Spec \C$ and $n=0$ recovers the moduli stack of $G$-bundles on $C$. For suitable numerical invariants $\mu$, the semistable locus $\cM^G_n(X)^{\mu\dash\ss}$ admits a good moduli space. The numerical invariant $\mu$ depends on a choice of norm on cocharacters of $G$, that then gives a norm on graded points of $\cM^G_n(X)$ similarly to the case of $G$-bundles \cite[Definition 2.21]{gauged-maps}. Therefore $\cM^G_n(X)^{\mu\dash\ss}$ is a normed good moduli stack and the balancing stratification of $\cM^G_n(X)^{\mu\dash\ss}$ is defined. The construction in \cite{gauged-maps} works over a general noetherian base $S$ of characteristic $0$ with affine diagonal. In this setting 
one still gets a noetherian normed good moduli stack $\cM^G_n(X)^{\mu\dash\ss}$ with affine diagonal, and its balancing stratification is thus defined.

\section{Chains of stacks}\label{section: chains of stacks}
We introduce the formalism of \emph{chains of stacks} (\Cref{definition: chain of stacks}) as a tool to compute the iterated balanced filtration. For every chain of stacks there is an associated sequential filtration (\Cref{definition: Qinfinity filtration of a chain}). We give two different constructions of chains. The first, the \emph{balancing chain} (\Cref{construction: balancing chain}), is very close to the balancing stratification and it computes the iterated balanced filtration. The second, the \emph{torsor chain} (\Cref{construction: torsor chain}) is more closely related to combinatorial versions of the iterated balanced filtration. The main theorem of this section states that the torsor chain also computes the iterated balanced filtration (\Cref{theorem: torsor computes the iterated balanced filtration}). This fact will be used to relate the iterated balanced filtration to convex geometry (\Cref{theorem: balancing chain of states corresponds to torsor chain} and \Cref{corollary: iterated balanced filtration for states equals that for stacks}) and, in a sequel to this paper, to artinian lattices (\Cref{theorem: correspondence HKKP filtration and iterated balanced filtration}).

\subsection{Chains}
Let $k$ be a field. A \emph{$k$-pointed stack} is an algebraic stack $\cX$ together with a $k$-point $x\colon \Spec k\to \cX$. $k$-Pointed stacks form a 2-category as follows. A morphism $(\cX,x)\to (\cY,y)$ of $k$-pointed stacks is a morphism $f\colon \cX\to \cY$ of stacks and a 2-isomorphism $\alpha\colon f\circ x\to y$. The composition of $(f,\alpha)\colon (\cX,x)\to (\cY,y)$ and $(g,\beta)\colon (\cY,y)\to (\cZ, z)$ is $(g\circ f, \beta\circ(1_g \ast \alpha))$, where $\ast$ denotes horizontal composition. If $(f,\alpha),(f',\alpha')\colon (\cX,x)\to (\cY,y)$ are morphisms of $k$-pointed stacks, a 2-morphism $(f,\alpha)\to (f',\alpha')$ is a 2-morphism $\gamma\colon f\to f'$ such that
\[\begin{tikzcd}[column sep=small]
x\circ f \arrow[dr,"\alpha" swap] \arrow[rr,"1_x\ast \gamma"] & & x\circ f' \arrow[dl, "\alpha'"] \\
& y \end{tikzcd}\]
commutes.

\begin{definition}[Chain]\label{definition: chain of stacks}
A \emph{chain} $(\cX_n,x_n,\gamma_n,u_n)_{n\in\N}$ of $k$-pointed stacks is data:
\begin{enumerate}
\item For each $n\in \N$, a $k$-pointed normed stack $(\cX_n,x_n)$, where $\cX_n$ is of finite presentation over $\Spec k$ with affine diagonal and such that $\cX_n\to \Spec k$ is a good moduli space.
\item For each $n\in \N$, a $\Q$-filtration $\gamma_n\in \qfilt(\cX_n,x_n)$ of $x_n$.
\item For each $n\in \N$, a representable, separated and norm-preserving morphism 
\[u_n\colon (\cX_{n+1}, x_{n+1})\to (\Grad_\Q(\cX_n),\gr \gamma_n).\]
\end{enumerate}

We say that the chain $(\cX_n,x_n,\gamma_n,u_n)_{n\in\N}$ is \emph{bounded} if there is $N\in \N$ such that, for all $n\geq N$, we have that $\gamma_n=0$ and $u_n$ induces an isomorphism between $\cX_{n+1}$ and $\cX_n$, seen as the closed and open substack of $\Grad_\Q(\cX_n)$ of “trivial gradings”.

A \emph{morphism $f\colon (\cX_n',x_n',\gamma_n',u_n')\to (\cX_n,x_n,\gamma_n,u_n)$ of chains} consists of morphisms $f_n\colon (\cX_n', x_n')\to (\cX_n,x_n)$ of pointed stacks, together with isomorphisms $f(\gamma_n')\to \gamma_n$ of filtrations and 2-commutative squares

\[ \begin{tikzcd}
(\cX_{n+1}',x_{n+1}') \arrow[r,"u_n'"]\arrow[d,swap,"f_{n+1}"] & (\Grad_\Q(\cX_n'),\gr \gamma_n') \arrow[d,"\Grad_\Q(f_n)"] \\
(\cX_{n+1},x_{n+1}) \arrow[r,"u_n"]& (\Grad_\Q(\cX_n), \gr \gamma_n)
\end{tikzcd}
\]
of pointed stacks.
\end{definition}

Suppose that $(\cX_n,x_n,\gamma_n,u_n)_{n\in\N}$ is a bounded chain. For $n\in \N$, we define a map $c_n\colon \cX_{n+1}\to \Grad_\Q^{n+1}(\cX_0)$ by 
\[c_n=\Grad_\Q^n(u_0)\circ \Grad_\Q^{n-1}(u_1)\circ \cdots \circ u_n.\]
By Proposition \ref{proposition: representable and separated induced injection on set of filtrations}, $c_n$ induces an injection 
\[\qfilt(\cX_{n+1},x_{n+1})\to \qfilt(\Grad_\Q^{n+1}(\cX_0),c_n(x_{n+1})).\]
Define $\lambda_{n+1}\in \qfilt(\Grad_\Q^{n+1}(\cX_0),c_n(x_{n+1}))$ to be the image of $\gamma_{n+1}$ under this injection. Define also $\lambda_0\coloneqq \gamma_0 \in \qfilt(\cX_0,x_0)$.

\begin{lemma}
There is a canonical isomorphism $c_n(x_{n+1})\simeq \gr \lambda_n$, for all $n\in \N$.
\end{lemma}
\begin{proof}
For $n=0$, $c_0(x_1)=u_0(x_1)\simeq \gr \gamma_0=\gr \lambda_0$ is given by $u_0$ as a pointed map. For $n>0$, we have
\begin{align*}
\gr \lambda_n=\gr\left(\Filt_\Q(c_{n-1})(\gamma_n)\right)=\Grad_\Q(c_{n-1})(\gr \gamma_n)\simeq \\
\Grad_\Q(c_{n-1})\left(u_n(x_{n+1})\right)=c_n(x_{n+1}),
\end{align*}
as desired.
\end{proof}

Note that there is $N\in \N$ such that $\lambda_n=0$ for $n\geq N$, because the chain is bounded. The lemma gives canonical isomorphisms
\[\qfilt(\Grad_\Q^{n+1}(\cX_0),c_n(x_{n+1}))\cong\qfilt(\Grad_\Q^{n+1}(\cX_0),\gr\lambda_n)\]
for all $n$.
Therefore, by \Cref{remark: description set of infty filtrations}, $(\lambda_n)_{n\in \N}$ defines an element of $\qinfilt(\cX_0,x_0)$.

\begin{definition}[Sequential filtration of a chain]\label{definition: Qinfinity filtration of a chain}
The \emph{$\Q^\infty$-filtration associated to} the bounded chain $(\cX_n,x_n,\gamma_n,u_n)_{n\in\N}$ is $(\lambda_n)_{n\in \N}\in \qinfilt(\cX_0,x_0)$.
\end{definition}

\subsection{The balancing chain}
We now construct a chain closely related to the balancing stratification.

\begin{construction}\label{construction: balancing chain}
Let $\cX$ be a normed noetherian algebraic stack with affine diagonal and a good moduli space $\pi\colon \cX\to X$. Let $x\colon \Spec k\to \cX$ be a $k$-point, with $k$ a field. We define, a chain $(\cX_n,x_n,\gamma_n,u_n)_{n\in\N}$ over $k$ as follows. 

We set $(\cX_0,x_0)=(\pi^{-1}(\pi(x)),x)$, which has good moduli space $\Spec k$. For $n\in \N$, assume $\cX_n$ and $x_n$ are defined. We now define $\cX_{n+1}$, $x_{n+1}$, $\gamma_n$ and $u_n$ in terms of $\cX_n$ and $x_n$. We consider two cases:

\emph{Case 1.} The point $x_n$ is closed in $\cX_n$. We then define $\cX_{n+1}=\cX_n$, $x_{n+1}=x_n$, $\gamma_{n}=0$ the trivial filtration in $\qfilt (\cX_n,x_n)$, and $u_n\colon (\cX_{n+1},x_{n+1})\to (\Grad_\Q(\cX_n),\gr \gamma_{n})$ given by the “trivial grading” map.

\emph{Case 2.} The point $x_n$ is not closed in $\cX_n$. Then we consider the blow-up $p\colon \cB=\Bl_{\cX_n^{\max}}\cX_n\to \cX_n$, where $\cX_n^{\max}$ is the closed substack of points with maximal dimension stabiliser \cite[Appendix C]{_Edidin_CanonicalreductionofstabilizersforArtinstackswithgoodmodulispaces}. In fact, $\abs{\cX_n^{\max}}$ is a singleton consisting of the unique closed point of $\abs{\cX_n}$ \cite[Proposition 9.1]{_Alper_GoodmodulispacesforArtinstacks}, which is different from $x_n$ by assumption. Thus the point $x_n$ lifts uniquely to a point of $\cB$ that we still denote $x_n$.

By \Cref{theorem: theta stratification proper over gms} and \Cref{example: stratification on a blow-up}, $\cB$ has a $\Theta$-stratification induced by the natural $p$-ample line bundle on $\cB$ and the norm. Let $\cS$ be the locally closed $\Theta$-stratum containing $x_n$ and let $\cZ$ be its centre (\Cref{definition: Theta-stratification my version}). Let $\gamma_n\in\qfilt(\cB,x_n)$ be the Harder-Narasimhan filtration of $x_n$ in $\cB$ (\Cref{definition: HN filtration})
and let $x_{n+1}=\gr \gamma_n\in \cZ(k)$. We identify $\qfilt(\cB,x_n)=\qfilt(\cX_n,x_n)$ by \Cref{proposition: proper map induces bijection of set of filtrations}, and under this identification $\gamma_n$ is the balanced filtration of $x_n$ in $\cX_n$ (\Cref{definition: balanced filtration} and \Cref{proposition: Kempfs intersection number and blowups}). By \Cref{theorem: theta stratification proper over gms}, $\cZ$ has a good moduli space
$\pi_\cZ\colon \cZ\to Z$. We set $\cX_{n+1}=\pi_\cZ^{-1}(\pi_\cZ(x_{n+1}))$. We define $u_n$ to be the composition
\[
u_n\colon\cX_{n+1}\rightarrow \cZ \hookrightarrow \Grad_\Q(\cB) \xrightarrow{\Grad_\Q(p)}\Grad_\Q(\cX_n),
\]
which is representable and separated, since applying $\Grad_\Q$ preserves representability and separatedness, and the first two maps are immersions. 

The stack $\cX_{n+1}$ inherits a norm from $\cX_n$ along the composition $\cX_{n+1}\to \Grad_\Q(\cX_n)\to \cX_n$. By commutativity of 
\[ \begin{tikzcd}
\Filt_\Q(\cB) \arrow[r,""]\arrow[d,swap,"\gr"] & \Filt_\Q(\cX_n) \arrow[d,"\gr"] \\
\Grad_\Q(\cB) \arrow[r,""]& \Grad_\Q(\cX_n)
\end{tikzcd}
\]
we get a pointed morphism $u_n\colon (\cX_{n+1},x_{n+1})\to (\Grad_\Q(\cX_n),\gr\gamma_n)$.
\end{construction}

\begin{definition}[The balancing chain]\label{definition: balancing chain stack}
Let $\cX$ be a normed noetherian good moduli stack with affine diagonal, let $k$ be a field and let $x\in \cX(k)$ be a $k$-point. The \emph{balancing chain} of $(\cX,x)$ is the chain $(\cX_n,x_n,\gamma_n,u_n)_{n\in \N}$ of $k$-stacks constructed in \Cref{construction: balancing chain}.
\end{definition}

\begin{remark}\label{remark: basic properties balancing chain}
The following properties of the balancing chain are clear from \Cref{theorem: Kempf}:
\begin{enumerate}
\item For every $n\in \N$, we have $\gamma_n=0$ if and only if $x_n$ is closed in $\cX_n$.
\item For every $n\in \N$, we have that $\ev_0(\gamma_n)$ is the unique closed point of $\cX_n$.
\end{enumerate}
\end{remark}

\begin{lemma}\label{lemma: about central ranks and balancing chain}
Assume the setup of \Cref{definition: balancing chain stack}. For every $n\in \N$ such that $x_n$ is not closed in $\cX_n$, we have that $z(\cX_n)\geq z(\cX_0)+n$, where $z(-)$ denotes central rank (\Cref{definition: central rank}).
\end{lemma}
Note as well that $z(\cX_0)\geq z(\cX)$.
\begin{proof}
Let $n\in \N$ and suppose that $x_n$ is not closed in $\cX_n$. Then $z(\cX_{n+1})>z(\cX_{n})$ by \Cref{lemma: centres of unstable strata have bigger central rank}. The results follows by induction.
\end{proof}

\begin{corollary}\label{corollary: balancing chain is bounded}
Assuming the setup of \Cref{definition: balancing chain stack}, the balancing chain of $(\cX,x)$ is bounded.
\end{corollary}
\begin{proof}
If it was not bounded, then for every $n\in \N$ we would have that $x_n$ is not closed in $\cX_n$ and thus that $z(\cX_n)\geq n$ by \Cref{lemma: about central ranks and balancing chain}. This would contradict the bound $z(\cX_n)\leq d(\cX_0)$, where $d(-)$ denotes maximal stabiliser dimension (\Cref{definition: maximal stabiliser dimension}).
\end{proof}

\begin{proposition}\label{proposition: balancing chain computes iterated balanced filtration}
Let $\cX$ be a normed noetherian good moduli stack with affine diagonal, and let $x\in \cX(k)$ be a field-valued point. Then the sequential filtration of $x$ associated to the balancing chain of $(\cX,x)$ (\Cref{definition: Qinfinity filtration of a chain,definition: balancing chain stack}) equals the iterated balanced filtration of $x$ (\Cref{definition: iterated balanced filtration}).
\end{proposition}
\begin{proof}
Let $(\cX_n,x_n,\gamma_n,u_n)_{n\in \N}$ be the balancing chain of $(\cX,x)$, and let \[\lambda_{\text{bc}}(x)\in \qfilt(\cX_0,x_0)=\qfilt(\cX,x)\] be its associated sequential filtration. Note that, by \Cref{corollary: compatibility of balancing stratification with fibres}, the iterated balanced filtration $\lambda_{\text{ib}}(x)$ of $(\cX,x)$ equals that of $(\cX_0,x)$. 

We prove the statement by induction on $N(\cX)=d(\cX)-z(\cX)$. If $N(\cX)=0$, then $x$ is closed in $\cX_0$ and therefore $\lambda_{\text{ib}}(x)=\lambda_{\text{bc}}(x)=0$. If $x$ is not closed in $\cX_0$ then, using the notation of \Cref{construction: balancing chain}, $\gamma_0\in \qfilt(\cX_0,x)$ is the balanced filtration of $x$ (\Cref{definition: balanced filtration}). From \Cref{proposition: iterative cartesian diagram sequential filtrations}, we get a map $\phi\colon \qinfilt(\cX_1,x_1)\hookrightarrow \qinfilt(\Grad_\Q(\cX_0),\gr \gamma_0)\to \qinfilt(\cX_0,x_0)$ that, in the notation of \Cref{remark: description set of infty filtrations}, can be written as
\[(\lambda_0,\lambda_1,\ldots)\mapsto (\gamma_0,\lambda_0,\lambda_1,\ldots).\]
By construction of the balancing chain, we have $\lambda_{\text{bc}}(x)=\phi(\lambda_{\text{bc}}(x_1))$, and by construction of the balancing stratification, we have $\lambda_{\text{ib}}(x)=\phi(\lambda_{\text{ib}})(x_1)$. By \Cref{lemma: centres of unstable strata have bigger central rank}, $N(\cX_1)<N(\cX_0)$. Therefore $\lambda_{\text{ib}}(x_1)=\lambda_{\text{bc}}(x_1)$ by induction, and hence $\lambda_{\text{ib}}(x)=\lambda_{\text{bc}}(x)$.
\end{proof}

\begin{remark}
One could define chain similar to the balancing chain but where one replaces $\cX_n$ with $\overline{\{x_n\}}$ (with reduced structure) at each step, and \Cref{proposition: balancing chain computes iterated balanced filtration} would still be true by the same reasons.
\end{remark}

\subsection{The torsor chain}
We introduce a second chain whose construction is similar to that of the balancing chain, but where at each step the exceptional divisor is replaced by the natural $\G_m$-torsor over it. Then we prove (\Cref{theorem: torsor computes the iterated balanced filtration}) that this new chain also computes the iterated balanced filtration.

\begin{construction}\label{construction: torsor chain}
Let $\cX$ be a normed noetherian algebraic stack with affine diagonal and with a good moduli space $\pi\colon \cX\to X$. Let $k$ be a field and let $x\in \cX(k)$ be a $k$-point.
We construct a chain $(\cY_n,y_n,\eta_n, v_n)_{n\in \N}$ inductively as follows. 

We set $(\cY_0,y_0)=(\pi(\pi^{-1}(x)),x)$. Suppose that $(\cY_n,y_n)$ is defined. We define $\eta_n$, $v_n$ and $(\cY_{n+1},y_{n+1})$ in terms of $(\cY_n,y_n)$.

\textit{Case 1.} The point $y_n$ is in $\cY_n^{\max}$. In that case, we set $\eta_n=0$, $(\cY_{n+1},y_{n+1})=(\cY_n,y_n)$ and $v_n\colon (\cY_{n+1},y_{n+1})\to (\Grad_\Q(\cY_n),\gr\eta_n)$ the “trivial grading” map.

\textit{Case 2.} The point $y_n$ is not in $\cY_n^{\max}$. Then $y_n$ lifts uniquely to $\cB\coloneqq \Bl_{\cY_n^{\max}}\cY_n$, which is canonically endowed with a linear form on graded points, coming from the exceptional divisor, and with the induced norm on graded points (\Cref{theorem: theta stratification proper over gms} and \Cref{example: stratification on a blow-up}). We let $\eta_n\in \qfilt(\cY_n,y_n)$ be the HN filtration of $y_n$ in $\cB$. Let $\cS$ be the locally closed $\Theta$-stratum of $\cB$ containing $y_n$, and let $\cZ$ be its centre. Let $\cE$ be the exceptional divisor and $\cN\to \cE$ the natural $\G_{m,k}$-torsor, that is, $\cN$ is the complement of the zero section inside the total space of the normal cone to $\cY_n^{\max}\to \cY_n$. If $\cL$ is the ideal sheaf of the exceptional divisor $\cE$, then $\cN=\A((\cL\vert_{\cE})^\vee)^*$, where $\A(-)^*$ denotes total space minus zero section.

\begin{lemma}\label{lemma: the centre lives over the exceptional divisor}
The “forget the grading” map $h\colon \cZ\to \cB$ factors through $\cE\to \cB$. As a consequence, the open immersion $\cZ\to \Grad_\Q(\cB)$ factors through the closed immersion $\Grad_\Q(\cE)\to \Grad_\Q(\cB)$. In particular, the induced map $\cZ\to \Grad_\Q(\cE)$ is an open immersion.
\end{lemma}
\begin{proof}
The centre $\cZ$ carries a canonical rational $B\G_m$-action, which endows every quasi-coherent sheaf $\cM$ on $\cZ$ with a $\Q$-grading $\cM=\bigoplus_{c\in \Q}\cM_c$. To see how this grading originates, let us assume for simplicity that the $B\G_m$-action is integral. Then for any morphism $f\colon T\to \cZ$ there is an associated map $g\colon T\times B\G_m\to \cZ$. The pullback $g^*\cM$ is a $\G_m$-equivariant sheaf on $T$, that is, a $\Z$-graded sheaf on $T$, and the underlying sheaf is $(T\to T\times B\G_m)^*g^*\cM=f^*\cM$. Thus $f^*\cM$ has a canonical $\Z$-grading for every $f$, and this gives the $\Z$-grading for $\cM$ itself.

Let $\cL$ be the ideal sheaf of $\cE$. Since $\cZ$ is the centre of a stratum, unstable with respect to the linear form $\langle -,\cL\rangle$, we must have that the $\Q$-grading on $h^*\cL$ is concentrated in degree $-1$. Indeed, we know that for every map $\Spec l\to \cZ$, with $l$ a field, the pullback $(\Spec l\times B\G_m \to \cZ)^*h^*\cL$ is concentrated in degree $-1$. Thus by Nakayama's lemma, $(h^*\cL)_c=0$ for all $c\in \Q\setminus \{-1\}$.

On the other hand, the structure sheaf $\cO_\cZ$ is concentrated in degree $0$. There is a map $\cL\to \cO_\cB$ because $\cL$ is an ideal sheaf. Pulling this map back along $h$, we get a homomorphism $h^*\cL\to \cO_\cZ$ that must be $0$ for degree reasons. Therefore $h$ factors through $\cE\to \cB$. Since $\Grad_\Q(\cE)=\cE\times_\cB \Grad_\Q(\cB)$, the lemma follows.
\end{proof}
\end{construction}

Let $\cM$ be the fibre product
\[\begin{tikzcd}[ampersand replacement=\&]
    {\cM} \& \cN \\
    \cZ \& \cE,
    \arrow[from=2-1, to=2-2]
    \arrow[from=1-2, to=2-2]
    \arrow[from=1-1, to=1-2]
    \arrow[from=1-1, to=2-1]
    \arrow["\ulcorner"{anchor=center, pos=0.125}, draw=none, from=1-1, to=2-2]
\end{tikzcd}\]
and let $y_{n+1}$ be any $k$-point of $\cM$ lying above $\gr \eta_n$. Any two choices of $y_{n+1}$ are related by a unique $\G_{m,k}$-torsor automorphism of $\cM\to \cZ$. Now $\cM$ has a good moduli space $\pi_{\cM}\colon \cM\to M$ because $\cM\to \cZ$ is affine and $\cZ$ has a good moduli space. Finally, we let $\cY_{n+1}\coloneqq\pi_{\cM}^{-1}(\pi_\cM(y_{n+1}))$ and $v_n\colon (\cY_{n+1},y_{n+1})\to (\cM,y_{n+1}) \to(\cZ,\gr\eta_n)\to (\Grad_\Q(\cX_n),\gr\eta_n)$.  We are abusively denoting $\eta_n$ both the filtration of $x_n$ in $\cX_n$ and in $\cB$, and by $\gr \eta_n$ the associated graded point in the two cases.

\begin{definition}[Torsor chain]\label{definition: torsor chain stack}
Let $\cX$ be a normed noetherian good moduli stack with affine diagonal, let $k$ be a field and let $x\in \cX(k)$ be a $k$-point. The \emph{torsor chain} of $(\cX,x)$ is the chain $(\cY_n,y_n,\eta_n,v_n)_{n\in \N}$ of $k$-stacks from \Cref{construction: torsor chain}.
\end{definition}

\begin{theorem}\label{theorem: torsor computes the iterated balanced filtration}
Let $k$ be a field and let $(\cX,x)$ be a $k$-pointed normed noetherian good moduli stack with affine diagonal. Then the torsor chain of $(\cX,x)$ is bounded and its associated $\Q^\infty$-filtration (\Cref{definition: torsor chain stack,definition: Qinfinity filtration of a chain}) equals the iterated balanced filtration of $(\cX,x)$ (\Cref{definition: iterated balanced filtration}).
\end{theorem}
\begin{proof}
Let $(\cX_n,x_n,\gamma_n,u_n)_{n\in\N}$ be the balancing chain of $(\cX,x)$ and let $(\cY_n,y_n,\eta_n,v_n)_{n\in\N}$ be the torsor chain. The iterated balanced filtration is well-behaved under field extension by \Cref{proposition: compatibility of balancing stratification with pullback}, and the torsor chain also commutes with base field extension by a similar argument.
Therefore, we may assume that $k$ is algebraically closed. 

 By \Cref{corollary: stack with good moduli space a point is a quotient algebraically closed case}, $\cX_0=R_0/G_0$, where $R_0$ is an affine scheme and $G_0$ is the stabiliser of a $k$-point $p_0$ of $R_0$ which is also the unique closed $G_0$-orbit of $R_0$. Let $q$ be the norm on graded points of $\cX_0$. There is a closed immersion $BG_0\to \cX_0$ given by the point $p_0$, so $BG_0$ inherits a norm from $\cX_0$ by pullback. On the other hand, the quotient presentation gives a representable map $\cX_0\to BG_0$, which then induces a new norm $q'$ on $\cX_0$ from the norm on $BG_0$. We claim that the balanced and the torsor chains do not depend on whether we endow $\cX_0$ with $q$ or with $q'$. Indeed, for any $\lambda\colon \Theta_k\to \cX_0$ such that $\langle\lambda,\cX_0^{\max}\rangle>0$, we have that $\ev_0(\lambda)$ lies in $\cX_0^{\max}$, and thus $\gr \lambda$ factors through $BG_0\to \cX_0$. Therefore $q(\lambda)=q'(\lambda)$.  Since these are the only norms that show up in the definition of the $\Theta$-stratification of $\cB_0=\Bl_{\cX_0^{\max}}\cX_0$, we have that $q$ and $q'$ define the same $\Theta$-stratification on $\cB_0$. If $\cZ_0$ is the centre of an unstable stratum of $\cB_0$, then for any filtration $\Theta_k\to \cZ_0$, the composition $\Theta_k\to\cZ_0\to \cX_0$ factors through $BG_0\to \cX_0$. Therefore $q$ and $q'$ induce the same norm on $\cZ_0$ and, consequently, also on the $\G_{m}$-torsor $\cN_0\to \cZ_0$ considered in the torsor chain. We may thus assume that $q=q'$ without loss of generality. The stack $\cX_n$ can be written as $\cX_n=R_n/G_n$, with $R_n$ affine and $G_n$ the stabiliser of a $k$-point $p_n$ of $R_n$ which is also the unique closed $G_n$-orbit; and the norm on $\cX_n$ is induced by a norm on cocharacters of $G_n$. This last statement follows by considering the representable map $\cX_n\to \cX_0$ and the fact for $n=0$.

We start by introducing some natural extra structure on the balancing chain. Let $N$ the smallest natural number such that $\gamma_N=0$. For every $n\leq N$, there are line bundles $\cL_{n,1},\cdots,\cL_{n,n}$ on $\cX_n$ constructed inductively as follows. If $n<N$, then $\cL_{n+1,i}=(\cX_{n+1}\to \cX_n)^*\cL_{n,i}$ for all $i=1,\ldots,n$. The map $\cX_{n+1}\to \cX_n$ we are considering is the composition of $u_n\colon \cX_{n+1}\to \Grad_\Q(\cX_n)$ and the “forget the grading” map $\Grad_\Q(\cX_n)\to \cX_n$. Using the notation of \Cref{construction: balancing chain}, we have a relatively ample line bundle $\cL_{\cB}=\cO_\cB(-\cE)$ on the blow-up $\cB=\Bl_{\cX_n^{\max}}\cX_n$, and we let $\cL_{n+1,n+1}=(\cX_{n+1}\to \cB)^*\cL_\cB$, where the morphism considered is the composition $\cX_{n+1}\to \cZ\to \Grad_\Q(\cB)\to \cB$.

There is also a natural non-degenerate rational $B\G_m^n$-action on $\cX_n$ for each $n\leq N$. Indeed, if this action has been constructed for $\cX_n$, with $n<N$, we construct it for $\cX_{n+1}$ as follows. Again, we use the notations of Case 2 in \Cref{construction: balancing chain}. By \Cref{lemma: lifting BG_m^n actions along blow-ups}, the rational $B\G_m^n$-action on $\cX_n$ lifts canonically to a rational $B\G_m^n$-action on $\cB$. Since $\cZ$ is an open substack of $\Grad_\Q(\cB)$, it has a natural $B\G_m^n\times B\G_m$-action, and it is nondegenerate by the argument in the proof of \Cref{lemma: centres of unstable strata have bigger central rank}. This action now restricts to the closed substack $\cX_{n+1}$ of $\cZ$.

Let $T_n$ be a maximal torus of $G_n$. It follows from \cite[Theorem 1.4.8]{_HalpernLeistner_Onthestructureofinstabilityinmodulitheory} that the rational $B\G_m^n$-action on $\cX_n$ is given by rational one-parameter subgroups $\beta_{1},\ldots,\beta_{n}$ of the centre $Z(G_n)$ that act trivially on $R_n$. Indeed, the $B\G_m^n$-action is given by an element $\beta\in \Hom(\G_{m,k}^n,T_n)\otimes_\Z \Q$ and a connected component $C$ of the fixed point locus $R_n^{\beta,0}$ such that the induced “forget the grading” map $C/L_{G_n}(\beta)\subset \Grad_\Q(\cX_n)\to \cX_n$ is an isomorphism. Since $R_n$ has a point fixed by $G_n$, it must be $L_{G_n}(\beta)=G_n$, and consequently $C=R_n$, from what we see that $(\beta_{1},\ldots,\beta_{n})=\beta$ has the desired property.

The blow-up $\cB_n=\Bl_{\cX_n^{\max}}\cX_n$can be written as $\cB_n=B_n/G_n$, with $B_n$ the blow-up of $R_n$ along a closed subscheme $R_n^{\max}$ given by $\cX_n^{\max}$. The centre of the locally closed $\Theta$-stratum of $\cB_n$ containing the lift of $x_n$ to $\cB_n$ is $\cZ_n=Z_n/L_{G_n}(\beta_{n+1})$, where $\beta_{n+1}\in \Gamma^\Q(T_n)$ corresponds to $\gr \gamma_n$ and $Z_n$ is an open subscheme of the fixed point locus $B_n^{\beta_{n+1},0}$. The group $G_{n+1}$ is the stabiliser of a point of $Z_n$, so it is identified with a subgroup of $L_{G_n}(\beta_{n+1})$ containing $\beta_1,\ldots,\beta_{n+1}$ in its centre.  With this identifications we regard the $\beta_i$ as independent of $n$.

\begin{claim}
For all $n< N$, for all $1\leq j\leq n$, we have the equalities $\langle \beta_j, \cL_{n+1,n+1}\rangle=0$ and $\langle \beta_{n+1},\cL_{n+1,n+1}\rangle =1$. Therefore $\langle \beta_i,\cL_{n,j}\rangle=\delta_{ij}$ for $n\leq N$ and $i\leq j\leq n$.
\end{claim}

The first equality follows from the fact that the $\beta_1,\ldots,\beta_n$ act trivially on $R_n$ and thus also on the normal cone to $R_n^{\max}$. The second follows from the definition of $\gamma_n$ as the minimiser of $\norm{-}$ subject to $\langle \gamma_n,\cL_{n+1,n+1}\rangle=1$, and the fact that $\gr \gamma_n$ is given by $\beta_{n+1}$.

\begin{claim}\label{claim: orthogonality of the betas}
For all $n\leq N$ and all $1\leq i\neq j\leq n$, we have $(\beta_i,\beta_j)=0$, where $(-,-)$ denotes the inner product in $\Gamma^\Q(T_n)$ defined by the norm on cocharacters of $G_n$.
\end{claim}
We prove the claim by induction. Let $n<N$. By the Linear Recognition Theorem \ref{theorem: linear recognition}, a point $z\in (B_n)^{\beta_{n+1},0}$ belongs to $Z_n$ if and only if it is semistable with respect to the shifted linear form, which in this case is
\[\ell_n=\langle -, \cL_{\cB_n}\rangle-\dfrac{1}{\norm{\beta_{n+1}}^2}(\beta_{n+1},-),\]
where $\cL_{\cB_n}$ is the relatively ample line bundle on $\cB_n$.
Since $\beta_i$ for $i\leq n$ all fix $z$, we must have $\ell_n(\beta_i)\leq 0$ and $\ell_n(-\beta_i)\leq 0$, but $\ell_n(\beta_i)=\dfrac{1}{\norm{\beta_{n+1}}^2}(\beta_{n+1},\beta_i)$ and $\ell_n(-\beta_i)=-\dfrac{1}{\norm{\beta_{n+1}}^2}(\beta_{n+1},\beta_i)$, so $(\beta_{n+1},\beta_i)=0$. This proves the claim.

Let us denote $\cV_n=\A(\cL^\vee_{1,n})^*\times_{\cX_n} \cdots \times_{\cX_n} \A(\cL^\vee_{n,n})^*$, where $\A(-)^*$ denotes total space minus zero section. There is a natural quotient presentation $\cV_n=V_n/G_n$, where $V_n=\A(\cL^\vee_{n,1}\vert_{R_n})^*\times_{R_n}\cdots \times_{R_n} \A(\cL^\vee_{n,n}\vert_{R_n})^*$, which is a $\G_m^n$-torsor over $R_n$, and in particular carries a $\G_{m,k}^n$-action. Let $T'_n=\im (\beta_1,\ldots,\beta_n)\subset Z(G_n)$. Note that this is well-defined even if the $\beta_i$ are rational cocharacters and not necessarily integral.

\begin{claim}
The torus $T_n'$ acts on $V_n$ via a homomorphism $\delta_n\colon T_n'\to \G^n_{m,k}$. Moreover, $\delta_n$ is an isogeny.
\end{claim}

Recall that the torus $T_n'$ acts trivially on $R_n$. The $\G_m^n$-torsor $V_n/T_n'\to R_n/T_n'=R_n\times BT_n'$ is given by a map $r\colon R_n\times BT'_n\to B\G_{m,k}^n$, which is in turn given by a map $R_0\to \uHom(BT'_n,B\G_{m,k}^n)$. We have the equality $\uHom(BT'_n,B\G_{m,k}^n)=\bigsqcup_{\alpha\in \Hom(T'_n,\G_{m,k}^n)} B\G_m^n$ because $T'_n$ is a split torus and by \cite[Theorem 1.4.8]{_HalpernLeistner_Onthestructureofinstabilityinmodulitheory}. Therefore, since $R_n$ is connected, $r$ corresponds to a pair $(R_0\xrightarrow{o} B\G_{m,k}^n,\delta_n)$, where $o$ corresponds to the $\G_m^n$-torsor $V_n\to R_0$ and $\delta_n\in \Hom(T',\G_{m,k}^n)$. We recover $r$ as the composition $R_0\times BT'\xrightarrow{(o,B\delta_n)} B\G_{m,k}^n\times B\G_{m,k}^n\to B\G_{m,k}^n$, the last map being multiplication. The homomorphism $\delta_n$ induces a map $D(\delta_n)_\Q\colon \Gamma_\Q(\G_{m,k}^n)\to \Gamma_{\Q}(T'_n)$. If we take $\beta_1,\ldots,\beta_n$ as a basis of $\Gamma_{\Q}(T'_n)$ and the standard basis of $\Gamma_\Q(\G_{m,k}^n)$, then $D(\delta_n)_\Q$ is given by the matrix $\langle \beta_i,\cL_{n,j}\rangle$, which we have shown is upper triangular with $1$'s in the diagonal. Therefore $D(\delta_n)_\Q$ is an isomorphism and $\delta_n$ is an isogeny.

\begin{claim}
Let $\overline\cX_n=R_n/(G_n/T_n')$ and let $D_n=\ker \delta_n$. Then $V_n/(G_n/D_n)\cong \overline\cX_n$.
\end{claim}

The group $(G_n/D_n)/(T'_n/D_n)$ acts on the stack $V_n/(T'_n/D_n)$ and there is a natural isomorphism
\[V_n/(G_n/D_n)\cong \left(V_n/(T'_n/D_n)\right)/ \left((G_n/D_n)/(T'_n/D_n)\right),\]
by \cite[Remark 2.4]{Romagny_groupactionsonstacks}. Since $T_n'/D_n=\G_{m,k}^n$ and $Y_n\to R_n$ is a $\G_m^n$-torsor, we have $V_n/(T'_n/D_n)=R_n$. Noting that $(G_n/D_n)/(T'_n/D_n)=G_n/T'_n$, we get the desired isomorphism.

As a consequence of the claim, we see that the good moduli space of $\cV_n$ is $\Spec k$.

\begin{claim}
Let $f_n\colon \cV_n\to \cX_n$ be the $\G_m^n$-torsor map. Then $f_n^{-1}(\cX_n^{\max})=\cV_n^{\max}$ and $d(\cX_n)=d(\cV_n)+n$.
\end{claim}

Let $p'_n$ be a $k$-point of $V_n$ mapping to $p_n$ along $\cV_n\to \cX_n$. Necessarily $G_np'_n$ is the unique closed orbit inside $V_n$. Let $H_n$ be the stabiliser of $p_n'$. By \cite[Theorem 10.4, (5) and (6)]{_Alper_Theetalelocalstructureofalgebraicstacks} together with the fact that the good moduli space of $\cV_n$ is $\Spec k$, there is a locally closed $H_n$-equivariant subscheme $S_n$ of $V_n$ containing $x_n'$ such that $S_n/H_n\to \cV_n$ is an isomorphism. The map $S_n/(H_n/D_n)\to V_n/(G_n/D_n)$ must also be an isomorphism. Since $S_n$ has a point fixed by $H_n$, the maximal dimension stabiliser locus is $\cV_n^{\max}=S_n^{(H_n)_\circ}/H_n$, where $(H_n)_\circ$ is the reduced identity component of $H_n$ (see \cite[Appendix C]{_Edidin_CanonicalreductionofstabilizersforArtinstackswithgoodmodulispaces}). Similarly, $(V_n/(G_n/D_n))^{\max}=(S_n)^{(H_n/D_n)_\circ}/(H_n/D_n)$. Since $D_n$ acts trivially on $S_n$, we have $(S_n)^{(H_n/D_n)_\circ}=S_n^{(H_n)_\circ}$. This proves that, denoting $\rho_n\colon \cV_n\to \overline \cX_n$ the obvious map, we have $\cV_n^{\max}=\rho_n^{-1}(\cX_n^{\max})$. On the other hand, since $R_n$ has a point fixed by $G_n$, we have $\cX_n^{\max}=R_n^{(G_n)_\circ}/G_n$. Since $T'_n$ acts trivially on $R_n$, $\overline\cX_n^{\max}=R_n^{(G_n/T'_n)_\circ}/(G_n/T'_n)=R_n^{(G_n)_\circ}/(G_n/T'_n)$ so, denoting $q_n\colon \cX_n\to \overline \cX_n$, we have $q_n^{-1}(\overline\cX_n^{\max})=\cX_n^{\max}$. Hence $f_n^{-1}(\cX_n^{\max})=q_n^{-1}(\overline\cX_n^{\max})=\rho_n^{-1}(\overline \cX_n^{\max})=\cV_n^{\max}$. For the last statement, just note that $d(\cV_n)=d(\overline \cX_n)$ and that $d(\cX_n)=d(\overline \cX_n)+n$. This proves the claim.

Suppose that $n<N$. Let $x_n'$ be a $k$-point of $V_n$ mapping to $x_n$ along $\cV_n\to \cX_n$.

\begin{claim}
The balanced filtration of $(\cV_n,x'_n)$ equals the balanced filtration of $(\cX_n,x_n)$ under the injection $\qfilt(\cV_n,x_n')\to \qfilt(\cX_n,x_n)$.
\end{claim}

We identify the balanced filtration of $(\cX_n, x_n)$ with $\beta_{n+1}\in \Gamma^\Q(G_n)$. Since $\beta_{n+1}$ is orthogonal to the cocharacters in $T_n'$, by \Cref{claim: orthogonality of the betas}, and since $\cX_n\to \overline\cX_n$ preserves the substacks of maximal stabiliser dimensions, the image $\overline\beta_{n+1}\in \Gamma^\Q(G_n/T_n')$ of $\beta_{n+1}$ inside $G_n/T_n'$ is the balanced filtration of $(\overline \cX_n,x_n)$. Since $\rho_n\colon \cV_n\to \overline \cX_n$ is a gerbe banded by the finite group $D_n$, $\qfilt(\cV_n,x_n')=\qfilt(\overline \cX_n,x_n)$, and this equality identifies the balanced filtrations of $\cV_n$ and $\overline \cX_n$ because $\rho_n^{-1}(\overline \cX_n^{\max})=\cV_n^{\max}$. Therefore $\beta_{n+1}$ is also the balanced filtration of $(\cV_n,x'_n)$.

Let $\cC_n=\Bl_{\cV_n^{\max}}\cV_n=C_n/G_n$. Since $\cV_n\to \cX_n$ is flat, being a $\G_m^n$-torsor, and since $f_n^{-1}(\cX_n^{\max})=\cV_n^{\max}$, we have that the blow-ups form a cartesian diagram
\[\begin{tikzcd}
    {\cC_n} & {\cV_n} \\
    {\cB_n} & {\cX_n.}
    \arrow["{f_n}", from=1-2, to=2-2]
    \arrow[from=2-1, to=2-2]
    \arrow[from=1-1, to=1-2]
    \arrow[from=1-1, to=2-1]
    \arrow["\ulcorner"{anchor=center, pos=0.125}, draw=none, from=1-1, to=2-2]
\end{tikzcd}\]
Let $\cZ_n'$ be the centre of the locally closed $\Theta$-stratum of $\cC_n$ containing the lift of $x_n'$ to $\cC_n$.
\begin{claim}
There is a natural cartesian square
\[\begin{tikzcd}
    {\cZ_n'} & {\cV_n} \\
    {\cZ_n} & {\cX_n.}
    \arrow[from=1-1, to=1-2]
    \arrow["{f_n}", from=1-2, to=2-2]
    \arrow[from=2-1, to=2-2]
    \arrow[from=1-1, to=2-1]
    \arrow["\ulcorner"{anchor=center, pos=0.125}, draw=none, from=1-1, to=2-2]
\end{tikzcd}\]
\end{claim}
We denote $\beta_{n+1}=\lambda$ for simplicity. As mentioned above, the centre $\cZ_n$ is the semistable locus inside $\cB_n^\lambda\coloneqq B_n^{\lambda,0}/L_{G_n}(\lambda)$ for the shifted linear form
\[\langle -, \cL\rangle - \dfrac{1}{\norm{\lambda}^2}(\lambda,-),\]
where here $\cL\coloneqq \cL_{\cB_n}\vert_{\cB_n^\lambda}$ is the standard relatively ample line bundle on $\cB_n$, pulled back to $\cB_n^\lambda$. Similarly, $\cZ_n'$ is the semistable locus inside $\cC_n^\lambda\coloneqq C_n^{\lambda,0}/L_{G_n}(\lambda)$ for the form
\[\langle -, \cL\vert_{\cC_n^\lambda}\rangle - \dfrac{1}{\norm{\lambda}^2}(\lambda,-),\]
since $\cL\vert_{\cC_n^\lambda}$ is the standard relatively ample line bundle on $\cC_n$. First, note that the natural square
\[\begin{tikzcd}
    {\cC_n^\lambda} & {\cC_n} \\
    {\cB_n^\lambda} & {\cB_n}
    \arrow[from=2-1, to=2-2]
    \arrow[from=1-2, to=2-2]
    \arrow[from=1-1, to=1-2]
    \arrow[from=1-1, to=2-1]
    \arrow["\ulcorner"{anchor=center, pos=0.125}, draw=none, from=1-1, to=2-2]
\end{tikzcd}\]
is cartesian. This follows from cartesianity of the two squares
\begin{equation}\label{diagram: two squares}\begin{tikzcd}
    {\cC_n^\lambda} & {\cC_n} && {\cB_n^\lambda} & {\cB_n} \\
    {\overline\cC_n^\lambda} & {\overline\cC_n} && {\overline\cB_n^\lambda} & {\overline\cB_n,}
    \arrow[""{name=0, anchor=center, inner sep=0}, from=1-2, to=2-2]
    \arrow[from=2-1, to=2-2]
    \arrow[from=1-1, to=1-2]
    \arrow[from=1-1, to=2-1]
    \arrow["\ulcorner"{anchor=center, pos=0.125}, draw=none, from=1-1, to=2-2]
    \arrow[from=1-4, to=1-5]
    \arrow[from=1-5, to=2-5]
    \arrow[from=2-4, to=2-5]
    \arrow[""{name=1, anchor=center, inner sep=0}, from=1-4, to=2-4]
    \arrow["\ulcorner"{anchor=center, pos=0.125}, draw=none, from=1-4, to=2-5]
    \arrow["{\text{\normalsize and}}"{description}, draw=none, from=0, to=1]
\end{tikzcd}\end{equation}
where $\overline\cC_n=C_n/(G_n/D_n)$, $\overline \cC_n^\lambda=C_n^{\lambda,0}/L_{G_n/D_n}(\lambda)$, $\overline\cB_n=B_n/(G_n/T_n')$ and $\overline \cB_n^\lambda=B^{\lambda,0}_n/L_{G_n/T_n'}(\lambda)$, together with the isomorphisms $\overline \cC_n\cong \overline \cB_n$ and $\overline \cC_n^\lambda\cong \overline \cB_n^\lambda$. From the compatibility of the standard relatively ample line bundles on $\cC_n$, $\cB_n$ and $\overline \cC_n=\overline \cB_n$ under pullbacks and from the shape of the shifted linear forms, it follows that $(\overline \cC^\lambda_n)^\ss\times_{\overline \cC_n}\cC_n=(\cC_n^\lambda)^\ss$ and $\overline\cB_n^\lambda\times_{\overline \cB_n}\cB_n=(\cB_n^\lambda)^\ss$, and therefore $\cZ_n'=(\cC_n^\lambda)^\ss=(\cB_n^\lambda)^\ss\times_{\cB_n}\cC_n=\cZ_n\times_{\cX_n}\cV_n$, as desired. This proves the claim.

We now show by induction that $(\cV_n,x_n')\cong (\cY_n,y_n)$. For $n=0$ there is nothing to prove. If $n<N$ and $(\cV_n,x_n')\cong (\cY_n,y_n)$, then $\cY_{n+1}$ is constructed as follows. We take the standard relatively ample line bundle $\cL_{\cC_n}$ on $\cC_n$ and let $\cM=\A(\cL^\vee_{\cC_n}\vert_{\cZ_n'})^*$. We choose a point $y_{n+1}$ in $\cM$ mapping to $z\coloneqq \lim_{t\to 0} \lambda(t)x_n'\in \cZ_n'$, and we let $(\cY_{n+1},y_{n+1})$ be the fibre of the good moduli space of $\cN$ containing $y_{n+1}$. Note that by the previous claim, 
\[\cZ'_n=\A(\cL^\vee_{n,1}\vert_{\cZ_n})^*\times_{\cZ_n}\cdots \times_{\cZ_n}\A(\cL^\vee_{n,n}\vert_{\cZ_n})^*,\]
so actually $\cM=\A(\cL^\vee_{n,1}\vert_{\cZ_n})^*\times_{\cZ_n}\cdots \times_{\cZ_n}\A(\cL^\vee_{n,n}\vert_{\cZ_n})^*\times_{\cZ_n}\A(\cL^\vee_{\cB_n}\vert_{\cZ_n})^*$. The stack $\cX_{n+1}$ is the fibre of the good moduli space of $\cZ_n$ containing $\lim_{t\to 0}\lambda(t)x_n$, and from the definitions we have $\cV_{n+1}=\cX_{n+1}\times_{\cZ_n}\cM$. Since we have seen that $\cV_{n+1}$ has for moduli space $\Spec k$, it must be the fibre of the good moduli space of $\cM$. Therefore $\cV_{n+1}\cong \cY_{n+1}$. By either choosing the lifts $x_{n+1}'$ appropriately or applying a torsor automorphism of $\cY_n$, we can arrange that $y_{n+1}=x_{n+1}'$. 

In particular, from these isomorphisms it follows that $y_N$ is closed in $\cY_N$, so the torsor chain is bounded. We have seen that the balanced filtration of $\cV_n$ maps to the balanced filtration of $\cX_n$. Therefore the maps $\cY_n\to \cX_n$ give a morphism of chains (the compatibility of link morphisms $\cY_{n+1}\to \Grad_\Q(\cY_n)$ and $\cX_{n+1}\to \Grad_\Q(\cX_n)$ naturally follows from the construction). Since $\cY_0=\cX_0$, the sequential filtration of the torsor chain equals the sequential filtration of the balancing chain, as we wanted to show.
\end{proof}

\section{The iterated balanced filtration and convex geometry}
Despite its seemingly convoluted definition in terms of repeated blow-ups and {$\Theta$-}stratifications, we illustrate in this section how the iterated balanced filtration can be explicitly computed in terms of convex geometry and convex optimisation when the groups involved are diagonalisable. For an account on diagonalisable algebraic groups, we refer the reader to \cite{_Milne_Algebraicgroupsthetheoryofgroupschemesoffinitetypeoverafield}.

Suppose one is interested in the iterated balanced filtration of a geometric point $x\colon \Spec k \to \cX$ in a normed good moduli stack $\cX$. By taking the fibre of the good moduli space at $x$, we may assume that $\cX\to \Spec k$ is the good moduli space. If $y$ is the unique closed $k$-point of $\cX$, and $G$ is the stabiliser of $y$, then $\cX$ is of the form $\cX=\Spec A/G$ by \Cref{corollary: fibres of good moduli space are quotient stacks}. Since closed immersions have no  effect in the iterated balanced filtration, we may assume after embedding $\Spec A$ in a representation of $G$ that $\Spec A=V$ is a finite dimensional vector space on which $G$ acts linearly, and the point $x$ is given by some vector $x\in V$. From now we assume that $G$ is diagonalisable (for example, $G=\G_{m,k}^n$ is a split torus). Then the structure of $G$-representation $V$ is determined by a direct sum decomposition $V=\bigoplus_{\chi\in \Gamma_\Z(G)}V_\chi$, where $G$ acts on $V_\chi$ by the character $\chi$. The \emph{state} of $x$ (named after \cite{_Kempf_InstabilityinInvariantTheory}) is the finite set $\Xi=\{\chi\in \Gamma_\Z(G)\st p_\chi(x)\neq 0\}$, where $p_\chi\colon V\to V_\chi$ denotes the projection. Writing $V=V'\oplus V_0$ and considering the closed immersion $V'/G\cong \left(V'\times \{p_0(x)\}\right)/G\to V/G$, we may assume that $0\notin \Xi$. We may simplify the situation further by considering $\A^\Xi_k$ to be the product of $\#(\Xi)$ many copies of $\A^1_k$, endowed with the action of $G$ via the characters in $\Xi$. There is a $G$-equivariant linear closed immersion $\A^\Xi\to V$ sending the point $(1,\ldots,1)$ to $x$, so we may replace $V/G$ and $x$ with $\A^\Xi/G$ and $(1,\ldots,1)$. This is the pointed stack \emph{associated} to the state $\Xi$ (\Cref{definition: stack associated to a state}).

We will be able to determine the balanced filtration of $x$ in terms of its state $\Xi$. To this aim, we develop a theory of polarised states (where, in addition to the set $\Xi$, we have the data of a character $\alpha\colon G\to \G_{m,k}$) purely in combinatorial terms. We introduce a notion of filtration (\Cref{definition: filtrations of a polarised state}) and $\Q^\infty$-filtration (\Cref{definition: sequential filtration of a state}) for states, and an analogue of the balancing chain and of the iterated balanced filtration (\Cref{definition: balancing chain of a state and the iterated balanced filtration}). We then define a functor (\Cref{definition: stack associated to a state}) from a category of normed semistable polarised states to the category of pointed normed stacks over a field $k$ and show that the functor maps the balancing chain of states to the torsor chain of stacks (\Cref{theorem: balancing chain of states corresponds to torsor chain}), concluding that the iterated balanced filtration of the pointed stack coincides with that of the state (\Cref{corollary: iterated balanced filtration for states equals that for stacks}). Computing the iterated balanced filtration of a state boils down to a simple convex optimisation problem.

\subsection{Polarised states}
We develop here a combinatorial analogue of the theory of iterated balanced filtrations.

\begin{definition}[Polarised state]\label{definition: polarised state}
A \emph{polarised state} $\bXi$ is a triple $\bXi=(M,\Xi,\alpha)$ where
\begin{enumerate}
\item $M$ is a finite type abelian group,
\item $\Xi\subset M\setminus\{0\}$ is a finite subset (the \emph{state}), and
\item $\alpha\in M_\Q$ (the \emph{polarisation}).
\end{enumerate}
We denote $M_\Q=M\otimes_\Z \Q$, $M^\vee=\Hom(M,\Z)$ and $M^\vee_\Q=M^\vee\otimes_\Z \Q$. A \emph{normed polarised state} is a polarised state $\bXi=(M,\Xi,\alpha)$ together with a rational inner product on $M^\vee_\Q$.
\end{definition}

\begin{remark}
For the theory of polarised states developed here, the rational numbers can be replaced by any subfield of $\R$, but this generality is not necessary for our purposes.
\end{remark}

Let $\bXi=(M,\Xi,\alpha)$ be a polarised state. We denote $\langle -,-\rangle\colon M^\vee\times M\to \Z$ the duality pairing. For $\lambda\in M^\vee_\Q$, we let $H_\lambda=\{\beta \in M_\Q\st \langle \lambda,\beta\rangle\geq 0\}$ be the half-space defined by $\lambda$ and $\partial H_\lambda=\{\beta\in M_\Q\st \langle \lambda,\beta\rangle =0\}$ the hyperplane defined by $\lambda$. We will need a few basic notions about convex geometry for which a sufficient source is \cite{_Fulton_IntroductiontoToricVarieties}. A \emph{cone} in $M_\Q$ is a subset of the form $H_{\lambda_1}\cap \cdots \cap H_{\lambda_n}$ for $\lambda_1,\ldots,\lambda_n\in M_\Q^\vee$ or, equivalently, of the form $(\Q_{\geq 0})\chi_1+\cdots+(\Q_{\geq 0})\chi_k$ for $\chi_1,\ldots,\chi_k\in M_\Q$. We include the degenerate cases $\{0\}$ and $M_\Q^\vee$. If $C$ is a cone, a \emph{face} of $C$ is a subset of the form $C\cap \partial H_\lambda$, where $\lambda\in M^\vee_\Q$ is such that $C\subset H_\lambda$.

\begin{definition}[Semistable and polystable polarised states]\label{definition: semistable and polystable polarised state}
We say that the polarised state $\bXi$ is \emph{semistable} if $\alpha$ is in the cone $\cone\left(\Xi\right)$ generated by $\Xi$ inside $M_\Q$. We say that $\bXi$ is \emph{polystable} if it is semistable and the smallest face of $\cone\left(\Xi\right)$ containing $\alpha$ is $\cone\left(\Xi\right)$ itself (that is, $\alpha$ is in the relative interior of $\cone\left(\Xi\right)$). 
\end{definition}

We are abusing notation by denoting the image of $\Xi$ inside $M_\Q$ also by $\Xi$. By $\cone\left(\Xi\right)$ we will always mean a subset of $M_\Q$.

\begin{definition}[Filtrations of a state]\label{definition: filtrations of a polarised state}
Suppose that $\bXi$ is semistable. The set of \emph{rational filtrations} (or \emph{$\Q$-filtrations}) of $\bXi$ is
\[\qfilt(\bXi)\coloneqq \{\lambda\in N_\Q\st \langle \lambda,\Xi\rangle\geq 0 \text{ and } \langle \lambda,\alpha\rangle=0\}.\]
\end{definition}

We are using the notation $\langle \lambda,\Xi\rangle\geq 0$ to mean that $\langle \lambda,\chi\rangle \geq 0$ holds for all $\chi \in \Xi$.

\begin{definition}[Associated graded state]\label{definition: associated graded state}
Suppose that $\bXi$ is semistable and let $\lambda\in \qfilt(\bXi)$. The \emph{associated graded state} $\Grad_\lambda(\bXi)$ is the semistable polarised state $\Grad_\lambda(\bXi)=(M,\Xi_{\lambda,0},\alpha)$, where $\Xi_{\lambda,0}\coloneqq \Xi\cap \partial H_\lambda$.
\end{definition}

\begin{proof}[Proof that {$\Grad_\lambda(\bXi)$} is semistable]
Since $\Xi\subset H_\lambda$, we have the equality $\cone(\Xi)\cap \partial H_\lambda=\cone(\Xi\cap \partial H_\lambda)$. Therefore $\alpha\in \cone(\Xi_{\lambda,0})$.
\end{proof}

\begin{definition}[Sequential filtrations of a state]\label{definition: sequential filtration of a state}
Let $\bXi=(M,\Xi,\alpha)$ be a semistable polarised state. The set $\qinfilt(\bXi)$ of \emph{sequential filtrations} (or \emph{$\Q^\infty$-filtrations}) of $\bXi$ is the set of those sequences $(\lambda_n)_{n\in\N}$ in $M_\Q^\vee$ satisfying
\begin{enumerate}
\item $\lambda_n=0$ for $n>>0$, and
\item $\lambda_0\in \qfilt(\bXi)$ and, for all $n\in \N_{>0}$, we have 
\[\lambda_n\in \qfilt\left(\Grad_{\lambda_{n-1}}\left(\cdots\Grad_{\lambda_1}\left(\Grad_{\lambda_0}(\bXi)\right)\cdots\right)\right).\]
\end{enumerate}
\end{definition}

\begin{definition}[Morphism between states]\label{definition: morphism between polarised states}
A \emph{morphism} $\phi\colon \bXi_1=(M_1,\Xi_1,\alpha_1)\to \bXi_2=(M_2,\Xi_2,\alpha_2)$ between semistable polarised states is a surjective homomorphism $\phi\colon M_2\to M_1$ such that
\begin{enumerate}
\item for all $\chi\in \Xi_2$, either $\phi(\chi)\in \Xi_1$ or $\phi(\chi)=0$; and
\item $\alpha_2\in \cone\left(\Xi_2\cap \ker \phi\right)$.
\end{enumerate}
If $\bXi_1$ and $\bXi_2$ are normed, we say that $\phi$ is a morphism between normed semistable polarised states if the inner product on $(M_1)^\vee_\Q$ is the restriction of that on $(M_2)^\vee_\Q$ along the inclusion $(M_1)^\vee_\Q\hookrightarrow (M_2)^\vee_\Q$ defined by $\phi$.

With the obvious composition and identity, semistable (normed) polarised states form a category.
\end{definition}

\begin{lemma}\label{lemma: morphism between states induces map in set of filtrations}
Let $\phi\colon \bXi_1=(M_1,\Xi_1,\alpha_1)\to \bXi_2=(M_2,\Xi_2,\alpha_2)$ be a morphism between semistable polarised states. Then the injection $\phi^\vee_\Q\colon (M_1)_\Q^\vee\to (M_2)_\Q^\vee$ induces a map $\qfilt(\bXi_1)\to\qfilt(\bXi_2)$ between sets of rational filtrations.
\end{lemma}
\begin{proof}
If $\lambda \in \qfilt(\bXi_1)$ and $\chi\in \Xi_2$, then $\langle \phi^\vee_\Q(\lambda),\chi\rangle=\langle \lambda,\phi(\chi)\rangle\geq 0$, since either $\phi(\chi)=0$ or $\phi(\chi)\in \Xi_1$. Likewise, $\langle \phi^\vee_\Q(\lambda),\alpha_2\rangle=\langle \lambda,\phi(\alpha_2)\rangle= 0$, since $\phi(\alpha_2)=0$ in $(M_1)_\Q$. Therefore $\phi^\vee_\Q(\lambda)\in \qfilt(\bXi_2)$.
\end{proof}

In the situation of \Cref{lemma: morphism between states induces map in set of filtrations}, we will use the simpler notation $\phi(\lambda)\coloneqq\phi^\vee_\Q(\lambda)$.

\begin{proposition}\label{proposition: map induced between associated graded states}
Let $\phi\colon \bXi_1=(M_1,\Xi_1,\alpha_1)\to \bXi_2=(M_2,\Xi_2,\alpha_2)$ be a morphism between semistable polarised states and let $\lambda\in \qfilt(\bXi_1)$. The homomorphism $\phi\colon M_2\to M_1$ induces a map $\Grad_\lambda(\phi)\colon \Grad_\lambda(\bXi_1)\to \Grad_{\phi(\lambda)}(\bXi_2)$ between associated graded states.
\end{proposition}
\begin{proof}
The first condition to be checked is that for all $\chi\in (\Xi_2)_{\lambda,0}$, we have $\phi(\chi)=0$ or $\phi(\chi)\in (\Xi_1)_{\lambda,0}$. If $\phi(\chi)$ is not $0$, then $\phi(\chi)\in \Xi_1$, and since $\langle \lambda,\phi(\chi)\rangle=\langle \phi(\lambda),\chi\rangle=0$, we have indeed $\phi(\chi)\in (\Xi_1)_{\lambda,0}$. For the second condition, we have 
\begin{align*}\alpha_2\in \cone\left(\Xi_2\cap \ker \phi\right)\cap \partial H_{\phi(\lambda)} &=\cone\left(\Xi_2\cap \ker \phi\cap \partial H_{\phi(\lambda)}\right)\\
&= \cone\left((\Xi_2)_{\phi(\lambda),0}\cap \ker \phi\right),\end{align*}
as desired.
\end{proof}

\begin{definition}[Chain of states]\label{definition: chain of states}
A \emph{chain of normed semistable polarised states} is a sequence $(\bXi_n,\lambda_n,u_n)_{n\in\N}$ where
\begin{enumerate}
\item for each $n\in \N$, $\bXi_n$ is a normed semistable polarised state; 
\item $\lambda_n\in \qfilt(\bXi_n)$ for each $n\in \N$;
\item $u_n\colon \bXi_{n+1}\to \Grad_{\lambda_n}\left(\bXi_n\right)$ is a morphism of normed semistable polarised states; and
\item for all $n>>0$, $\lambda_n=0$ and $u_n$ is an isomorphism (boundedness).
\end{enumerate}
\end{definition}

Suppose that $(\bXi_n,\lambda_n,u_n)_{n\in\N}$ is a chain of normed semistable polarised states. We denote $c_n\colon \bXi_{n+1}\to \Grad_{\lambda_n}\left(\cdots \Grad_{\lambda_1}\left(\Grad_{\lambda_0}(\bXi_0)\right)\cdots \right)$ the map defined by
\[c_n=u_n\circ \Grad_{\lambda_n}(u_{n-1})\circ \cdots \circ \Grad_{\lambda_n}\left(\cdots \Grad_{\lambda_1}\left(\Grad_{\lambda_0}(u_0)\right)\cdots \right).\]
We are abusively denoting by $\lambda_n$ the image of $\lambda_n$ along the relevant maps between sets of filtrations, but note as well that all these maps are injective. The sequence $(c_{n-1}(\lambda_{n}))_{n\in \N}$ (where $c_{-1}(\lambda_0)\coloneqq \lambda_0$) is a $\Q^\infty$-filtration of $\bXi_0$.

\begin{definition}[Associated sequential filtration]\label{definition: sequential filtration associated to a chain of states}
The \emph{sequential filtration associated to the chain} $(\bXi_n,\lambda_n,u_n)_{n\in\N}$ is $(c_{n-1}(\lambda_{n}))_{n\in \N}\in \qinfilt\left(\bXi_0\right)$. We will simply denote it $(\lambda_n)_{n\in \N}$.
\end{definition}

\begin{definition}[Slice of a state]\label{definition: slice of a state}
Let $\bXi=(M,\Xi,\alpha)$ be a semistable polarised state. Let $F$ be the smallest face of $\cone(\Xi)$ containing $\alpha$, let $K$ be the subgroup of $M$ generated by $F\cap \Xi$, let $M'=M/K$, let $q\colon M\to M'$ be the quotient map, and set $\Xi'=q(\Xi\setminus F)$. We define the \emph{slice} of $\bXi$ to be the state $\bXi'=(M',\Xi',0)$ together with the morphism $q\colon \bXi'\to \bXi$ given by $q\colon M\to M'$. If $\bXi$ is normed, we regard $\bXi'$ as a normed polarised state, where $(M')_\Q^\vee$ inherits an inner product along the inclusion $q^\vee_\Q\colon(M')_\Q^\vee\to M_\Q^\vee$.
\end{definition}

Note that $\bXi'$ is semistable, since $0$ belongs to any cone. We are again abusing notation: the expression $F\cap \Xi$ denotes the subset of $\Xi$ consisting of the elements $\chi$ whose image in $M_\Q$ is contained in $F$.

\begin{proposition}\label{proposition: slicing does not change set of filtrations of state}
Let $\bXi=(M,\Xi,\alpha)$ be a semistable polarised state and let $q\colon \bXi'\to \bXi$ be its slice. Then the map $\qfilt(\bXi')\to \qfilt(\bXi)$ induced by $q$ is a bijection.
\end{proposition}
\begin{proof}
We use notation from \Cref{definition: slice of a state}. Let $\lambda\in \qfilt(\bXi)$. Then $\cone(\Xi)\cap \partial H_\lambda$ is a face of $\cone(\Xi)$ containing $\alpha$, so $F\subset \cone(\Xi)\cap \partial H_\lambda$. In particular, $\langle \lambda,\chi\rangle=0$ for all $\chi\in F\cap \Xi$. Therefore $\lambda$ has a preimage $\lambda'$ along $q^\vee_\Q\colon (M')_\Q^\vee\hookrightarrow M^\vee_\Q$, and $\langle \lambda',q(\chi)\rangle=\langle \lambda,\chi\rangle\geq 0$ for all $\chi\in \Xi\setminus F$, so $\lambda'\in \qfilt(\bXi')$. Hence the map $\qfilt(\bXi')\to \qfilt(\bXi)$ is surjective. Since it is also injective, it is a bijection.
\end{proof}

\begin{definition}[Complementedness of a filtration]\label{definition: complementedness of a filtration state}
Let $\bXi=(M,\Xi,\alpha)$ be a semistable polarised state and let $\bXi'=(M',\Xi',0)$ be its slice. Let $\lambda\in \qfilt(\Xi)=\qfilt(\Xi')$. We define the \emph{complementedness} $\langle \lambda,\fl\rangle\in \Q_{\geq 0}\cup\{\infty\}$ of $\lambda$ to be
\[\langle \lambda,\fl\rangle=\inf_{\chi\in \Xi'} \langle \lambda,\chi\rangle.\]
\end{definition}

\begin{proposition}\label{proposition: existence and uniqueness of balanced filtration for states}
Let $\bXi=(M,\Xi,\alpha)$ be a semistable normed polarised state. There is a unique element $\lambda\in \qfilt(\bXi)$ such that 
\begin{enumerate}
\item $\langle \lambda,\fl\rangle \geq 1$, and 
\item for all $\gamma\in \qfilt(\bXi)$ such that $\langle \gamma,\fl\rangle\geq 1$, we have $\norm{\lambda}\leq\norm{\gamma}$.
\end{enumerate}
\end{proposition}
\begin{proof}
We may assume $\bXi=\bXi'$. If $\Xi=\varnothing$, then $\lambda=0$ is the unique rational filtration satisfying the conditions. Otherwise, $P=\{\gamma\in M^\vee_\R\st \langle \gamma,\chi\rangle \geq 1,\ \forall\chi\in \Xi\}$
is a nonempty closed convex set (actually, an intersection of translated half-spaces) inside $M^\vee_\R=M^\vee_\Q\otimes_\Q \R$, so there is a unique element $\lambda\in P$ minimising the norm $\norm{-}$. To see that $\lambda\in M^\vee_\Q$, note that $\lambda$ lies in the relative interior of a face $F$ of $P$. The affine space $F$ generates is of the form $V_\R+v$, where $V$ is a vector subspace of $M^\vee_\Q$, $V_\R=V\otimes_\Q \R$ and $v\in M^\vee_\Q$. Since $\lambda$ is also the closest point to the origin in $V_\R+v$, we must have that $\lambda=v-p(v)$, where $p$ is the orthogonal projection $p\colon M^\vee_\R\to V_\R$. Since the inner product on $M^\vee_\Q$ is rational, $p$ is defined over $\Q$, and thus $\lambda$ is rational.
\end{proof}

\begin{definition}[Balanced filtration of a state]\label{definition: balanced filtration of a state}
Let $\bXi$ be a semistable normed polarised state. The \emph{balanced filtration} $\lambda_{\text{b}}(\bXi)$ of $\bXi$ is the unique $\lambda\in \qfilt(\bXi)$ satisfying the conditions of \Cref{proposition: existence and uniqueness of balanced filtration for states}.
\end{definition}

\begin{remark}\label{remark: state polystable if and only if balanced filtration is 0}
For a normed semistable polarised state $\bXi=(M,\Xi,\alpha)$, we have $\Xi'=\varnothing$ if and only if $\bXi$ is polystable, if and only if the balanced filtration $\lambda_{\text{b}}(\bXi)=0$.
\end{remark}

\begin{definition}
Let $\bXi$ be a normed semistable polarised state, and let $\lambda\in\qfilt(\bXi)$ be such that $\langle \lambda,\fl\rangle\geq 1$. We define a normed polarised state $\Lambda_\lambda(\bXi)$ as follows. Let $\bXi'=(M',\Xi',0)$ be the slice of $\bXi$, let 
\[(\Xi')_{\lambda,1}=\{\chi\in M'\st \langle \lambda,\chi\rangle=1\}\subset M',\]
and let $\lambda^\vee$ be the unique element of $(M')_\Q$ satisfying $\langle \gamma,\lambda^\vee \rangle=(\gamma,\lambda)$ for all $\gamma\in (M')_\Q^\vee$, where $(-,-)$ denotes the inner product on $(M')_\Q^\vee$. Finally, we set
\[\Lambda_\lambda(\bXi)=(M',(\Xi')_{\lambda,1},\lambda^\vee).\]
\end{definition}

The following theorem is an analogue for states of \Cref{theorem: linear recognition} in the case of algebraic stacks and \cite[Theorem 4.9]{_Haiden_Semistabilitymodularlatticesanditeratedlogarithms} in the case of artinian lattices.

\begin{theorem}[Recognition of the balanced filtration for states]\label{theorem: Recognition of the balanced filtration for states}
Let $\bXi$ be a normed semistable polarised state. Then the balanced filtration of $\bXi$ is the unique filtration $\lambda\in\qfilt(\bXi)$ satisfying
\begin{enumerate}
\item $\langle \lambda,\fl\rangle\geq 1$ and
\item the polarised state $\Lambda_\lambda(\bXi)$ is semistable.
\end{enumerate}
\end{theorem}
\begin{proof}
We use the notation of the proof of \Cref{proposition: existence and uniqueness of balanced filtration for states}. By definition, the balanced filtration $\lambda$ is the unique element in $P$ minimising the function $1/2\norm{-}^2$, whose differential at a point $\gamma\in (M')_\R^\vee$ is precisely $\gamma^\vee\in (M')_\R$. By the Karush-Kuhn-Tucker conditions, for $\gamma\in P$ we have $\gamma=\lambda$ if and only if there are numbers $u_\chi\geq 0$ for each $\chi\in \Xi'$ such that $\gamma^\vee=\sum_{\chi\in \Xi'} u_\chi \chi$ and $u_w=0$ if $\langle \gamma,\chi\rangle >1$. But this precisely means that $\gamma^\vee\in \cone\left((\Xi')_{\gamma,1}\right)$, i.e. that $\Lambda_{\gamma}(\bXi)$ is semistable.
\end{proof}

If $\lambda$ is the balanced filtration of $\bXi$, then $\id_{M'}$ defines a morphism $\Lambda_\lambda(\bXi)\to \Grad_\lambda(\bXi')$ (note that $(\Xi')_{\lambda,0}=\varnothing$), and the quotient map $M\to M'$ defines a morphism $\Grad_\lambda(\bXi')\to \Grad_\lambda(\bXi)$. Therefore $\Lambda_\lambda(\bXi)$ is equipped with a canonical map $\Lambda_\lambda(\bXi)\to \Grad_\lambda(\bXi)$.

\begin{definition}[Balancing chain of a state and the iterated balanced filtration]\label{definition: balancing chain of a state and the iterated balanced filtration}
Let $\bXi$ be a normed semistable polarised state. The \emph{balancing chain} of $\bXi$ is the chain $(\bXi_n,\lambda_n,u_n)_{n\in\N}$ of normed semistable polarised states defined inductively as follows:
\begin{enumerate}
\item $\bXi_0\coloneqq \bXi'$;
\item for every $n\in \N$, $\lambda_n$ is the balanced filtration of $\bXi_n$, $\bXi_{n+1}\coloneqq \left(\Lambda_{\lambda_n}(\bXi_n)\right)'$, and $u_n\colon \bXi_{n+1}\to \Grad_{\lambda_n}(\bXi_n)$ is the composition of $\left(\Lambda_{\lambda_n}(\bXi_n)\right)'\to \Lambda_{\lambda_n}(\bXi_n)$ and the canonical map $\Lambda_{\lambda_n}(\bXi_n)\to \Grad_{\lambda_n}(\bXi_n)$.
\end{enumerate}
The \emph{iterated balanced filtration} of $\bXi$ is the sequential filtration $\lambda_{\text{ib}}(\bXi)\in \qinfilt(\bXi)$ associated to the balancing chain of $\bXi$.
\end{definition}

\begin{proof}[Proof that the balancing chain is well-defined]
We need to check the boundedness condition in \Cref{definition: chain of states}. We observe that, denoting $\bXi=(M,\Xi,\alpha)$, if $\alpha\neq 0$ and $\bXi$ is not polystable, then $\#(\Xi')<\#(\Xi)$. Therefore, eventually $\lambda_n=0$, and hence $\Xi_n'=\varnothing$ by \Cref{remark: state polystable if and only if balanced filtration is 0}, from where the chain stabilises.
\end{proof}

\begin{remark}[Torsion]
The torsion of $M$ does not affect the iterated balanced filtration, but we allow for any finite-type abelian group $M$ in the definition of polarised state in order to make the correspondence with stacks cleaner. Even if we are only interested in the case of a torus action, other diagonalisable groups that are not tori will show up as stabilisers of points, so we need to consider that case too.
\end{remark}

\subsection{From states to good moduli stacks}

We now define a functor from states to pointed stacks, and prove that both theories of iterated balanced filtrations coincide. We fix a field $k$ for the rest of this section.

\begin{definition}[Stack associated to a state]\label{definition: stack associated to a state}
Let $\bXi=(M,\Xi,\alpha)$ be a (normed) polarised state. We define a $k$-pointed (normed) good moduli stack $(\cX_\bXi,x_\bXi)$ over $k$ and a line bundle $\cL_\Xi$ on $\cX_\bXi$ as follows.

First, we denote $G_\bXi=D(M)$ the diagonalisable algebraic group over $k$ Cartier dual to $M$. The group of characters $\Gamma_\Z(G_\bXi)$ of $G_\bXi$ is identified with $M$. 
Let $\A^{\Xi}_k$ be the product of $\#(\Xi)$ many copies of $\A^1_k$. 
Consider the action of $G_\bXi$ on $\A^{\Xi}_k$ through the characters in $\Xi$, that is, $g\cdot(x_\chi)_{\chi\in \Xi}=(\chi(g)x_\chi)_{\chi\in \Xi}$ for $g$ in $G_\bXi$ and $(x_\chi)_{\chi\in \Xi}$ in $\A^\Xi_k$. Let $\cX_\bXi\coloneqq\A^\Xi_k/G_\bXi$ and $x_\bXi=(1,\cdots,1)\in \A_k^\Xi(k)$. We also denote $x_\bXi$ the composition $\Spec k\xrightarrow{x_\bXi} \A^\Xi_k\to \cX_\bXi$. We identify $M_\Q=\Pic(BG_\bXi)\otimes_\Z \Q$, and thus the polarisation $\alpha$ defines a rational line bundle $\cL_\bXi\coloneqq\cO_{\cX_\bXi}(\alpha)=(\cX_\bXi\to BG_\bXi)^*\alpha$ on $\cX_\bXi$. We denote $\ell_\bXi$ the linear form on $\cX_\bXi$ associated to $\cL_\bXi$. If $\bXi$ is normed, the data of the inner product on $M_\Q^\vee$ is equivalent to that of a norm on cocharacters of $G_\bXi$. It thus defines a norm on graded points of $BG_\bXi$ and also a norm on $\cX_\bXi$ by pullback along $\cX_\bXi\to BG_\bXi$.
\end{definition}

We fix a polarised state $\bXi=(M,\Xi,\alpha)$, and denote $\bXi^\circ=(M,\Xi,0)$ the associated “unpolarised” state. We abbreviate $(\cX,x)=(\cX_\bXi,x_\bXi)$ and $G=G_\bXi$.

\begin{proposition}\label{proposition: characterisation of set of filtrations in terms of states}
There is a canonical bijection
\[\qfilt(\cX,x)=\qfilt(\bXi^\circ).\]
\end{proposition}
\begin{proof}
By \Cref{remark: filtrations of a quotient stack}, the set $\qfilt(\cX,x)$ of rational filtrations of $x$ in $\cX$ is identified with the set of those rational cocharacters $\lambda\in M^\vee_\Q=\Gamma^\Q(G)$ such that $\lim_{t\to 0}\lambda(t)x$ exists in $\A^\Xi_k$. Since $\lambda(t)x=(t^{\langle \lambda,\chi\rangle})_{\chi\in \Xi}$, the limit exists precisely when $\langle \lambda,\chi \rangle\geq 0$ for all $\chi\in \Xi$.
\end{proof}

For $y=(y_\chi)_{\chi\in\Xi}\in \A^\Xi(k)$, we define the $\emph{state}$ of $y$ to be the set $\Xi_y=\{\chi\in \Xi\st y_\chi\neq 0\}$. If we let $\bXi_y=(M,\Xi_y,\alpha)$, then we have a closed immersion $f\colon \cX_{\bXi_y}\to \cX$ given by 
\[f\left((z_\chi)_{\chi\in \Xi_y}\right)_\chi=\begin{cases}z_\chi y_\chi,\ \chi\in \Xi_y \\ 0, \ \text{else,}\end{cases}\]
and $f$ maps $x_{\bXi_y}=(1,\cdots,1)$ to $y$.

\begin{proposition}\label{proposition: state of a limit}
Let $\lambda\in \qfilt(\cX,x)$. Then the state of $y=\lim_{t\to 0}\lambda(t)x$ is $\Xi_y=\Xi\cap \partial H_\lambda$.
\end{proposition}
\begin{proof}
Again, $\lambda(t)x=(t^{\langle \lambda,\chi\rangle})_{\chi\in \Xi}$ and thus 
\[y_\chi=\begin{cases}0,\quad \langle \lambda,\chi \rangle>0;  \\ 1,\quad \langle \lambda,\chi\rangle =0;\end{cases}\]
which implies the claim.
\end{proof}

The state also determines the stabiliser of $x$.

\begin{proposition}\label{proposition: stabiliser of a point in terms of state}
Let $K$ be the subgroup of $M$ generated by the elements of $\Xi$, and let $C=M/K$. Let $S=D(C)$ be the Cartier dual of $C$, which is equipped with an injection $S\to G$. Then $S$ is the stabiliser group of $x$.
\end{proposition}
\begin{proof}
The group $G$ acts on $\A^\Xi$ via the characters $\chi\colon G\to \G_{m,k}$, $\chi\in \Xi$, each of which can be seen as the Cartier dual of the map $\Z\to \Gamma_\Z(G)\colon 1\mapsto \chi$. If $\Xi=\{\chi_1,\ldots,\chi_n\}$, then the stabiliser of $x$ is the kernel of $(\chi_1,\ldots,\chi_n)\colon G\to \G_m^n$, and its group of characters is, by Cartier duality, the cokernel of the map $\Z^n\to M\colon e_i\mapsto \chi_i$, which is $C$.
\end{proof}

\begin{proposition}\label{proposition: characterisation of the semistable locus in terms of states}
The point $x$ is semistable in $\cX$ for the linear form $\ell_\bXi$ if and only if the polarised state $\bXi$ is semistable.
\end{proposition}
\begin{proof}
It follows from \Cref{proposition: characterisation of set of filtrations in terms of states}, together with the fact that $\langle \lambda,\cO_\cX(\alpha)\rangle = -\langle \lambda,\alpha\rangle$ with our sign conventions (\Cref{remark: sign convention on weights}), that $x$ is semistable if and only if for all $\lambda\in M^\vee_\Q$ such that $\langle \lambda,\Xi\rangle \geq 0$ (that is, $\Xi\subset H_\lambda$) we have $\langle \lambda,\alpha\rangle \geq 0$ (that is, $\alpha\in H_\lambda)$. The result now follows from this and the fact that $\cone\left(\Xi\right)$ is the intersection of those half-spaces $H_\lambda$ such that $\Xi\subset H_\lambda$ (where we are again abusing notation by identifying $\Xi$ with its image inside $M_\Q$).
\end{proof}

\begin{proposition}\label{proposition: when is the limit semistable in terms of states}
Suppose that $x$ is semistable (equivalently, $\bXi$ is semistable) and let $\lambda\in \qfilt(\cX,x)$. Then the limit $y=\lim_{t\to 0}\lambda(t)x$ is semistable if and only if $\langle \lambda,\alpha\rangle=0$.
\end{proposition}
\begin{proof}
We have equalities $\partial H_\lambda \cap \cone\left(\Xi\right)=\cone \left(\partial H_\lambda \cap \Xi \right)=\cone \left(\Xi_y\right)$, the second of which follows from \Cref{proposition: state of a limit}. The result follows from \Cref{proposition: characterisation of the semistable locus in terms of states} applied to $\cX_{\bXi_y}$, which is a closed substack of $\cX$ containing $y$.
\end{proof}

Note that $\cX$ admits a good moduli space, and $\ell_\bXi$ is trivially $\id_\cX$-positive (\Cref{definition: positive linear form}). Therefore \Cref{theorem: theta stratification proper over gms} implies that the semistable locus $\cX^\ss=\left(\A^\Xi_k\right)^\ss/G$ with respect to $\ell_\bXi$ has a good moduli space $\pi\colon \cX^\ss\to X$.

\begin{proposition}\label{proposition: filtrations of a state and filtrations of the stack of semistables}
There is a canonical bijection $\qfilt(\bXi)\cong \qfilt(\cX^\ss,x)$.
\end{proposition}
\begin{proof}
$\qfilt(\cX^\ss,x)$ is the subset of those $\lambda\in \qfilt(\cX,x)$ such that $\lim_{t\to 0} \lambda(t)x$ is semistable, while $\qfilt(\bXi)$ is the subset of those $\lambda\in \qfilt(\bXi^\circ)$ such that $\langle \lambda,\alpha\rangle=0$. Thus the result follows from \Cref{proposition: characterisation of set of filtrations in terms of states,proposition: when is the limit semistable in terms of states}.
\end{proof}

We can also characterise polystability in terms of the state:
\begin{proposition}\label{proposition: polystability in terms of the state}
The point $x$ is polystable inside $\cX^\ss$ if and only if the state $\bXi$ is polystable.
\end{proposition}
\begin{proof}
The point $x$ being polystable is equivalent to it being semistable and, for all $\lambda\in M^\vee_\Q$ such that $\langle \lambda,\Xi\rangle \geq 0$ and $y=\lim_{t\to 0}\lambda(t)x$ is semistable, having $x=y$. By \Cref{proposition: state of a limit,proposition: characterisation of the semistable locus in terms of states,proposition: when is the limit semistable in terms of states}, this condition is equivalent to having that $\bXi$ is semistable and, for all $\lambda\in M^\vee_\Q$, having that the conditions $\Xi\subset H_\lambda$ and $\langle \lambda,\alpha\rangle =0$ imply that $\Xi\subset \partial H_\lambda$. This means that the smallest face of $\cone\left(\Xi\right)$ containing $x$ is $\cone\left(\Xi\right)$ itself, that is, that $\bXi$ is polystable.
\end{proof}

From now, we do not fix a particular polarised state $\bXi$ and we stop abbreviating $(\cX,x)=(\cX_\bXi,x_\bXi)$ and $G=G_\bXi$.

\begin{definition}\label{definition: map of states defines map between stacks}
Let $\phi\colon \bXi_1=(M_1,\Xi_1,\alpha_1)\to \bXi_2=(M_2,\Xi_2,\alpha_2)$ be a morphism between semistable polarised states. We define pointed morphisms $f_\phi\colon (\cX_{\bXi_1},x_{\bXi_1})\to (\cX_{\bXi_2},x_{\bXi_2})$ and $f_\phi^\ss\colon(\cX_{\bXi_1}^\ss,x_{\bXi_1})\to (\cX_{\bXi_2}^\ss,x_{\bXi_2})$ as follows.

First, the homomorphism $\phi\colon M_2\to M_1$ defines a Cartier dual group homomorphism $D(\phi)\colon G_{\bXi_1}\to G_{\bXi_2}$. The $k$-algebra homomorphism
\[k[t_\chi;\  \chi\in \Xi_2]\to k[t_\psi; \ \psi\in \Xi_1]\colon \chi\mapsto \begin{cases}t_{\phi(\chi)},\ \phi(\chi)\in \Xi_1, \\ 1,\ \phi(\chi)=0.\end{cases}\]
defines, after taking $\Spec$, a $D(\phi)$-equivariant map
\[\A^{\Xi_1}_k\to \A^{\Xi_2}_k\colon (y_\psi)_{\psi\in \Xi_1}\mapsto \left(\begin{cases}y_{\phi(\chi)},\ \phi(\chi)\in \Xi_1, \\ 1,\ \phi(\chi)=0.\end{cases}\right)_{\chi\in\Xi_2}\]
sending $(1,\ldots,1)$ to $(1,\ldots,1)$, and thus a pointed morphism of stacks $f_\phi\colon (\cX_{\bXi_1},x_{\bXi_1})\to (\cX_{\bXi_2},x_{\bXi_2})$. For all geometric points $y=(y_\psi)_{\psi\in \Xi_1}\in \A^{\Xi_1}_k\left(\overline k\right)$, the state of $f_\phi(y)$ contains all elements of $\Xi_2\cap \ker \phi $ and thus $f_\phi(y)$ is semistable. Therefore $f_\phi$ restricts to a morphism $f_\phi^\ss\colon(\cX_{\bXi_1}^\ss,x_{\bXi_1})\to (\cX_{\bXi_2}^\ss,x_{\bXi_2})$.
\end{definition}

The assignments $\phi\mapsto f_\phi$ and $\phi\mapsto f_\phi^\ss$ respect composition. If $\phi$ is a morphism of \emph{normed} semistable polarised states, then $f_\phi$ and $f_\phi^\ss$ are normed morphisms of stacks. Therefore, the assignments $\bXi\mapsto (\cX_\bXi^\ss,x_\bXi)$ and $\phi\mapsto f^\ss_\phi$ define a functor from the category of normed semistable polarised states to the category of $k$-pointed normed good moduli stacks with affine diagonal and finitely presented over $k$.

\begin{proposition}\label{proposition: Grad commutes with passing from states to stacks}
For $\bXi$ a semistable polarised state and $\lambda\in \qfilt(\bXi)$, there is a natural pointed isomorphism $\left(\Grad_\Q(\cX_\bXi^\ss)_{\gr \lambda}, \gr \lambda\right)\cong \left(\cX_{\Grad_\lambda(\bXi)}^\ss, x_{\Grad_\lambda(\bXi)}\right)$. Here, $\Grad_\Q(\cX_\bXi^\ss)_{\gr \lambda}$ is the connected component of $\Grad_\Q(\cX_\bXi^\ss)$ containing $\gr \lambda$.
\end{proposition}
\begin{proof}
By \cite[Theorem 1.4.8]{_HalpernLeistner_Onthestructureofinstabilityinmodulitheory}, $\Grad_\Q(\cX_\bXi^\ss)_{\gr \lambda}=\left(\A^{\Xi}_k\right)^{\lambda,0,\ss}/G_{\bXi}$, where $\left(\A^{\Xi}_k\right)^{\lambda,0}$ is the fixed point locus for the (rational) $\G_{m,k}$ action on $\A^{\Xi}_k$ given by $\lambda$ and $\left(\A^{\Xi}_k\right)^{\lambda,0,\ss}$ is the semistable locus. On the other hand, looking at the weights one gets the equality $\left(\A^{\Xi}_k\right)^{\lambda,0}=\A^{\Xi_{\lambda,0}}_k$, and both $\gr \lambda$ and $x_{\Grad_\lambda(\bXi)}$ are the point $(1,\cdots,1)$, giving the desired isomorphism of pointed stacks.
\end{proof}

\begin{definition}\label{definition: chain of stacks associated to a chain of states}
Let $(\bXi_n,\lambda_n,u_n)_{n\in \N}$ be a chain of normed semistable polarised states. Applying the functor $\bXi\mapsto (\cX_\bXi^\ss,x_\bXi)$ gives an \emph{associated chain of $k$-stacks} $(\cX_{\bXi_n}^\ss,x_{\bXi_n},\lambda_n,h_n)$, where each $h_n$ is the composition of
\[f_{u_n}^\ss\colon \left(\cX^\ss_{\bXi_{n+1}},x_{\bXi_{n+1}}\right)\to \left(\cX_{\Grad_\lambda(\bXi_n)}^\ss, x_{\Grad_\lambda(\bXi_n)}\right),\]
the isomorphism $\left(\cX_{\Grad_\lambda(\bXi_n)}^\ss, x_{\Grad_\lambda(\bXi_n)}\right)\xrightarrow{\sim}\left(\Grad_\Q\left(\cX_{\bXi_n}^\ss\right)_{\gr \lambda}, \gr \lambda\right)$ from \Cref{proposition: Grad commutes with passing from states to stacks}, and the open and closed immersion $\left(\Grad_\Q\left(\cX_{\bXi_n}^\ss\right)_{\gr \lambda}, \gr \lambda\right)\to \left(\Grad_\Q\left(\cX_{\bXi_n}^\ss\right), \gr \lambda\right)$.
\end{definition}

\begin{proposition}\label{proposition: set of sequential filtrations state coincides with stack}
Let $\bXi$ be a semistable polarised state. Then there is a canonical bijection $\qinfilt(\bXi)\cong \qinfilt(\cX_\bXi^\ss,x_\bXi)$.
\end{proposition}
\begin{proof}
The bijection follows from the description of $\qinfilt(\cX_\bXi^\ss,x_\bXi)$ in \Cref{remark: description set of infty filtrations} and an iterated application of \Cref{proposition: Grad commutes with passing from states to stacks,proposition: filtrations of a state and filtrations of the stack of semistables}.
\end{proof}

\begin{proposition}\label{proposition: the stack associated to the slice is the fibre of the good moduli space}
Let $\bXi$ be a semistable polarised state, and let $\bXi'\to \bXi$ be its slice. Then the associated morphism $\cX_{\bXi'}^\ss=\cX_{\bXi'}\to \cX_{\bXi}^\ss$ identifies $\cX_{\bXi'}$ with the fibre of the good moduli space $\cX_{\bXi}^\ss\to X_\bXi^\ss$ containing $x_\bXi$.
\end{proposition}
\begin{proof}

We use the notations of \Cref{definition: slice of a state}. Let $\chi_1,\ldots,\chi_n$ be the different elements of $\Xi$ and assume, after reordering, that $\{\chi_1,\ldots,\chi_l\}=F\cap \Xi$ (where $F$ is the smallest face of $\cone(\Xi)$ containing $\alpha$), the equality to be interpreted inside $M_\Q$ (that is, modulo torsion). 
The state of $\bXi'$ is $\Xi'=q\left(\{\chi_{l+1},\dots,\chi_n\}\right)\subset M'$, where $q\colon M\to M'$ is the quotient map. 

We remark that for any $\lambda\in M^\vee_\Q$ such that $\Xi\subset H_\lambda$ and $\partial H_\lambda\cap \cone\left(\Xi\right)=F$, the limit $y=\lim_{t\to 0}\lambda(t)x$ is polystable. Note also that $y$ does not depend on the choice of $\lambda$. By \Cref{proposition: stabiliser of a point in terms of state} applied to $\Xi_y$, the stabiliser of $y$ is $H\coloneqq G_{\bXi'}$.

We identify $\A^{\Xi}_k=\A^n_k$, the action of $G=G_\bXi$ on $\A^n_k$ being via the characters $\chi_1,\ldots,\chi_n$. The $G$-equivariant open subscheme $\G_{m,k}^l\times \A_k^{n-l}\subset \left(\A^n_k\right)^\ss$ is saturated with respect to the good moduli space $\cX_{\bXi}^\ss=\left(\A^n_k\right)^\ss/G\to \left(\A^n_k\right)^\ss\git G=X^\ss_\bXi$, and thus the fibre of $\cX_\bXi^\ss\to X_\bXi^\ss$ containing $x_\bXi=(1,\ldots,1)$ equals the fibre of the good moduli space $\left(\G_{m,k}^l\times \A_k^{n-l}\right)/G\to \left(\G_{m,k}^l\times \A_k^{n-l}\right)\git G$ containing $(1,\ldots,1)$. Indeed, for every $\overline k $-point $(a,b)$ in $\G_{m,k}^l\times \A_k^{n-l}$, $G_{\overline k}(a,0)$ is the associated polystable orbit, and conversely if a semistable $\overline k$-point $z$ in $\A^n_k$ has associated polystable orbit of the form $G_{\overline k}(a,0)$, then it should lie in $\G_{m,k}^l\times \A_{k}^{n-l}$.

Consider the $(H\to G)$-equivariant map
\[h\colon \A_k^{n-l}\to \G_{m,k}^l\times \A_k^{n-l}\colon (z_1,\ldots,z_{n-l})\mapsto (1,\ldots,1,z_1,\ldots,z_{n-l}),\]
and the associated morphism $\A_k^{n-l}/H\to \left(\G_{m,k}^l\times \A_k^{n-l}\right)/G$, which is the restriction on the codomain of the map $f^\ss_q\colon \cX_{\bXi'}^\ss\to \cX_{\bXi}^\ss$. Let $I$ be the image of the homomorphism $G\to \G_{m,k}^l$ given by $\chi_1,\ldots,\chi_k$.
From a explicit computation of the ring of invariants it follows that $\left(\G_{m,k}^l\times \A_k^{n-l}\right)/G\to\left(\G_{m,k}^l\times \A_k^{n-l}\right)\git G=\G_{m,k}^l/I$ is the good moduli space, and its fibre over $e\in \G_{m,k}^l/I(k)$ is $\left(I\times \A^{n-l}_k\right)/G$.
To conclude, we note the isomorphism
\[\left(I\times \A_k^{n-l}\right)/G\cong\left((G/H)\times \A_k^{n-l}\right)/G\cong \A^{n-l}_k/H\]
induced by $h$.
\end{proof}

From \Cref{proposition: the stack associated to the slice is the fibre of the good moduli space} we get bijections
\[\qfilt(\bXi')\cong \qfilt(\cX_{\bXi'},x_{\bXi'})\cong \qfilt(\cX_{\bXi},x_\bXi)\cong\qfilt(\bXi)\]
between sets of filtrations. This is consistent with \Cref{proposition: slicing does not change set of filtrations of state}.

\begin{proposition}\label{proposition: Kempf number equals complementedness states}
Let $\bXi$ be a semistable polarised state and let $\lambda\in \qfilt(\cX_\bXi,x_\bXi)$. Then the Kempf number $\langle \lambda,\cX_{\bXi'}^{\max}\rangle$ equals the complementedness of $\lambda$:
\[\langle \lambda,\cX_{\bXi'}^{\max}\rangle=\langle \lambda,\fl\rangle.\]
\end{proposition}
\begin{proof}
Write $\cF=\cX_{\Xi'}=\A^m_k/H$ and $x=(1,\ldots,1)$, using notation from the proof of \Cref{proposition: the stack associated to the slice is the fibre of the good moduli space} and where $m=n-l$. The maximal stabiliser locus is $\cF^{\max}=\left(\A^m_k\right)^{H_\circ}/H=\{0\}/H$, where $H_\circ=\left(H^\circ\right)_{\red}$ is the reduced identity component \cite[Proposition C.5]{_Edidin_CanonicalreductionofstabilizersforArtinstackswithgoodmodulispaces}. Let $r=\langle \lambda,\cF^{\max}\rangle$. There is a cartesian square
\[\begin{tikzcd}
    {\Spec \left(k[t]/(t^r)\right)} & {\A^1_k} & t \\
    {\{0\}} & {\A^m_k} & {\left(t^{\langle \lambda,\chi\rangle}\right)_{\chi\in \Xi'}.}
    \arrow["{\lambda(t)x}"', from=1-2, to=2-2]
    \arrow[from=1-1, to=1-2]
    \arrow[from=2-1, to=2-2]
    \arrow[from=1-1, to=2-1]
    \arrow["\lrcorner"{anchor=center, pos=0.125}, draw=none, from=1-1, to=2-2]
    \arrow[maps to, from=1-3, to=2-3]
\end{tikzcd}\]
Taking global sections, we get a cocartesian square
\[\begin{tikzcd}
    {k[t]/(t^r)} & {k[t]} & {t^{\langle \lambda,\chi\rangle}} \\
    k & {k[t_\chi, \chi\in \Xi']} & {t_\chi.}
    \arrow[from=2-2, to=1-2]
    \arrow[from=2-2, to=2-1]
    \arrow[from=2-1, to=1-1]
    \arrow[from=1-2, to=1-1]
    \arrow["\ulcorner"{anchor=center, pos=0.125}, draw=none, from=1-1, to=2-2]
    \arrow[maps to, from=2-3, to=1-3]
\end{tikzcd}\]
Therefore $k[t]/(t^r)=k[t]/(t^{\langle\lambda,\chi\rangle}, \chi\in \Xi')$ and thus $r=\inf \{\langle\lambda,\chi\rangle\st \chi\in \Xi'\}=\langle \lambda,\fl\rangle$.
\end{proof}

\begin{proposition}\label{proposition: balanced filtration in terms of states}
Let $\bXi$ be a normed semistable polarised state. Then the balanced filtration of $(\cX_\bXi^\ss,x_\bXi)$ (\Cref{definition: balanced filtration}) equals, under the bijection $\qfilt(\cX_\bXi^\ss,x_\bXi)\cong\qfilt(\bXi)$, the balanced filtration of $\bXi$ (\Cref{definition: balanced filtration of a state}).
\end{proposition}
\begin{proof}
This follows directly from \Cref{proposition: the stack associated to the slice is the fibre of the good moduli space,proposition: Kempf number equals complementedness states}.
\end{proof}

\begin{theorem}\label{theorem: balancing chain of states corresponds to torsor chain}
Let $\bXi$ be a normed semistable polarised state and let $(\bXi_n,\lambda_n,u_n)_{n\in \N}$ be its balancing chain (\Cref{definition: balancing chain of a state and the iterated balanced filtration}). Then the chain of stacks associated to $(\bXi_n,\lambda_n,u_n)_{n\in \N}$ (\Cref{definition: chain of stacks associated to a chain of states}) is isomorphic to the torsor chain of $(\cX_\bXi^\ss,x_\bXi)$ (\Cref{definition: torsor chain stack}).
\end{theorem}
\begin{proof}
Let $(\cY_n,y_n,\eta_n,v_n)_{n\in \N}$ be the torsor chain of $(\cX^\ss_{\bXi},x_\bXi)$. We will provide, for all $m\in \N$, isomorphisms
\[i_m \colon (\cX^\ss_{\bXi_m},x_{\bXi_m})\cong (\cY_m,y_m)\]
such that,
\begin{enumerate}
\item\label{item 1 theorem states and iterated balanced filtration} under the identification $\qfilt(\bXi_m)=\qfilt(\cY_m, y_m)$, we have $\eta_m=\lambda_m$; and
\item\label{item 2 theorem states and iterated balanced filtration} the square
\begin{equation}\label{diagram: square for having isomorphism of chains states and torsor}
\begin{tikzcd}
    {(\cX^\ss_{\bXi_{m+1}},x_{\bXi_{m+1}})} & {(\cX_{\Grad_{\lambda_m}(\bXi_m)}^\ss,x_{\Grad_{\lambda_m}(\bXi_m)})} \\
    {(\cY_{m+1},y_{m+1})} & {(\Grad_\Q(\cY_m)_{\gr \eta_m}, \gr \eta_m)}
    \arrow["{f^\ss_{u_m}}", from=1-1, to=1-2]
    \arrow["{v_n}", from=2-1, to=2-2]
    \arrow[from=1-2, to=2-2]
    \arrow["{i_{m+1}}"', from=1-1, to=2-1]
\end{tikzcd}\end{equation}
commutes, where the arrow on the right comes from \Cref{proposition: Grad commutes with passing from states to stacks}.
\end{enumerate}
For $n=0$, we have that $y_0$ is the fibre of the good moduli space of $\cX_\bXi^\ss$ containing $x_\bXi$, and that $\bXi_0=\bXi'$. Therefore $(\cY_0,y_0)=(\cX_{\bXi_0}^\ss,x_{\bXi_0})$, by \Cref{proposition: the stack associated to the slice is the fibre of the good moduli space}. Let $n\in \N$ and suppose that isomorphisms $i_m$ as above have been provided in such a way that conditions (\ref{item 1 theorem states and iterated balanced filtration}) and (\ref{item 2 theorem states and iterated balanced filtration}) above hold for all $m<n$. Since $i_n$ is an isomorphism, we have the equality $\eta_n=\lambda_n$ by \Cref{proposition: balanced filtration in terms of states}, so condition (\ref{item 1 theorem states and iterated balanced filtration}) also holds for $m=n$.

If $\lambda_n=0$, then $\bXi_n$ and $(\cY_n,y_n)$ are polystable. Therefore $\bXi_{n+1}=\bXi_n$ and $(\cY_{n+1},y_{n+1})=(\cY_n,y_n)$, so there is nothing to prove.

Assume $\lambda_n\neq 0$. We freely use the notation of Case 2 in \Cref{construction: torsor chain}. Let $\chi_1,\ldots,\chi_l$ be the different elements in $\bXi_n$, and let $G=G_{\bXi_n}$. We denote $V=k^l$ the $G$-representation given by the characters $\chi_1,\ldots,\chi_l$, so that $\A^{\Xi_0}_k=\A(V)$ (where $\A$ denotes \emph{total space}). We have $\cY_n=\A(V)/G$, and $\cY_n^{\max}=\{0\}/G$. Therefore the relevant blow-up is $\cB=\left(\Bl_0 \A(V)\right)/G$, the exceptional divisor is $\cE=\P(V)/G$, and the $\G_m$-torsor over it is $\cN=\left(\A(V)\setminus \{0\}\right)/G$. We denote $y_n$ also the unique lift of $y_n=(1,\ldots,1)$ to $\Bl_0 \A(V)$. The rational one-parameter subgroup $\lambda_n$ has the property that $\langle \lambda_n,\chi_j\rangle \geq 1$ for all $j$ and that equality holds for at least one $j$. Therefore, if we set $z=\lim_{t\to 0}\lambda_n(t)y_n$, the limit taken inside $\Bl_0 \A(V)$, and if we write $z=[z_1,\ldots,z_l]$ in projective coordinates, noting that $z$ lies on the exceptional divisor $\P(V)$, then we have $z_j=0$ if $\langle \lambda_n,\chi_j\rangle >1$ and $z_j=1$ if $\langle \lambda_n,\chi_j\rangle =1$. We denote $z^*=(z_1,\ldots,z_l)$ this lift of $z$ to $\A(V)$.

The limit $z$ lifts to the connected component $\overline \cZ$ of $\Grad(\cB)$ containing $\gr \eta_n$. If we set $V_1=\bigoplus_{\langle \lambda,\chi\rangle=1} V_\chi$, where $V_\chi$ is the subrepresentation of $V$ where $G$ acts via the character $\chi$, then $\overline \cZ=\P(V_1)/G$. Thus $\gr \eta_n$ is identified with $z\in \P(V_1)(k)$, and the lift $z^*\in \A(V_1)$ to $\A(V_1)$ can be written in coordinates as $z=(1,\ldots,1)$ when seen inside $V_1$. The centre $\cZ$ of the locally closed $\Theta$-stratum of $\cB$ containing $y_n$ is the semistable locus for the shifted linear form $\ell_c$, by the Linear Recognition Theorem \ref{theorem: linear recognition}. 

In this case, $c=\dfrac{1}{\norm{\lambda_n}^2}$, and the shifted linear form is
\[\ell_c=\ell-\dfrac{1}{\norm{\lambda_n}^2}\langle \lambda_n^\vee,-\rangle,\]
where $\ell$ is the linear form associated to the ample line bundle $\cO_{\P(V_1)/G}(1)$. Let $x\in \P(V_1)(\overline k)$ and let $x^*\in \A(V_1)(\overline k)$ be a lift of $x$ to $\A(V_1)$ with state $\Xi_x$. The point $x$ is semistable for $\ell_c$ if for all $\gamma\in \Gamma^\Q(G)$ we have
\[0\geq \ell_c(\gamma)=\min \langle \gamma, \Xi_x\rangle - \dfrac{1}{\norm{\lambda_x}^2}(\gamma,\lambda_n),\]
which holds if and only if $0\in \conv\left(\Xi_x-\dfrac{\lambda_n^\vee}{\norm{\lambda_n}^2}\right)$. Since all elements of $\Xi_x$ are in the hyperplane $\langle\lambda_n,-\rangle =1$, the condition that $0$ is in the convex hull of $\Xi_x-\dfrac{\lambda_n^\vee}{\norm{\lambda_n}^2}$ is equivalent to the condition that $\lambda_n^\vee$ is in the cone generated by $\Xi_x$, by \Cref{lemma: convex geometry lemma comparison cone and conv} below. This is in turn equivalent to the lift $x^*$ of $x$ to $\A(V_1)$ being semistable for the linear form given by $\lambda_n^\vee$, by \Cref{proposition: characterisation of the semistable locus in terms of states}. Therefore we have a cartesian square
\[\begin{tikzcd}
    {\cM=\A(V_1)^{\ss(\lambda_n^\vee)}/G} & {\cN=\left(\A(V)\setminus 0\right)/G} \\
    {\cZ=\P(V_1)^{\ss(\ell_c)}/G} & {\cE=\P(V)/G,}
    \arrow[from=1-2, to=2-2]
    \arrow[""{name=0, anchor=center, inner sep=0}, from=2-1, to=2-2]
    \arrow[from=1-1, to=2-1]
    \arrow[from=1-1, to=1-2]
    \arrow["\ulcorner"{anchor=center, pos=0.125}, draw=none, from=1-1, to=0]
\end{tikzcd}\]
using the notations of \Cref{construction: torsor chain}. Just from the definitions, we see that $\cM=\cX^\ss_{\Lambda_{\lambda_n}(\bXi_n)}$. We choose $y_{n+1}=z^*$ as the preimage of $z$ along $\cM\to \cZ$ needed for the construction of the torsor chain. The stack $\cY_{n+1}$ is by definition the fibre of the good moduli space of $\cM$, and hence by \Cref{proposition: the stack associated to the slice is the fibre of the good moduli space} we have the desired isomorphism $i_{n+1}\colon \cX^\ss_{\bXi_{n+1}} = \cX^\ss_{\Lambda_{\lambda_n}(\bXi_n)'}\cong \cY_{n+1}$, which sends $x_{\bXi_{n+1}}$ to $y_{n+1}$ by our choice of $z^*$. Both $\cX_{\Grad_{\lambda_n}(\bXi_n)}^\ss$ and $\Grad_\Q(\cY_n)_{\gr \eta_n}$ are naturally identified with $BG$, so the square (\ref{diagram: square for having isomorphism of chains states and torsor}) commutes for $m=n$. Since the torsor chain is bounded, repeating this process we eventually reach the case $\eta_n=0$, getting the desired isomorphism of chains.
\end{proof}

In the proof of \Cref{theorem: balancing chain of states corresponds to torsor chain} we used the following fact in convex geometry.

\begin{lemma}\label{lemma: convex geometry lemma comparison cone and conv}
Let $N$ be a finite dimensional $\Q$-vector space endowed with a rational inner product $(-,-)$. Let $\Xi\subset N$ be a nonempty finite set and let $\gamma\in N$ be an element such that $(\gamma,\chi)=1$ for all $\chi\in \Xi$. Then we have $0\in \conv\left(\Xi - \dfrac{\gamma}{\norm{\gamma}^2}\right)$ if and only if $\gamma \in \cone(\Xi)$.
\end{lemma}
\begin{proof}
Each $\chi\in \Xi$ can be written as $\chi=\dfrac{\gamma}{\norm{\gamma}^2}+\beta_\chi$ with $(\gamma,\beta_\chi)=0$. Note that the condition $0\in \conv\left( \Xi - \dfrac{\gamma}{\norm{\gamma}^2}\right)$ is equivalent to $0\in \cone\left( \Xi - \dfrac{\gamma}{\norm{\gamma}^2}\right)$. If this is satisfied, then $0=\sum_\chi c_\chi\beta_\chi$ with $c_\chi\geq 0$. After rescaling we may assume that $\sum c_\chi=\norm{\gamma}^2$. Then
\[\gamma=\left(\sum c_\chi\right)\dfrac{\gamma}{\norm{\gamma}^2}+\sum_\chi c_\chi\beta_\chi=\sum c_\chi \chi,\]
and thus $\gamma\in\cone (\Xi)$.

Conversely, if $\gamma=\sum_\chi c_\chi \chi$ with the $c_\chi\geq 0$, then after applying $(\gamma,-)$ to both sides we get $\norm{\gamma}^2=\sum_\chi c_\chi$, so 
\[\gamma=\sum_\chi c_\chi\left(\dfrac{\gamma}{\norm{\gamma}}+\beta_\chi\right)=\gamma +\sum_\chi c_\chi\beta_\chi\]
and thus  $0=\sum c_\chi \beta_\chi$ is in $\cone\left( \Xi - \dfrac{\gamma}{\norm{\gamma}^2}\right)$.
\end{proof}

\begin{corollary}\label{corollary: iterated balanced filtration for states equals that for stacks}
Let $\bXi$ be a normed semistable polarised state. Then the iterated balanced filtration (\Cref{definition: balancing chain of a state and the iterated balanced filtration}) of $\bXi$ equals, under the bijection $\qinfilt(\bXi)\cong \qinfilt(\cX_\bXi^\ss,x_\bXi)$, the iterated balanced filtration of $(\cX_\bXi^\ss,x_\bXi)$ (\Cref{definition: iterated balanced filtration}).
\end{corollary}
\begin{proof}
By \Cref{theorem: balancing chain of states corresponds to torsor chain}, the iterated balanced filtration of $\bXi$ equals the sequential filtration associated to the torsor chain of $(\cX_\bXi^\ss,x_\bXi)$. The results then follows from \Cref{theorem: torsor computes the iterated balanced filtration}.
\end{proof}

\begin{example}\label{example: simple case of conjecture on flows}
Consider the normed polarised state 
\[\bXi=\left(\Z^2,\left\{ (1,0),(1,1)\right\},0\right),\] 
where $(\Z^2)^\vee=\Z^2$ has the standard inner product ($(1,0),(0,1)$ is an orthonormal base). The associated stack over $\C$ is $\cX_{\bXi}=\C^2/(\C^\times)^2$, where $(\C^\times)^2$ acts by $(t_1,t_2)(v_1,v_2)=t_1v_1,t_1t_2v_2)$. The linearisation $\cL_{\bXi}$ is trivial, so every point is semistable.

Let us compute the iterated balanced filtration of the state $\bXi$. We have that $\bXi$ is its own slice $\bXi'=\bXi$, and the balanced filtration $\lambda_0=(a,b)$ of $\bXi$ is the minimiser of $a^2+b^2$ subject to the condition that $\langle \lambda_0, (1,0)\rangle =a\geq 1$ and $\langle \lambda_0, (1,1)\rangle =a+b\geq 1$. Therefore $\lambda_0=(1,0)$. The iterated state is 
\[\Lambda_{\lambda_0}(\bXi)= \left(\Z^2,\left\{(1,0),(1,1)\right\},(1,0)\right),\]
whose slice is $\bXi_1=\Lambda_{\lambda_0}(\bXi)'=\left(\Z(0,1), \left\{(0,1)\right\},0\right)$, and the balance filtration of $\bXi_1$ is $\lambda_1=(0,1)$. Since $\Lambda_{\lambda_1}(\bXi_1)$ is polystable, the balancing chain of $\bXi$ terminates here, and we conclude that the iterated balanced filtration of $\bXi$ is $\lambda_0=(1,0),\lambda_1=(0,1)$. By \Cref{corollary: iterated balanced filtration for states equals that for stacks}, we deduce that the iterated balanced filtration of $x_\bXi=(1,1)\in \C^2/(\C^\times)^2$ is the sequence $(1,0),(0,1)$ in $\Q^2=\Gamma^\Q((\C^\times)^2)$.

We now analyse \Cref{conjecture: on asymptotics of gradient flows affine GIT} in this case. We endow $\C^2$ with the standard hermitian metric. The associated Kempf-Ness potential for the point $(1,1)\in \C^2$ is 
\[f\colon \R^2\to \R\colon x\mapsto e^{x\cdot \binom{1}{0}}+e^{x\cdot \binom{1}{1}}\]
up to the addition of a constant, where we are identifying $\R^2\xrightarrow{i\cdot-}\Lie((S^1)^2)\xrightarrow[\sim]{exp(\frac{1}{2i}-)}(\C^\times)^2/(S^1)^2$, the unit circle $S^1$ being the maximal compact subgroup of $\C^\times$. In this case, the exponential map is a global isometry between $\Lie((S^1)^2)$ and $(\C^\times)^2/(S^1)^2$. The gradient is
\[\nabla f(x)=e^{x_1}\binom{1}{0}+e^{x_1+x_2}\binom{1}{1}.\]
The equation for $h\colon (0,\infty)\to \R^2$ to be a flow line for $-\nabla f$ is
\[h'(t)=-e^{h_1(t)}\binom{1}{0}-e^{h_1(t)+h_2(t)}\binom{1}{1}.\]
We write $h(t)=-\log(t)\binom{1}{0}-\log\log(t)\binom{1}{0}+z(u)$, where $u=\log\log(t)$. The equation for $z$ becomes
\[z'(u)=\binom{e^u(1-e^{z_1(u)})-e^{z_1(u)+z_2(u)}}{1-e^{z_1(u)+z_2(u)}}.\]
For $N_1>0$ real, we have that
\begin{enumerate}
\item if $z_2=N_1$ and $z_1\in (-N_1,N_1)$, then $z_2'<0$.
\item if $z_2=-N_1$ and $z_1\in (-N_1,N_1)$, then $z_2'>0$.
\end{enumerate}
For $N_2>0$ big enough so that $-1\leq \log(1-e^{-N_2})$ and $u\geq 2$, we have that
\begin{enumerate}
\item if $z_1=-N_2$ and $z_2\in (-N_2-1,N_2+1)$, then $z_1'>0$.
\item if $z_1=N_2$ and $z_2\in (-N_2-1,N_2+1)$, then $z_1'<0$.
\end{enumerate}
Therefore, an appropriate choice of $N_1$ and $N_2$ gives a rectangle that $z$ cannot leave, because $z'$ points inwards at the boundary. Therefore $z$ is bounded when $t>>0$. We have verified \Cref{conjecture: on asymptotics of gradient flows affine GIT} in this example.
\end{example}

\footnotesize{\bibliography{bibliographyall}{}
\bibliographystyle{abbrv}}

 \textsc{Mathematical Institute, University of Oxford,
Andrew Wiles Building,
Radcliffe Observatory Quarter (550),
Woodstock Road,
Oxford,
OX2 6GG, United Kingdom}\par\nopagebreak
  \textit{E-mail address}, \texttt{ibaneznunez@maths.ox.ac.uk}
\end{document}